%% file: main.tex
\documentclass[a4paper,11pt]{article}
\input{preamble}

\input{newcommands}
\input{abbreviation}

\graphicspath{{../04_siam/figures/}}

\begin{document}


\title{\textbf{A Hybrid semi-Lagrangian Flow Mapping Approach for Vlasov Systems}\\ 
{\Large Combining Iterative and Compositional Flow Maps}}
 \date{\today}

\author{
     Philipp Krah\thanks{First authors: \href{mailto:philipp.krah@univ-amu.fr}{\nolinkurl{philipp.krah@univ-amu.fr}}, \href{mailto:zetao.lin@etu.univ-amu.fr}{\nolinkurl{zeato.lin@edu.univ-amu.fr}}} \thanks{Aix-Marseille Université, I2M, CNRS, UMR 7373, 3 place Victor Hugo, 13331 Marseille cedex 3, France.} \thanks{IRFM - CEA Cadarache,
13108 Saint-Paul-lez-Durance,
France }
    \and Zetao Lin\footnotemark[1] \footnotemark[2]
    \and R.-Paul Wilhelm\thanks{Centre for mathematical Plasma Astrophysics, Department of Mathematics, KU Leuven, B-3001 Leuven, Belgium} 
    \and Fabio Bacchini\footnotemark[4] \thanks{Royal Belgian Institute for Space Aeronomy, Solar-Terrestrial Centre of Excellence, B-1180 Uccle, Belgium}
    \and Jean-Christophe Nave\thanks{McGill University, Department of Mathematics and Statistics, Montreal, Canada}
    \and Virginie Grandgirard\footnotemark[3]
    \and Kai Schneider\footnotemark[2]
}

\maketitle


\begin{abstract}
We propose a hybrid semi-Lagrangian scheme for the Vlasov--Poisson equation that combines the Numerical Flow Iteration (NuFi) method with the Characteristic Mapping Method (CMM). Both approaches exploit the semi-group property of the underlying diffeomorphic flow, enabling the reconstruction of solutions through flow maps that trace characteristics back to their initial positions.
NuFi builds this flow map iteratively, preserving symplectic structure and conserving invariants, but its computational cost scales quadratically with time. Its advantage lies in a compact, low-dimensional representation depending only on the electric field. In contrast, CMM achieves low computational costs when remapping by composing the global flow map from explicitly stored submaps. 
The proposed hybrid method merges these strengths: NuFi is employed for accurate and conservative local time stepping, while CMM efficiently propagates the solution through submap composition. This approach reduces storage requirements, maintains accuracy, and improves structural properties. Numerical experiments demonstrate the effectiveness of the scheme and highlight the trade-offs between memory usage and computational cost.  
We benchmark against a semi-Lagrangian predictor--corrector scheme used in modern gyrokinetic codes, evaluating accuracy, conservation properties.
\end{abstract}
 
    
\input{content}

\section*{Author Contribution Statement}

The following outlines the authors' contributions to this work. Z. Lin and P. Krah contributed equally as first authors. 

{\small
\noindent
\begin{tabular}{@{}p{0.17\linewidth} p{0.77\linewidth}}
\textbf{P. Krah:} & Conceptualization, Methodology, Investigation, Implementation, Software, Validation, Visualization, Writing original draft,\\
\textbf{Z. Lin:} & Implementation, Visualization, Numerical Simulations, Writing – original draft,\\
\textbf{F. Bacchini:} & Review \& editing, Funding acquisition.\\
\textbf{P. Wilhelm:} & Numerical Simulations, Writing original draft, Review and Editing, \\
\textbf{J.-C. Nave} & Review \& Editing, Funding acquisition,\\
\textbf{V. Grandgirard} & Review \& Editing, Supervision of P.Krah, Funding acquisition,\\
\textbf{K. Schneider:} & Writing original draft, Review \& Editing, Funding acquisition, supervision.
\end{tabular}
}

\section*{Code Availability}
The code will be made publicly available upon publication.

\section*{Acknowledgments}
The authors thank Xi Yuan Yin and William Holman-Besinger for fruitful discussions.
This work was supported by the Institute for Fusion Sciences and Instrumentation in Nuclear Environments (ISFIN) at Aix-Marseille Université, funded by the French government under the France 2030 program managed by the A*MIDEX initiative (AMX-19-IET-013). Additional support was provided by the French Federation for Magnetic Fusion Studies (FR-FCM) and the Eurofusion Consortium, funded by the Euratom Research and Training Programme (Grant Agreement No. 633053). The views and opinions expressed herein do not necessarily reflect those of the European Commission. Support was provided by the Agence Nationale de la Recherche (ANR), project CM2E (ANR-20-CE46-0010-01). 
The authors acknowledge the Grand Équipement National de Calcul Intensif (GENCI) for granting access to the HPC resources of IDRIS (Allocation No. AD012A01664R1) and thank the Centre de Calcul Intensif d’Aix-Marseille for its computing resources. Z.L. thanks KU Leuven for its hospitality. The work of J-C.N. was partially supported by the NSERC Discovery Grant program. 

\bibliographystyle{siamplain}
\bibliography{references}


\input{appendix}

\end{document}

%% file: preamble.tex
\usepackage[margin=0.8in]{geometry}
\usepackage{adjustbox}
\usepackage{algorithm}
\usepackage{algorithmic}
\usepackage{amsmath,amssymb,amsfonts,amsthm}
\usepackage{array}
\usepackage{blkarray}
\usepackage{bm}
\usepackage{booktabs}
\usepackage{cite}
\usepackage[numbers]{natbib}  
\usepackage{hyperref}
\hypersetup{
    colorlinks,
    linkcolor={red!50!black},
    citecolor={blue!50!black},
    urlcolor={blue!80!black}
}
\usepackage{cleveref}
\usepackage{diagbox}
\usepackage{dsfont}
\usepackage{enumerate}
\usepackage[inline]{enumitem}
\usepackage{graphicx}
\usepackage{mathtools}
\usepackage{multirow}
\usepackage{nth}
\usepackage{siunitx}
\usepackage{soul}
\usepackage{stmaryrd}
\usepackage{subcaption}
\usepackage{textcomp}
\usepackage{xcolor}
\usepackage{parskip}
\usepackage{url}
\usepackage{tikz}
\usepackage{tikzsymbols}
\usepackage{pgfplots}
\pgfplotsset{compat=newest} 
\usetikzlibrary{arrows,shapes}
\usetikzlibrary{decorations.markings}
\usetikzlibrary{decorations.pathmorphing}
\usetikzlibrary{arrows.meta}
\usetikzlibrary{shadows,arrows,positioning,shapes.geometric}
\tikzset{cross/.style={cross out, draw=black, inner sep=0pt, outer sep=0pt},cross/.default={1pt}}
\usepackage{tikzsymbols}
\usepackage{pgfplots}
\pgfplotsset{compat=newest} 
\usetikzlibrary{plotmarks}
\usetikzlibrary{arrows.meta}
\usepgfplotslibrary{patchplots}
\usepackage{grffile}
\usepackage{pdflscape}
\usepackage{fontawesome}
\pgfplotsset{every axis/.append style={
                    label style={font=\scriptsize},
                    tick label style={font=\scriptsize},
                    legend style={font=\scriptsize}
                    }}
\pgfplotsset{compat=newest}
\pgfplotsset{plot coordinates/math parser=false}
\pgfplotsset{grid style={dotted,gray}}
\pgfplotsset{
compat=1.11,
legend image code/.code={
\draw[mark repeat=2,mark phase=2]
plot coordinates {
(0cm,0cm)
(0.15cm,0cm)        
(0.3cm,0cm)         
};%
}
}
\usepgfplotslibrary{groupplots}
\newlength\figureheight
\newlength\figurewidth
\pgfdeclarelayer{background}
\pgfdeclarelayer{foreground}
\pgfsetlayers{background,main,foreground}

\newtheorem{theorem}{Theorem}[section]

\newtheorem{remark}{Remark}
\newtheorem{example}{Example}

%% file: newcommands.tex
\renewcommand{\vec}[1]{\boldsymbol{#1}}

\renewcommand{\d}{\mathrm{d} }
\newcommand{\flowmap}{\vec{\Phi}}
\newcommand{\submap}{\flowmap}
\newcommand{\scalar}{f}
\newcommand{\abs}[1]{|#1|}
\newcommand{\Ord}[1]{\mathcal{O}(#1)}

\newcommand{\dds}{\ensuremath{\tfrac{{\mathrm d}}{{\mathrm d}s}}}

\newcommand{\R}{\ensuremath{\mathbb{R}}} 

\newcommand{\dv}{\ensuremath{{\mathrm d}v}} 

\newcommand{\interp}{\mathcal{P}}
\newcommand{\norm}[1]{\Vert #1 \Vert}
\newcommand{\Nmaps}{M}
\newcommand{\Nremap}{N_\text{remap}}
\newcommand{\mapcmm}{\vec{\chi}}
\newcommand{\mapgrid}{\mathcal{G}_{\mapcmm}}
\newcommand{\Nmapgrid}{N_{\mapcmm}}
\newcommand{\Nsample}{N_{f}}

\newcommand{\samplegrid}{\mathcal{G}_{f}}

\newcommand{\NuFiMap}{\boldsymbol{\Psi}}
\newcommand{\Nnufi}{N}
\newcommand{\NnuFi}{\Nnufi}

\newcommand{\energy}{\mathcal{E}}
\newcommand{\mass}{\mathcal{M}}
\newcommand{\density}{\rho}
\newcommand{\flux}{j}
\newcommand{\momentum}{\mathcal{P}}
\newcommand{\Etot}{\energy}                 
\newcommand{\Epot}{\energy_\mathrm{pot}}                 
\newcommand{\Ekin}{{\energy}_\mathrm{kin}}               

\newcommand{\highlight}[1]{\textit{#1}}
\newcommand{\define}[1]{\textit{#1}}                            

\usepackage{xcolor}

%% file: abbreviation.tex
\usepackage[acronym,automake,nonumberlist]{glossaries}
\newcommand{\KE}{{KE}\xspace}
\newcommand{\VP}{{VP}\xspace}

\makeglossaries

%% file: content.tex
\section{Introduction} 

This work addresses the efficient numerical resolution of fine-scale structures in collisionless kinetic plasma systems, a long-standing challenge in Vlasov-based simulations. Such structures can emerge from various physical effects: phase mixing of small perturbations, as in Landau damping; resonant interactions between counter-streaming electron populations leading to the two-stream instability; 
and many others.  
Accurately capturing these multiscale features requires numerical representations that remain both memory-efficient and computationally affordable. To this end, we introduce a new approach that exploits the semigroup properties of the characteristic flow. By representing the solution through a composition of flowmaps, rather than through increasingly dense phase space grids, we achieve a more scalable resolution of fine-scale dynamics than with conventional methods.

In this work, we focus on the one-dimensional, one-velocity (1D+1V) Vlasov--Poisson (VP) system,
\begin{equation}
\label{eqn:vlasov}
\partial_t f + v \, \partial_{x} f 
+ E \, \partial_{v} f = 0\,,
\end{equation}
where $f$ denotes the particle distribution function and $E$ is the electric field
\begin{align}
\label{eqn:poisson}
-\Delta_{x} \varphi = \rho = 1 - \int_\R f \, \mathrm{d}v, \\
\label{eqn:Edef}
E = -\nabla_{x}\varphi\,.
\end{align}
The numerical treatment of the Vlasov equation is particularly demanding: The distribution function is in general 
defined over a six-dimensional (in this work, two-dimensional) phase space, which gives rise to the curse of dimensionality and thus forces practical simulations to rely on reduced models or coarse discretizations. In addition, the nonlinear and collisionless character of the VP dynamics promotes the formation of fine-scale filamentary structures in phase space. These filaments are notoriously difficult to resolve accurately, yet they are essential for understanding the onset of kinetic instabilities and may play a significant role in plasma heating and dissipation processes~\cite{case_for_electron_astrophysics,multi_scale_solar_wind,lapenta_mercury}.

\subsubsection*{State of the Art Plasma Simulations}
The numerical solution of the Vlasov equation generally employs methods falling into one of three classes: Eulerian (discrete velocity), Semi-Lagrangian, or Lagrangian, as discussed e.g.\ in reviews by Filbet and Sonnendrücker \cite{review_num_kin_filbet} and Palmroth et al. \cite{Palmroth2018}.
The first two approaches are grid-based, meaning that the solution $f$ is represented on a discretization of the 
phase-space. Eulerian discretizations compute the evolution of the distribution by evaluating its local 
temporal derivative at each grid node, thereby recasting the Vlasov equation as a large coupled system of 
ODEs defined on the computational mesh \cite{ArberVann2002, gene_1, gene_2, JUNO2018110, Mandell_Hakim_Hammett_Francisquez_2020, fijalkow_numerical_1999, filbet_conservative_2001}. 
While Eulerian techniques are widely used in fluid dynamics, they are less common for the  
Vlasov equation. Although the Eulerian framework provides a general procedure for 
constructing numerical solvers for PDEs, it typically does not exploit the additional structural 
properties of the system. A well-known plasma gyrokinetic code that is based on an Eulerian description is 
GENE \cite{jenko_gene_2011,JenkoDorlandKotschenreutherRogers2000}, which uses a mix between pseudospectral 
and finite-difference approximations to solve the Vlasov part. It has been used in studies of existing experimental devices 
\cite{bahner2026magnetic} and in the design of upcoming devices \cite{LionStealris2025}.

Because the phase space velocity field of the Vlasov--Poisson equations is divergence-free, characteristics do not cross. Semi-Lagrangian (SL) methods make use of this: they retain a phase-space grid but 
evolve grid points backward in time for one time step along the characteristic equations, 
reconstructing the solution at the new time step via interpolation at the characteristic foot 
\cite{ChenKnorr1976,SonnendrückerRocheBertrandGhizzo1999}, which is generally referred to as the  
\emph{backward Semi-Lagrangian} (BSL) approach. Alternatively, it is also possible to evolve 
points forward in time, which is called the \emph{forward Semi-Lagrangian} (FSL) method. FSL  
somewhat simplifies the extension to higher order schemes, however, most Semi-Lagrangian codes still implement the more accurate and stable BSL \cite{CrouseillesRespaudSonnendruecker2009}.
The Semi-Lagrangian approach improves conservation properties compared to pure Eulerian schemes, particularly when combined with symplectic time integration. 
Cubic B-splines are a common choice for interpolation due to their balance between accuracy and 
efficiency \cite{Riishøjgaard1998, kormann_massively_parallel}, though some codes use nodal interpolation \cite{sldg2015, Einkemmer2019, CrouseillesMehrenbergerVecil2011, cottet_etancelin_perignon_picard_2014, COTTET2018362, RossmanithSeal2011, PalmrothGansePfau-KempfBattarbeeTurcBritoGrandinHoilijokiSandroosVonAlfthan2018a, muphy_I, muphy_II}.
In this work, we will refer to an SL predictor-corrector scheme using B-splines. It is the basis of Gysela \cite{GrandgirardSarazinAngelinoBottinoCrouseillesDarmetDif-PradalierGarbetJolliet2007} that 
simulates 5D gyrokinetic equations. It is currently rewritten and generalized for exascale computations \cite{bourne2025gyselalib++} using a multi-patch semi-Lagrangian approach \cite{VidalBourneGrandgirardMehrenbergSonnendruecker2025}.

However, there are several drawbacks when using grid-based methods to solve the Vlasov equation. Maintaining a grid of up to six dimensions is extremely costly, especially with adaptivity or non-rectangular domains. 
The high memory requirements not only severely restrict the maximum possible resolution but also lead to significant
communication overhead, which in turn limits parallel efficiency and scalability, especially on modern accelerator hardware such as GPUs \cite{kormann_massively_parallel, Einkemmer2019, sl_456d_einkemmmer_moriggl}. 
Sparse-grid and adaptive refinements become difficult to deploy in high-dimensional phase space \cite{sparse_kormann_sonnendruecker, muphy_II}. Moreover, any grid-based discretization introduces numerical 
diffusion, which reintroduces artificial dissipation into the collisionless Vlasov dynamics and undermines 
conservation properties. Although several correction techniques exist \cite{filbet_conservative_2001, CrouseillesMehrenbergerVecil2011}, the effect cannot be completely removed. Grid methods also struggle with 
irregular geometries, nontrivial boundary conditions, and the need to accommodate growing velocity support.

An alternative to grid-based schemes are particle-based, or Lagrangian, solvers. The most common 
of these are the \emph{Particle-In-Cell} (PIC) methods, though other approaches exist, such as 
\emph{smoothed particle hydrodynamics} (SPH) and the \emph{Reproducing Kernel Hilbert Space 
Particle Method} (RKHS-PM) \cite{cottet_particle_1984, HockneyEastwood2021, BirdsallLangdon2018, WILHELM2023111720}. 
The idea behind Lagrangian methods is to approximate the distribution function with a collection of 
marker particles carrying weights or function values, which are advanced along the phase flow of the Vlasov equation. Since these trajectories coincide with the characteristics of the Vlasov equation, the particle weights remain constant in time—the same principle exploited by Semi-Lagrangian methods. 
In PIC, a grid in physical space is additionally introduced to solve the Poisson equation for the electric field
(or the Maxwell equations for the electro-magnetic case). The charge density, which enters the right-hand side of 
\eqref{eqn:poisson}, is then obtained either by counting particles per cell or through higher-order deposition 
schemes \cite{HarlowEvans2, cottet_raviart_pic, wang_particle--cell_2011, myers_4th-order_2016}. PIC-based approaches have been widely used in modern gyrokinetic simulations, for example, simulations of Tokamaks like ITER using GTC \cite{GTC-code} or Stellarators like W7X using EUTERPE \cite{EUTERPE}.

Like grid-based solvers, PIC is subject to the curse of dimensionality, but it offers important advantages. 
Particles naturally adapt to evolving dynamics, and the unstructured nature of the data simplifies treatment of 
complex domains, boundary conditions, and parallelization to distributed computing systems. The main drawback is 
numerical noise: to achieve the same accuracy as grid-based methods, significantly higher particle counts are 
required. This becomes problematic in low-density regions or when high accuracy is essential. The majority of 
implementations are explicit PIC codes \cite{Kraus_Kormann_Morrison_Sonnendrücker_2017, pic_on_gpu, piclas, DEROUILLAT2018351, 10046112, PARODI2025102244, HATZKY2007568, Hatzky2002}. 
However, due to the strong mass disparity between electrons and ions as well as between the size of individual-particle scales and the 
domain size, plasma dynamics are highly multi-scale, and explicit schemes demand extremely fine resolution to 
capture all relevant scales. To address this, fully and semi-implicit PIC methods have been developed \cite{BRACKBILL1982271, MASON1981233, Lapenta2006, MARKIDIS20101509, LAPENTA2017349, bacchini_relsim, ARSHAD2025109806, Bacchini_2019}.

This article combines the versatility and numerical efficiency of Lagrangian solvers with the accuracy of semi-Lagrangian approaches by using a flow map approach. 
The knowledge of the flow map allows us to trace the characteristics back up to the initial condition, where we can evaluate the distribution function to any given resolution on variable grids.
The flow map is generally evolved in a PIC-like manner using the Numerical Flow Iteration method (NuFI) \cite{KirchhartWilhelm2024}. NuFI has shown high parallel efficiency \cite{KirchhartWilhelm2024} and versatility concerning boundary conditions, even with low rank approximation \cite{WilhelmKormann2025} and multiple species simulations \cite{Wilhelm_2025}. Nevertheless, tracing the characteristics backward in time in an iterative manner becomes very costly for long simulations, as for each step forward in time, one additional backward step must be performed in the iterative approach. 

An alternative to NuFI is the Characteristic Mapping Method (CMM) \cite{YinMercierYadavSchneiderNave2021}, which performs remapping by the composition of submaps stored as polynomials. The approach allows for larger time steps backward in time and therefore increases the efficiency when backtracing characteristics. The CMM uses a gradient augmented level set method \cite{NaveRosalesSeibold2010} to advance the submaps further in time. It was originally introduced for the 2D incompressible Euler equations \cite{YinMercierYadavSchneiderNave2021}, extended to 3D Euler \cite{YinSchneiderNave2023} and then has been further generalized for Magneto-Hydro-Dynamics (MHD) \cite{YinKrahNaveSchneider2024}, shallow water equations \cite{TaylorNave2023}, and Vlasov--Poisson equations \cite{KrahYinBergmannNaveSchneider2024}.

\subsubsection*{Contribution}

To avoid the time complexity of NuFI, restarting was employed in \cite{WilhelmKormann2025}, which, however 
deteriorates the high-resolution properties of the method. To maintain this, we propose to combine NuFI with CMM. 
This work therefore contributes: a new flow mapping methodology that reduces the time complexity of NuFI,
memory and CPU-time estimations to demonstrate the computational advantage of the combined scheme, and benchmarking 
on classical test cases, including a comparison with an SL predictor-corrector scheme.

\subsubsection*{Organization of the Manuscript}
The remainder of the manuscript is organized as follows. In \cref{sec:nummeth} we introduce numerical methods for computing the flow maps, summarizing NuFI in \cref{subsec:NuFi} and CMM in \cref{subsec:CMM}. The proposed combination of the two methods is presented in \cref{sec:CMM-NuFi}.
Numerical results for Landau damping and the two-stream instability are shown in \cref{sec:results}. \Cref{sec:concl} concludes and gives some perspectives for future work.

\section{Numerical Methods for Calculating Flow Maps}
\label{sec:nummeth}
A kinetic equation (\KE) describes the evolution of a particle distribution function $f(x,v,t)$ in the phase space $(x,v)\in\mathbb{R}^2$, where $x$ and $v$ denote the canonical position and momentum. In general, \KE can be written in the (collisionless) advective form
\begin{equation}
    \frac{\mathrm{d} f}{\mathrm{d} t}(\hat x(t), \hat v(t), t)
    = \frac{\partial f}{\partial t}
    + \dot{\hat x}(t)\frac{\partial f}{\partial x}
    + \dot{\hat v}(t)\frac{\partial f}{\partial v}
    = 0\,.
    \label{eq:advect_form}
\end{equation}
Here, the variables $\hat x(t)$ and $\hat v(t)$ denote the \emph{characteristic trajectories} in phase space, i.e., the solutions of the Hamiltonian equations
\begin{equation}
    \dot{\hat x}(t) = \partial_v H(\hat x(t), \hat v(t)),
    \qquad
    \dot{\hat v}(t) = -\partial_x H(\hat x(t), \hat v(t)),
    \label{eq:Hamilton-eqn}
\end{equation}
with Hamiltonian $H(x,v)$. The initial conditions for these characteristic curves are $(\hat x(0), \hat v(0)) = (x, v)$.
In this work, we consider the 1D+1V Vlasov--Poisson system (\VP) with the Hamiltonian $H(x,v)=\tfrac{1}{2}\lvert v\rvert^2 + q\phi(x)$, $q=1$. Equation~\eqref{eq:advect_form} expresses that the distribution function $f$ remains constant along the characteristic curves $(\hat x(t), \hat v(t))$ defined by~\eqref{eq:Hamilton-eqn}.
The Hamilton flow under which all these trajectories move $\flowmap_{t}^0(x,v)=(\hat x(t),\hat v(t))$ is volume preserving as a consequence of \cref{eq:Hamilton-eqn} divergence-free property (Liouville's Theorem).
Thus the solution to \cref{eq:advect_form} can be described by
\begin{equation}
    f(x,v,t) = f_0(\vec{\Phi}_0^{t}(x,v))\,,
    \label{eq:comp_with_inicond}
\end{equation}
where $f_0$ is the initial distribution and $\flowmap_0^{t}$ the backwards characteristic map.
In other words, the knowledge of $\vec{\Phi}_0^t$ inherits all information needed to reconstruct the solution at any phase space position $(x,v)$.

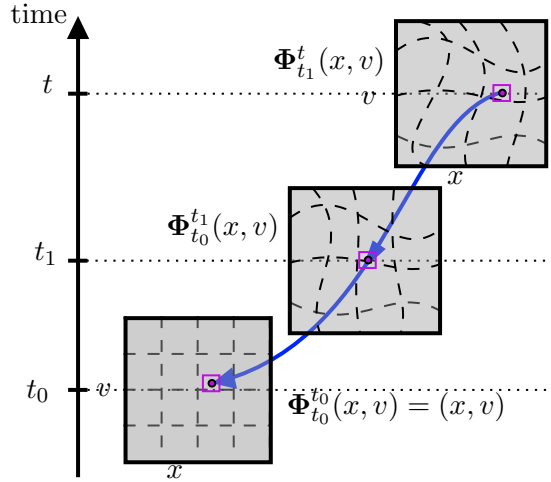
\begin{figure}[htp!]
    \centering
    \input{figures/submaps.tex}
    \caption{Illustration of the backward map, as a composition of two submaps. The blue line illustrates one characteristic curve.}
    \label{fig:submap}
\end{figure}

We remark that the flow map inherits semi-group properties, namely
\begin{equation}
    \vec{\Phi}_{t}^t = \text{Id}_{\mathbb{R}^2}, \quad {\flowmap}_{t_0}^{t_1}\circ \vec{\Phi}_{t_1}^{t_2} = \vec{\Phi}_{t_0}^{t_2}, \qquad (\vec{\Phi}_0^t)^{-1}=\vec{\Phi}_t^0\,,
    \label{eq:semi-group-properties}
\end{equation}
which means that the underlying flow $\vec{\Phi}_t^{0}$ can be subdivided in the so called submaps:
\begin{equation}
    \vec{\Phi}_0^t = \submap_0^\tau\circ\submap_\tau^{2\tau}\circ \dots \circ\submap_{(n-1)\tau}^{t=n\tau}\,.
    \label{eq:comp_structure}
\end{equation}
The subdivision into two maps is illustrated in \cref{fig:submap}.
The submaps obey the same properties as the backward map. However, they are only defined on subintervals.
It was shown previously \cite{KrahYinBergmannNaveSchneider2024} for the Vlasov--Poisson (\VP) equation, that composing the map into submaps and advecting the individual submaps leads to a more efficient representation compared to advecting $f$ itself (see fig. 11 in \cite{KrahYinBergmannNaveSchneider2024}). This is because $\submap$ does not inherit all the fine-scale properties as present in $f$.
Furthermore, in previous works \cite{YinMercierYadavSchneiderNave2021,YinSchneiderNave2023,KrahYinBergmannNaveSchneider2024} the submaps $\submap$ have been stored as Hermite polynomials of $p=3$ order. Thus, making use of the fact that the potentially high-order polynomial representation of $\flowmap$ can be stored more efficiently as compositions of low-order polynomial $\submap$.
It is essentially this \highlight{compositional refinement} that, in contrast to $h$ or $p$ refinement, allows exponential resolution in linear time, enabling the zoom properties shown in \cref{fig:zoom_2stream}.

\subsection{The Numerical Flow Iteration (NuFI)}
\label{subsec:NuFi}

A detailed derivation and discussion of the algorithm can be found in previous work 
\cite{KirchhartWilhelm2024, Wilhelm_2025, nufi_pamm}. In this section, we want to recap the important ideas and highlight the challenges, which we try to tackle in the following chapters.

The core idea of the Numerical Flow Iteration (NuFI) is to approximate the flow map of the 
Vlasov--Poisson equation iteratively instead of directly solving for and storing the distribution functions 
$f$. Therefore, NuFI solves the Hamilton 
\begin{equation}
\label{eqn:ode_char}
\begin{aligned}
\dds \hat{x}(s) &= -\hat{v}(s), & \hat{x}(t) &= x, \\
\dds \hat{v}(s) &=  E\left(s,\hat{x}\left(s\right)\right), & \hat{v}(t) &= v.
\end{aligned}
\end{equation}
using the Strömer--Verlet method to advance the maps iteratively backwards in time
\begin{align}
    \label{eq:Nufiapprox}
    &f(x,v,\NnuFi\tau) = f_0(\flowmap_0^{\NnuFi\tau}(x,v))\approx f_0(\NuFiMap_0^{\NnuFi\tau}(x,v))\\&\text{where}\quad \NuFiMap_0^{\NnuFi\tau}= \NuFiMap_0^\tau \circ \NuFiMap_\tau^{2\tau}\circ \cdots \circ \NuFiMap^{\NnuFi\tau}_{(\NnuFi-1)\tau}
\end{align}
where each submap $\NuFiMap_{k\tau}^{(k-1)\tau}(x,v)$ is defined by:
 \begin{align}
		{v}_{k- 1/2} &= {v} {+} \tfrac{\tau}{2}  E( {x},\tau k), \label{eq:half_step1}\\ 
		{x}_{k-1} &= {x} {-} \tau {v}_{k-1/2}, \\
		{v}_{k-1} &= {v}_{k-1/2} {+} \tfrac{\tau}{2} 
 E({x}_{k-1},\tau(k-1))\\
 \NuFiMap_{(k-1)\tau}^\tau&(x,v):=(x_{k-1},v_{k-1})\,.
	\end{align} 
According to \cite{KirchhartWilhelm2024} the introduced error is
\begin{equation}
    \norm{\flowmap_0^{\Nnufi\tau}(x,v)-\NuFiMap_0^{\Nnufi\tau}(x,v)}_\infty=\Ord{\tau^2}\,.
\end{equation}
Thus, for a smooth initial condition, the error of \cref{eq:Nufiapprox} is given by
\begin{equation}
\label{eqn:eval_f_through_char_num}
\norm{f(x,v,\NnuFi\tau) - f_0(\NuFiMap_0^{\NnuFi\tau}(x,v))}_\infty = \Ord{\tau^2}
\end{equation}
Note that to compute $\rho(x,k\tau) \approx 1-\int f_0(\NuFiMap_0^{k\tau}(x,v))\dv$ via \eqref{eqn:poisson} we can skip the first half-step 
\eqref{eq:half_step1} as it can be recast as a transformation $\tilde{v}=v + \frac\tau2E(x,k\tau)$:
\begin{equation}
    \int f_0(\tilde{\NuFiMap}_0^{k\tau}(x,v+\frac\tau2 E(x,\tau k))\d v =\int f_0(\tilde{\NuFiMap}_0^{k\tau}(x,\tilde{v})) \d \tilde{v} 
\end{equation}
with functional determinant 1. The corresponding NuFI time-stepping algorithm is shown in \cref{algo:Nufi_iteration}. 
\begin{algorithm}[htp!]
    \small
    \begin{algorithmic}[1]
        \renewcommand{\COMMENT}[2][.7\linewidth]{%
        \leavevmode\hfill\makebox[#1][l]{//~#2}}
        \renewcommand{\algorithmicrequire}{\textbf{Input:}}
        \renewcommand{\algorithmicensure}{\textbf{Output:}}
        \REQUIRE{Position and velocity arrays $(x, v)$, electric field history $[E^{(k)}]_{k=0}^{\Nnufi-1}$, time step $\tau$, iteration count $\Nnufi$.}
        \ENSURE{Updated position and velocity arrays $(x, v)$ after $n$ steps with half steps.}
        \IF{$n = 0$}
            \STATE{\textbf{return} $(x, v)$}
        \ENDIF
        \WHILE{$n > 0$}
            \STATE{$n \leftarrow n - 1$}
            \STATE{$x \leftarrow x - \tau \cdot v$ \COMMENT{Position update (inverse signs, going backwards in time)}}
            \STATE{$v \leftarrow v + \tau \cdot \interp[E^n](x)$} \COMMENT{Velocity update, electric field is interpolated}
        \ENDWHILE
        \STATE{$x \leftarrow x - \tau \cdot v$ \COMMENT{Final position step}}
        \STATE{$v \leftarrow v + \frac{\tau}{2} \cdot \interp[E^0](x)$ \COMMENT{Final half velocity step}}
        \RETURN{$(x, v)$}
    \end{algorithmic}
    \caption{NuFI iteration $\texttt{NuFI}((x,v),[E^{(n)}]_n,\tau,\Nnufi)$ to evaluate $\NuFiMap_0^t(x,v)$}
    \label{algo:Nufi_iteration}
\end{algorithm}

\begin{figure}[h!]
    \centering
    \input{figures/NuFi-algo.tex}
    \caption{Graphical sketch of the NuFI algorithm. The red arrows indicate the NuFI iteration \cref{algo:Nufi_iteration} at which all $(x,v)$ are evaluated in a streamline fashion, without storing $f$.}
    \label{fig:nufi_vp}
\end{figure}
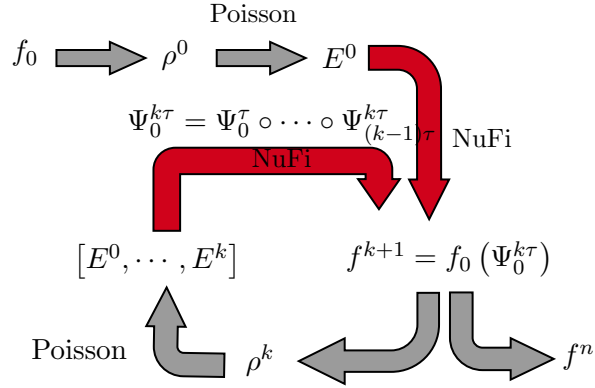

With these components in place, we can outline the full NuFI algorithm for the one-dimensional Vlasov--Poisson system, as depicted in \cref{fig:nufi_vp}. 
After computing the density~$\rho$, the algorithm obtains the electric field by solving $\partial_x E = \rho$ in Fourier space. Using the electric fields stored from all previous time steps, NuFI then applies the iterative procedure described in \cref{algo:Nufi_iteration} to trace characteristics back to their initial positions.

To evaluate the electric field at arbitrary spatial locations, an interpolation operator $\interp[E^n](x)$ is required. Its construction is detailed in \cref{subsec:CMM}, \cref{eq:cmm-interpolation-operator}. Unless stated otherwise, we employ cubic spline interpolation.

For each query point $(x,v)$, the iteration in \cref{algo:Nufi_iteration} is carried out, and the corresponding values of $f$ are accumulated to assemble the density field at the next time step.

Since for each iteration, only the knowledge of the electric fields $E(t_{k},x)$ is necessary, NuFI stores a time series of electric fields to reconstruct the flow map. 
This turns out to be more memory-efficient than storing $f$ itself, especially in higher dimensions $d=3$ on stores a set of $\Nnufi\times$3D interpolants $[E_1,\cdots,E^{\Nnufi}]$ instead of a 6D distribution function $f$ \cite{KirchhartWilhelm2024}. 
Nevertheless, it shifts the memory complexity to computational complexity, as for each step forward in time, the number of NuFI iterations to perform backwards in time increases. The complexity in time is therefore increasing quadratically. To reduce this, we propose to combine it with the characteristic mapping method in a hybrid approach that is detailed in \cref{sec:CMM-NuFi}. 

\subsubsection*{Conservation Properties}
We recall the outstanding conservation properties of NuFI proven in \cite{KirchhartWilhelm2024}.
\begin{theorem}
Let $f(x,v,\tau \NnuFi)=f_0(\NuFiMap_0^{\NnuFi\tau})$ be the discrete solution of NuFI. Then it preserves at all $t=\Nnufi\tau$:
\begin{itemize}
    \item the volume of the flow
    \begin{equation}
        \det \left( \nabla_{(x,v)} \NuFiMap_0^{\Nnufi\tau}(x,v) \right) =1 \qquad \forall (x,v) \in \mathbb{R}\times \mathbb{R}\,,
    \end{equation}
    \item the maximum-principle: $0\le f(x,v,\Nnufi\tau)\le \norm{f_0}_{L^\infty}$\,,
    \item kinetic entropy
    \begin{equation}
        \iint f(x,v,\tau \Nnufi)\ln(f(x,v,\tau \Nnufi)) \d x\d v= \iint f_0 \ln f_0 \d x \d v\,,
    \end{equation}
    \item all integrals for any function $g\colon \mathbb{R}\rightarrow \mathbb{R}$ of the following form {(Liouville's theorem)}:
    \begin{equation}
        \iint g(f(x,v,\tau \Nnufi))\d x \d v = \iint g(f_0(x,v)) \d v \d x\,.
        \label{eq:integral_conserv}
    \end{equation}
\end{itemize}
\label{theo:conserv_nufi}
\end{theorem}

As a consequence of \cref{eq:integral_conserv} implies the conservation of mass.

We further note that NuFI conserves momentum up to errors made in the spatial discretization of the electric field $E$ \cite{KirchhartWilhelm2024}.


\subsection{Characteristic Mapping Method (CMM)}
\label{subsec:CMM}

In contrast to NuFI, the characteristic mapping method (CMM) \cite{MercierYinNave2019} employs compositions based on interpolation to reference back to the initial condition. This more general flow-map framework is not limited to plasma systems but extends to any diffeomorphic flow, including the Euler equations, magnetohydrodynamics and hydrodynamic flows on manifolds \cite{YinKrahNaveSchneider2024,YinMercierYadavSchneiderNave2021,YinSchneiderNave2023,TaylorNave2023}. 

For our numerical experiments we use either Lagrange or spline interpolation on 2D periodic grids $\mapgrid$ with an equal number of grid points in each direction $\Nmapgrid\times \Nmapgrid$. Thus, given the grid data points of a map $\flowmap_{s}^{t}(x_i,v_j), i,j=1,\dots,\Nmapgrid$, we may represent it in polynomial form: 
\begin{equation}
			\vec{\Phi}_{s}^{t}(x,v) \approx \mapcmm_{s}^{t}(x,v)\equiv P_{\mapgrid}[\flowmap_s^t](x,v) = \sum_{i,j} \vec{a}_{ij} \cdot b_i(x) \cdot b_j(v)\,,
            \label{eq:cmm-interpolation-operator}
\end{equation}
where $b_i(x)$  are either barycentric Lagrange or spline basis elements. Respectively for a Lagrange representation, the coefficients are either $\vec{a}_{ij}=\flowmap_s^t(x_i,v_j)$, where $x_i$ and $v_j$ are $m+1$ nearby nodes around the evaluation points $(x,v)$ or spline coefficients $\vec{a}_{ij}$ as a result of solving a linear spline systems.
Depending on the polynomial order $m$ the approximation error introduced is thus: 
		\begin{equation}
			\norm{\Phi_s^t-\mapcmm_s^t}_\infty \leq c \max_{(x,v)} \frac{\norm{\nabla^{(m+1)}_{(x,v)}\Phi_s^t(x,v)}}{\left(m+1\right)!} h^{m+1}\,.
		\end{equation}

In earlier works \cite{KrahYinBergmannNaveSchneider2024}, the backwards map was advanced by a semi-Lagrangian scheme known as the gradient-augmented level set method (GALS) \cite{NaveRosalesSeibold2010}.
GALS solves 
		\begin{align}
			\mapcmm_0^t\colon 
            \begin{cases}
			  \partial_t \mapcmm + {\vec{u}}\cdot \nabla \mapcmm &=0\\
			  \mapcmm(x,v,0) &= (x,v) 
			  \end{cases}\qquad \text{with} \qquad \vec{u}(x,v)=(v,{E}({x},t))
		  \end{align}
backwards in time, which is equivalent to solving \cref{eq:Hamilton-eqn}. Note, to update the electric field $\partial_x E=\int f \d v$ we always have to map back to the initial condition: $f(x,v,t) = f_0(\mapcmm(x,v,t))$.
However, unlike NuFI, GALS is not symplectic and therefore introduces an incompressibility error,
\begin{equation}
\label{eq:incompressibility_error}
\epsilon_\text{incomp} = \abs{\det (\nabla \mapcmm) - 1},.
\end{equation}
As demonstrated for 2D Euler \citep[p.22]{Yadav2015}, this error can induce overshooting, wherein characteristics become too close or cross, compromising the map's injectivity. To control such deviations, a remapping procedure is employed: whenever a map $\mapcmm_{t_{i-1}}^{t_{i}}$ violates a prescribed incompressibility threshold $\delta_\epsilon$, i.e., $\epsilon_\text{incomp}>\delta_\epsilon$, a new map is generated and composed with the preceding maps whenever one needs to access the full map.
We thus exploit the semi-group property \cref{eq:semi-group-properties} of the characteristic map:
    \begin{align}
        \mapcmm_{t_{i-1}}^{t_i} &\colon 
        \begin{cases}
            &\partial_t \mapcmm + {\vec{u}}\cdot \nabla \mapcmm=0 \quad\text{for } t\in[t_{i-1},t_i] \\
            &\mapcmm(x,v,t_{i-1}) =(x,v)
        \end{cases}\quad \text{with}\\
  f(x,v,t)&=f_0(\mapcmm_{t_0}^{t}(x,v)) \text{ and }\mapcmm_{t_0}^{t}=\mapcmm_{t_0}^{t_1}\circ\mapcmm_{t_1}^{t_2}\circ \dots \circ \mapcmm_{t_{m-2}}^{t_{m-1}}\circ\mapcmm_{t_{m-1}}^{t} .\label{eq:cmm-composition}
    \end{align}

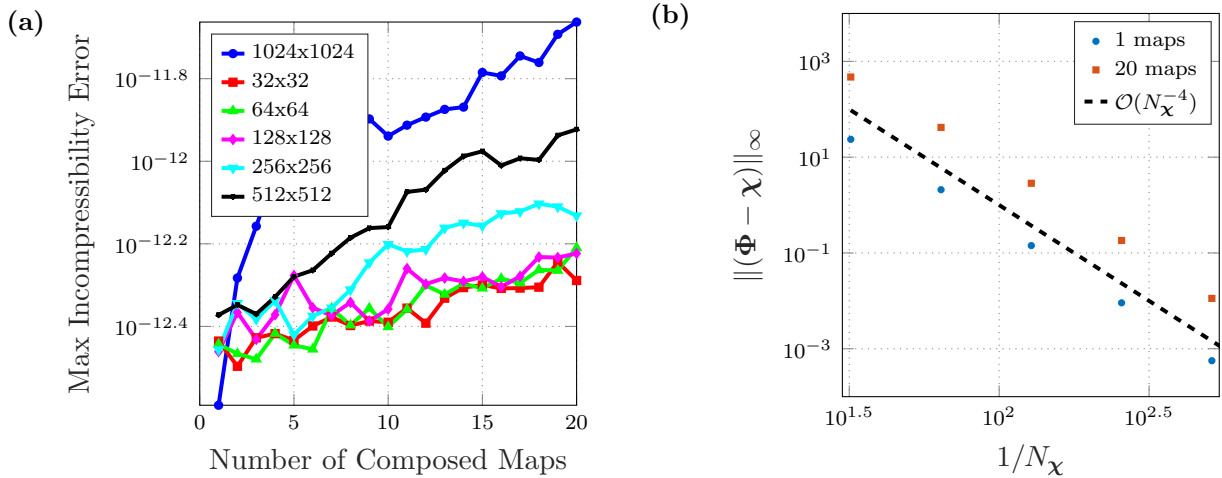
\begin{figure}[htp!]
  \setlength\figureheight{0.3\linewidth}
  \setlength\figurewidth{0.3\linewidth}
    \centering
\begin{subfigure}[b]{0.5\textwidth}
     \caption{}
     \centering
     \label{subfig:volumeperserving_upsample}
     \vspace{-0.3cm}
    \input{figures/volumepreserving_upsample.tex}
    \end{subfigure}\hfill
\begin{subfigure}[b]{0.5\textwidth}
     \caption{}
     \centering
     \label{subfig:l2error_volumeperserving}
     \vspace{-0.3cm}
    \input{figures/volumepreserving_l2error}
    \end{subfigure}
    \caption{Incompressibility error using Lagrange interpolation. a) maximum error as a function of the number of maps. b): Decay of the $L_\infty$-error as a function of the map grid resolution $\Nmapgrid$.}    \label{fig:incomp_error_flow_map}
\end{figure}

\begin{remark}
We remark that even for a symplectic scheme (such as NuFI) \cref{eq:incompressibility_error} will not vanish in most cases, since the finite representation errors of the map will show up in the higher order terms of the approximated derivatives. In other words, \cref{eq:incompressibility_error} not only measures the compressibility but the representation quality of the map in general.
This is intuitively clear as two arbitrarily close characteristics may not be representable on a finite grid. 
This fact is visualized for the consecutive compositions of the volume-preserving submap
\begin{equation}
    \flowmap(x,v) = (x + \sin(\pi v), v)
\end{equation}
in \cref{fig:incomp_error_flow_map}. Although the approximation error decreases with fourth order for our chosen 4th order Lagrangian interpolation scheme. The incompressibility error remains and increases with the increasing number of map decompositions.
\end{remark}

Similarly to NuFI, the characteristic map allows provides exponential resolution in linear memory, because the compositional structure allows for representing exponentially high order polynomials with a linear increase in the number of coefficients.
\begin{example} 
If each map is a cubic polynomial $\mapcmm_{t_{i-1}}^{t_i}\in P_3$, then
	\begin{equation}
		\mapcmm_{t_0}^{t}=\mapcmm_{t_0}^{t_1}\circ\dots\circ\mapcmm_{t_{m-1}}^{t}\in P_3\circ \dots \circ P_3\subset P_{3^m}
	\end{equation}
	Represents an element of $P_{3^m}$ ($3^{m}+1$ coefficients), by storing $(3+1)\cdot m$ coefficients of $P_3$.
\end{example}

As a consequence, the maps can be stored on a coarse grid $\mapgrid$ with $N_{\mapcmm}\times N_{\mapcmm}$ coefficients. To circumvent the high dimensionality, one can later evaluate $\rho(x)=\int f(\mapcmm_{t_0}^{t}(x,v))\mathrm{d}v$ pointwise in a fast way on any sample grid $\samplegrid$ that is fine enough to represent the solution well. This effectively avoids storing $f$ in the active memory, as $f$ can be evaluated in a streamlined fashion using the compositional structure \cref{eq:cmm-composition}.

\begin{algorithm}[!t]
    \small
    \begin{algorithmic}[1]       
        \renewcommand{\COMMENT}[2][.5\linewidth]{%
        \leavevmode\hfill\makebox[#1][l]{//~#2}}
        \renewcommand{\algorithmicrequire}{\textbf{Input:}}
        \renewcommand{\algorithmicensure}{\textbf{Output:}}
        \REQUIRE{Maps $[\mapcmm^{(k)}]_{k=1,\dots,\Nmaps}$ (4D array $N_{\mapcmm}\times N_{\mapcmm}\times 2\times \Nmaps$), coarse grid $\mapgrid$ at which map is stored, 
        $(\mathbf{X},\mathbf{V})$ query coordinates to evaluate the maps on}
        \ENSURE{Composed map $\mapcmm_\text{cmm}(x,v)=\mapcmm_{t_0}^{t_1}\circ\dots\circ \mapcmm_{t_{\Nmaps-1}}^{t_{\Nmaps}}(x,v)$}
                
        \FOR{$i = \Nmaps,\dots, 1$}
            \STATE{$\Delta\mapcmm^{(i)} \leftarrow \mapcmm^{(i)} - (\mathbf{X}_{\mapgrid}, \mathbf{V}_{\mapgrid})$}
            \COMMENT{Compute deviation of map $i$ around the identity}    
            
            \STATE{$(\Delta\mathbf{X},\Delta\mathbf{V}) \leftarrow \mathcal{P}_{\Delta \mapcmm^{i}}(\mathbf{X},\mathbf{V})$}
            \COMMENT{Interpolate displacement fields at query points}

            \STATE{$\mathbf{X} \leftarrow \mathbf{X} + \Delta\mathbf{X}$}
            \COMMENT{Apply transformation:}
            
            \STATE{$\mathbf{V} \leftarrow \mathbf{V} + \Delta\mathbf{V}$}
            \COMMENT{new point = old point + displacement}
        \ENDFOR
        
        \RETURN{$(\mathbf{X},\mathbf{V})$}
        \COMMENT{Return final composed map}
    \end{algorithmic}
    \caption{Numerical composition of maps \texttt{compose\_maps}($[\mapcmm^{(k)}]_{k=1}^{\Nmaps},(\mathbf{X},\mathbf{V}))$}
    \label{alg:compose_maps_numerical}
\end{algorithm}

\subsection{Hybrid CMM-NuFI method} 
\label{sec:CMM-NuFi}

The proposed semi-Lagrangian flow map approach combines NuFI (\cref{subsec:NuFi}) and CMM (\cref{subsec:CMM}). 
The idea is visualized in \cref{fig:flowmapdecomposition} and the algorithm detailed in \cref{algo:Nufi-halfv2,alg:CMM-NuFi}.
Depending on the remapping strategy, the solution is locally evolved using NuFI iterations and after a certain number of iterations $\Nnufi$ the generated iterative map $\NuFiMap^{t_i+\Nnufi\tau}_{t_i}$ is replaced by a stored submap $\mapcmm^{t_{i+1}}_{t_{i}}, t_{i+1}=t_i+\Nnufi\tau$. The number of NuFI iterations $\Nnufi$ can be either chosen statically by a remapping frequency or adapted to the solution itself as suggested in \cite{YinMercierYadavSchneiderNave2021}. In this work, we only make use of manual remapping, as we are interested in optimizing the time complexity of NuFI.
\begin{figure}[htp!]
    \centering
    \input{figures/hybrid-scheme.tex}
    \caption{Schematic of the backward flow $(x^*,v^*)=\flowmap_0^t(x,v)$ for two consecutive iterations with the hybrid CMM-NuFI method. 
    Given a point in phase space $(x,v)$, the characteristics are traced back to their footpoints $(x^*,v^*)$ using an iterative flowmaps $\NuFiMap$ (green arrows) and compositional flowmaps $\mapcmm$ (blue arrows). The two consecutive iterations visualize the remapping. After four NuFI steps (green), the resulting map is replaced in the next iteration by a single large submap (dark blue) and concatenated with the rest of the submaps.}
    \label{fig:flowmapdecomposition}
\end{figure}
Hence, to represent the solution $f(x,v,t)$ at time $t = t_\Nmaps + \Nnufi\tau$ after $\Nmaps$ remappings and $\Nnufi$ NuFI iterations, we use the following decomposition:
	\begin{align}
		f(x,v,t) &= f_0(\flowmap_0^t(x,v)),\qquad \text{where}\\
		\flowmap^t_0 &\approx \tilde{\flowmap}^t_0 :=\underbrace{\mapcmm_0^{t_1} \circ \mapcmm^{t_2}_{t_1} \circ \dots\circ \mapcmm^{t_{\Nmaps}}_{t_{\Nmaps-1}}}_{=\mapcmm_\text{CMM}}
		\circ 
		\underbrace{\NuFiMap^{t_{\Nmaps}+\tau}_{t_{\Nmaps}}\circ \NuFiMap^{t_{\Nmaps}+2\tau}_{t_{\Nmaps}+\tau}\circ\dots\circ \NuFiMap^{t_{\Nmaps}+\Nnufi\tau}_{t_{\Nmaps}+(\Nnufi-1)\tau}}_{=\NuFiMap_\text{Nufi}}\,.
		\label{eq:compose_hybrid}
	\end{align}
The convergence order of the hybrid method $\tilde{\flowmap}^t_0 :=\mapcmm_\text{CMM}\circ{\NuFiMap_\text{Nufi}}$ is given by the time error $\Ord{\tau^2}$ coming from $\NuFiMap_\text{NuFi}$ and the interpolation errors $\Ord{1/\Nmapgrid^\alpha}$ of $\mapcmm_\text{CMM}$:
	\begin{equation}
		\Vert \flowmap^{t}_0 - \tilde{\flowmap}^t_0 \Vert_\infty = \mathcal{O}(\tau^2) + \mathcal{O}(1/\Nmapgrid^\alpha)\,.
        \label{eq:flowmap_error}
	\end{equation}
	Same holds for the solution $\norm{f_0(\flowmap^t_0)-f_0(\tilde{\flowmap}^t_0)}_\infty$ if the initial condition $f_0$ is sufficiently smooth. The spatial convergence of our scheme is visualized for the Landau damping in \cref{fig:convergence_study}.

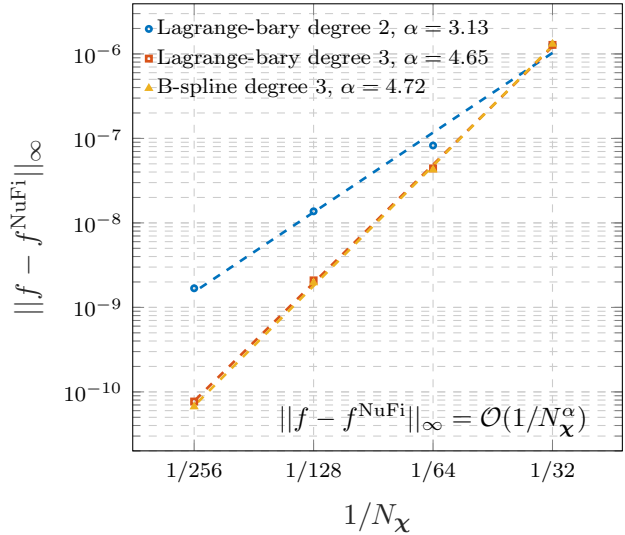
\begin{figure}
    \centering
    \setlength\figureheight{0.35\linewidth}
  \setlength\figurewidth{0.4\linewidth}
    \input{figures/convergence_comparison_landau_damping.tikz}    \caption{Spatial convergence for the Landau damping test case at $T=2$ with $\tau=0.1$ using different interpolation schemes. The reference solution is pure NuFI with a spatial resolution $\Nsample=512$.}
    \label{fig:convergence_study}
\end{figure}

In the proposed scheme, we use NuFI for local time stepping with step size $\tau$ and the CMM submap composition \cref{alg:compose_maps_numerical} to relate back the solution to the initial footpoints. Note, that the NuFI \cref{algo:Nufi-halfv2} has to be adapted for $\Nmaps>1$, avoiding the final half step for consistency with the previous maps.

\begin{algorithm}[!t]
    \small
    \begin{algorithmic}[1]
        \renewcommand{\COMMENT}[2][.5\linewidth]{%
        \leavevmode\hfill\makebox[#1][l]{//~#2}}
        \renewcommand{\algorithmicrequire}{\textbf{Input:}}
        \renewcommand{\algorithmicensure}{\textbf{Output:}}
        \REQUIRE{Position and velocity arrays $(x, v)$, electric field history $[E^{(k)}]_{k=0}^{\Nnufi-1}$, time step $\tau$, iteration count $\Nnufi$.}
        \ENSURE{Footpoints $(x, v)$ going back in time $\Nnufi\tau$.}
        \IF{$n = 0$}
            \STATE{\textbf{return} $(x, v)$}
        \ENDIF
        \IF{$N_{\text{maps}} = 0$}
            \STATE{$(x,v) = \texttt{NuFI}((x,v),[E^{(n)}]_n,\tau,\Nnufi)$}
        \ELSE
            \WHILE{$n > 1$}
                \STATE{$n \leftarrow n - 1$}
                \STATE{$x \leftarrow x - \tau v$}
                \STATE{$v \leftarrow v + \tau \interp[E^n](x)$}
            \ENDWHILE
        \ENDIF
        \RETURN{$(x, v)$}
    \end{algorithmic}
    \caption{NuFI half step adjusted \texttt{NuFI\_adj}$((x,v),[E^{(n)}]_n,\tau,\Nnufi)$}
    \label{algo:Nufi-halfv2}
\end{algorithm}

\begin{algorithm}[!t]
    \small
    \begin{algorithmic}[1]       
        \renewcommand{\COMMENT}[2][.5\linewidth]{%
        \leavevmode\hfill\makebox[#1][l]{//~#2}}
        \renewcommand{\algorithmicrequire}{\textbf{Input:}}
			\renewcommand{\algorithmicensure}{\textbf{Output:}}
			\REQUIRE{Submaps $[\mapcmm^{(k)}]_{k=1,\dots,M}$Electric field $[{E}^{(n)}]_{n=1,\dots,N}$, time step $\tau$.}
			\ENSURE{Updated map list $[\mapcmm^{(k)}]_{k}$ and electric field list $[{E}^{(n)}]_{n}$.}
			
			\STATE{$\Nnufi\leftarrow \Nnufi + 1$}
			\COMMENT{increase NuFI iteration counter}
				\STATE{Initialize:\\ $X_{ij}={i\Delta x}$,\\$V_{ij}=j\Delta{v}-L_v/2$ for $i,j=1,\dots,\Nsample$}
				\COMMENT{Sample Grid $\Nsample$}
				\STATE{$(X, V) =\texttt{NuFI\_adj}((X, V), [{E}^{(n)}]_{n}, \tau, \Nnufi)$}
				\COMMENT{implements $\NuFiMap_{t_{\Nmaps}}^{t_{\Nmaps}+\Nnufi\tau}$}
				\STATE{$(X^*, V^*) = \texttt{MapsCompose}([\mapcmm^{(k)}],(X, V))$}
				\COMMENT{Implements composition of maps}
				\STATE{$\scalar^{\text{new}}(X, V) = \scalar_{0}(X^*, V^*)$}
				\COMMENT{Compose with initial condition}
				\STATE{Compute Density $\rho=\int f^{(\text{new})} \mathrm{d}v -1$}
			\IF{$\Nnufi \bmod \Nremap = 0$}
				\STATE{$\Nmaps\leftarrow \Nmaps + 1$}
				\STATE{$\mapcmm^{\Nmaps} = \texttt{down\_sample}(X,V,{\mapgrid})$} \label{alg:hybrid-downsample}
				\COMMENT{add coarse submap to $[\mapcmm^{(k)}]$}
				\STATE{ $\Nnufi \leftarrow 1$}
				\COMMENT{reset NuFI iteration counter}
			\ENDIF
			\STATE{Solve Poisson: $\nabla^2 \phi = \rho$}
			\STATE{Set ${E}^{(\Nnufi)} = -\nabla \phi$}
		\end{algorithmic}

    \caption{Hybrid CMM-NuFI time step \texttt{CMM-NuFI}($[\mapcmm^{(k)}]_k, [E^{(n)}]_{n}, \tau$)}
    \label{alg:CMM-NuFi}
\end{algorithm}
Summarizing multiple NuFI steps into a big CMM step with step size $\Delta t = \Nnufi \tau$ allows us to reduce NuFI's quadratic complexity in time by several orders of magnitude.
This is shown in \cref{fig:cpu_time_analysis} for the cumulative CPU time, where we compare the hybrid scheme against NuFI using different remapping frequencies $\Nremap$. The typical behavior of the hybrid method can be observed, where the complexity per iteration of the method decreases, as fewer NuFI iterations are required per time step. This resetting of complexity creates a sawtooth function in \cref{fig:cpu_time_analysis}, resulting in a reduction of the overall CPU time by several orders of magnitude. This complexity reduction enables us to run much longer plasma simulations. 
Further note that a sweet spot between iterative and compositional mapping is achieved for $\Nremap=20$ iterations between remappings. This sweet spot is reached when the complexity of the map iteration and composition is equal, $t_\text{CPU}(\flowmap_\text{CMM}) = t_\text{CPU}(\NuFiMap_\text{NuFi})$. It might be dependent on the parallelization strategy of the NuFI iterations and CMM compositions, the problem size, and the length of the simulation, which affects the total number of CMM map compositions. 
As in the original CMM method \cref{subsec:CMM} we store $\mapcmm$ on a coarse grid $\mapgrid$ for computational efficiency. Therefore we only store $\Nmapgrid\ll \Nsample$ sample points when performing remapping in \cref{alg:CMM-NuFi} line \ref{alg:hybrid-downsample}. 

\subsubsection*{Memory and Computational Complexity}  
\begin{figure}[hbp!]
    \centering
    
    \includegraphics[width=0.9\linewidth]{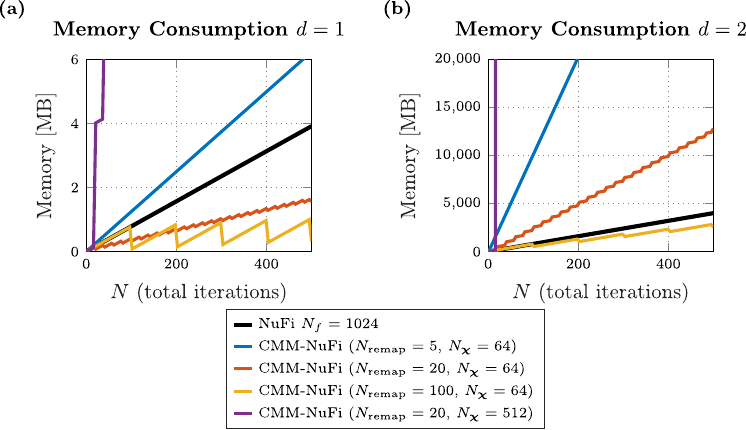}
    \caption{Memory vs iteration for various grid resolutions $\Nmapgrid$, remapping frequencies $\Nremap$ and dimensionality $d=1,2$. For all examples, the sample grid resolution is $\Nsample=1024$.}
\label{fig:memory_vs_iteration}
\end{figure}

For estimating the computational costs, we assume that the initial condition $f_0$ is given analytically and therefore does not require memory. Furthermore, we assume that the sampling and map grid are resolved with the same number of grid points in $x$ and $v$ direction $N_{v}=N_x$. However, the resolution of the sampling grid and map grid are different $\Nmapgrid\le\Nsample$. Furthermore, we assume that the electric field is stored on the fine sample grid $\samplegrid$, which might not be necessary, as it is a smooth function.
Following \cite{KirchhartWilhelm2024}, the memory usage is given by:
\begin{equation}
    \text{Memory NuFI} = (\Nnufi +1) \Nsample^d \times 8 \text{ bytes}\,.
\end{equation}
in terms of the dimension $d$ (in this study $d=1$).
Similarly, the storage for $\Nmaps$ CMM maps is
\begin{equation}
    \text{Memory CMM} = (\Nmaps) 2d\Nmapgrid^{2d} \times 8 \text{ bytes}\,.
\end{equation}
Choosing $\Nmaps=\lfloor\frac{\Nnufi}{\Nremap} \rfloor$ where $\Nnufi$ is assumed to be the total iteration number, we can deduce the hybrid memory storage as:
\begin{equation}
    \text{Memory CMM-NuFI} = \left[(\Nmaps) 2d\Nmapgrid^{2d} + (\Nnufi-\Nmaps\Nremap+1) \Nsample^d\right]\times 8 \text{ bytes}\,.
\end{equation}
In \cref{fig:memory_vs_iteration} we plot the memory consumption of NuFI vs. the number of iterations and compare it to the hybrid approach for different remapping frequencies $\Nremap$, map grid resolutions $\Nmapgrid$.
It can be seen that the memory consumption of both schemes increases per iteration. The slope of the memory increase is influenced by the dimensionality and grid resolution. The typical sawtooth structure of the memory consumption is best visible for the configuration ($\Nremap=100,\Nmapgrid=64,d=1$) in \cref{fig:memory_vs_iteration}. The sawtooth is due to the remapping, which shifts the storage of electric fields $E$ to the storage of maps $\mapcmm$.

In the same manner for the computational complexity, NuFI's cost scales as \cite{KirchhartWilhelm2024}:
\begin{equation}
    \text{Cost NuFI} = \Nsample^{2d}\sum_{k=1}^{\Nnufi}(k+1)=\Nsample^{2d}(\Nnufi+3)\frac{\Nnufi}{2}\,.
\end{equation}
assuming that NuFI is evaluated at the sampling grid $\samplegrid$.
The costs for interpolating submaps depend on the interpolation scheme. Here we assume that the costs for evaluating the polynomial of order $\alpha$ is dominant over the setup of the coefficients and compose $\Nmaps$ CMM maps
\begin{equation}
    \text{Cost CMM composition} = [(\alpha+1) \Nsample]^{2d}\times \Nmaps\,. 
\end{equation}
Choosing as before $\Nmaps=\lfloor\frac{\Nnufi}{\Nremap} \rfloor$, with $\Nnufi$ the total iteration number, we can deduce the hybrid computational complexity as:
\begin{equation}
    \text{Cost CMM-NuFI} = \Nsample^{2d}[(\Nnufi-\Nmaps\Nremap+3)\frac{\Nnufi-\Nmaps\Nremap}{2} + (\alpha+1)^{2d}\times \Nmaps]\,.
\end{equation}
For a typical example like Landau damping, the computational costs are plotted in \cref{fig:cpu_time_analysis}.
It can be seen that the cumulative CPU time scales quadratically both for NuFI and CMM-NuFI, however, the slope is reduced in the hybrid method.


\section{Results} 
\label{sec:results}

In the following section, we benchmark the proposed hybrid method using two classical test cases: linear Landau Damping in \cref{subsec:LandauDamping} and two-stream instability in \cref{subsec:two-stream}.
For our comparison, we recall the conserved properties that we studied in our test cases in a solution domain $\Omega_x\times\Omega_v$.
These are the \define{mass (0. moment)}
\begin{equation}
    \mass(t) =\int_{\Omega_x} \density(x,t) \d x = \iint_{\Omega_x\times \Omega_v}  f(x,v,t) \d x\, \d v\,,
\end{equation}
and its \define{moment (1. moment)}
\begin{equation}
    \momentum(t) =\int_{\Omega_x} \flux(x,t) \d x = \iint_{\Omega_x\times \Omega_v}  vf(x,v,t) \d x\, \d v\,.
\end{equation}
Furthermore, the \define{total energy} is conserved
\begin{equation}
    \energy(t) =\Epot(t)+ \Ekin(t)\,,
\end{equation}
which is the sum of \define{kinetic energy}
\begin{equation}
    \Ekin(t) =\frac{1}{2}\int_{\Omega_v} \int_{\Omega_x}  f(x,v,t) \abs{v}^2 \d x\, \d v\,,
\end{equation}
and \define{potential energy}
\begin{equation}
    \Epot(t) =\frac{1}{2}\int_{\Omega_x} \abs{E(x,t)}^2 \d x = \frac{1}{2}\int\int_{\Omega_x} \rho(t,x) \phi(t,x) \d{x}\,.
\end{equation}

\subsection{Landau Damping}
\label{subsec:LandauDamping}

Here we study linear Landau damping for 1D+1V Vlasov--Poisson equation
with a perturbed Maxwellian distribution
as initial condition:
\begin{equation}
   f_0(x,v) = (1+\epsilon\cos(kx))\frac{1}{\sqrt{2\pi}} e^{-v^2/2}\quad\text{with }  \epsilon = 0.01,k=0.5\,.
\end{equation}
For time stepping, we use $\tau=0.1$ and a constant remapping frequency of $\Nremap=20$. For the numerical discretization of the CMM maps, we use Lagrange interpolation of 3rd order with an equal number of $\Nmapgrid$  grid points in each direction. The distribution function $f$ is sampled using $\Nsample$ grid points in each dimension. Thus, the electric fields and the density are represented on $\Nsample$ points.

\subsubsection*{Damping rate}
\begin{figure}[htp!]
		\centering
		\setlength\figureheight{0.4\linewidth}
	  \setlength\figurewidth{0.65\linewidth}
		\input{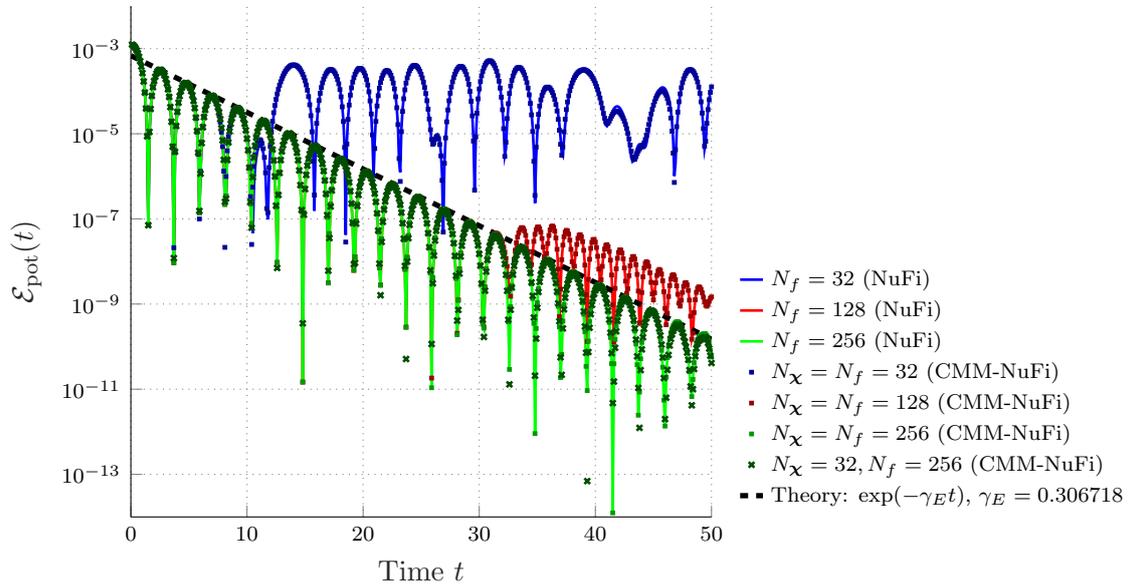}  
        \caption{Linear Landau damping for NuFI and CMM-NuFI studied for different sample grid resolutions $\Nsample$ and map grid resolutions $\Nmapgrid$, time step $\tau=0.1$ and remapping frequency $\Nremap=20$ for CMM-NuFI.}
		\label{fig:linear_landau_damping}
	\end{figure}

In \cref{fig:linear_landau_damping}, we compare the previously developed NuFI method \cite{KirchhartWilhelm2024} with the proposed method concerning the decay of the potential energy and compare it to the theoretical prediction. 
For all different resolutions studied, we see perfect overlap of the decay of the potential energy for both methods. Similarly to the observations in \cite{KrahYinBergmannNaveSchneider2024,KirchhartWilhelm2024} we observe that increasing the number of sampling grid points allows us to resolve lower potential energies. This is due to the filamentations that need to be resolved when solving the Poisson equation. Moreover, we note that the resolution of the map grid does not have visible influence on the evolution of $\Epot$ shown in \cref{fig:linear_landau_damping}. A map grid with $\Nmapgrid = 32$ points combined with a sample grid with $\Nsample = 256$ produces indistinguishable results from a simulation in which both grids use the same resolution, $\mapgrid = \samplegrid$ with $\Nmapgrid = \Nsample = 256$. This behavior is consistent with the findings reported in \cite{KrahYinBergmannNaveSchneider2024}.

\subsubsection*{Conservational Properties}
\begin{figure}[htp!]
    \centering  
  \setlength\figureheight{0.8\linewidth}
  \setlength\figurewidth{0.6\linewidth}
  \input{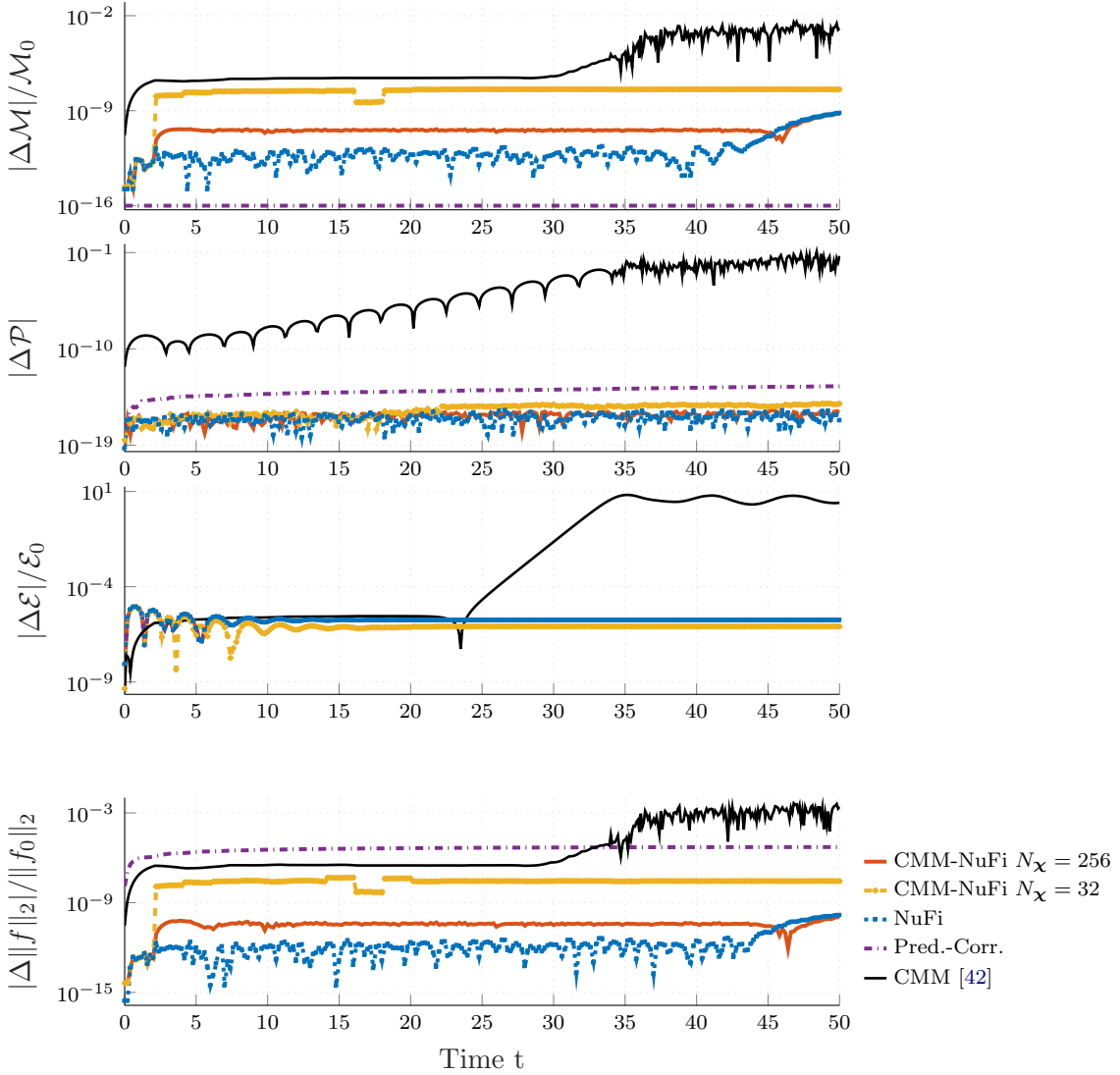}
    \caption{Conservation for Landau damping using a sample grid with $\Nsample=256$ and time step $\tau = 0.1$. For the CMM method, we use adaptive remapping with incompressibility threshold (incomp. threshold) $\delta_\epsilon=10^{-2}$, i.e., $\epsilon_\text{incomp}<\delta_\epsilon$, and CMM-NuFI is remapped after $\Nremap=20$ NuFI iterations. }
    \label{fig:conservation_landau_damping}
\end{figure}

The conservation of mass, momentum, total energy and the $L_2$-norm of the distribution function have been tested and compared in \cref{fig:conservation_landau_damping} against the previously developed CMM method \cite{KrahYinBergmannNaveSchneider2024}, that is, using the gradient augmented level-set method \cite{NaveRosalesSeibold2010} to solve \cref{eq:cmm-composition}, NuFI \cite{KirchhartWilhelm2024} and a classical backward semi-Lagrangian predictor-corrector method detailed in \cite{BourneMunschyGrandgirardMehrenbergerGhendrih2023}.

\begin{figure}[htp!]
    \centering  
  \setlength\figureheight{0.8\linewidth}
  \setlength\figurewidth{0.6\linewidth}
  \input{figures/incomp_err_landau_damping}
    \caption{Incompressibility error for Landau damping using a sample grid with $\Nsample=256$ and time step $\tau = 0.1$. For the CMM method, we use adaptive remapping $\epsilon_\text{incomp}<10^{-2}$ and CMM-NuFI is remapped after $\Nremap=20$ NuFI iterations. }
    \label{fig:incomp_err_landau_damping}
\end{figure}
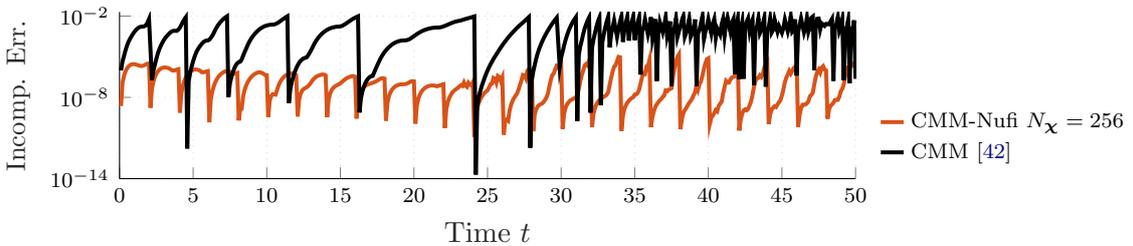

Additionally, we show the time evolution of the incompressibility error \cref{eq:incompressibility_error} in \cref{fig:incomp_err_landau_damping}.
If not otherwise stated, we use Lagrange interpolation of 3rd order for all three methods, time step $\tau=0.1$, and a sampling grid size of $\Nsample=256$ in each dimension $x,v$.

The time evolution of the difference with the initial condition shows that the hybrid CMM-NuFI and the classical NuFI methods overlap before the first remapping. After introducing the first CMM map, small deviations between the two methods are visible in the mass and $L_2$ conservation, which can be explained as the interpolation scheme is not volume preserving. Similarly, increasing the introduced interpolation error with a coarser map resolution (here $\Nmapgrid=32$) leads to larger errors in the volume preservation, which can be seen in the mass and $L_2$ conservation. These errors jump when a CMM map is introduced, as can be seen when comparing the remapping intervals in \cref{fig:incomp_err_landau_damping} with the jumps in mass conservation,\cref{fig:conservation_landau_damping} around $t=15$.
Nevertheless, we observe an advantage over the classical CMM method, which is expected as the GALS is not volume preserving. Comparing the flow-mapping methods with the predictor-corrector scheme, we see that the 3rd order predictor-corrector outperforms the others concerning mass conservation. On the other side, the flow mapping methods show several orders of magnitude better preservation of the $L_2$-norm.

\subsubsection*{Performance Improvements}

\begin{figure}[htp!]
    \centering
    \begin{subfigure}[b]{0.48\textwidth}
  \setlength\figureheight{0.58\linewidth}
  \setlength\figurewidth{0.75\linewidth}
     \caption{}
     \centering
     \label{subfig:cpu_per_iter}
     \vspace{-0.3cm}
    \input{figures/cpu_timing_evolution}
\end{subfigure}
\begin{subfigure}[b]{0.48\textwidth}
  \setlength\figureheight{0.58\linewidth}
  \setlength\figurewidth{0.75\linewidth}
     \caption{}
     \centering
     \label{subfig:cummulative_cpu_time}
     \vspace{-0.3cm}
    \input{figures/cpu_timing_cumulative}
\end{subfigure}

    \caption{CPU time per iteration (a) and cumulative CPU time (b) using a sample grid with $\Nsample=256$ and time step $\tau = 0.1$. Simulations have been performed on ARM64.}
    \label{fig:cpu_time_analysis}
\end{figure}
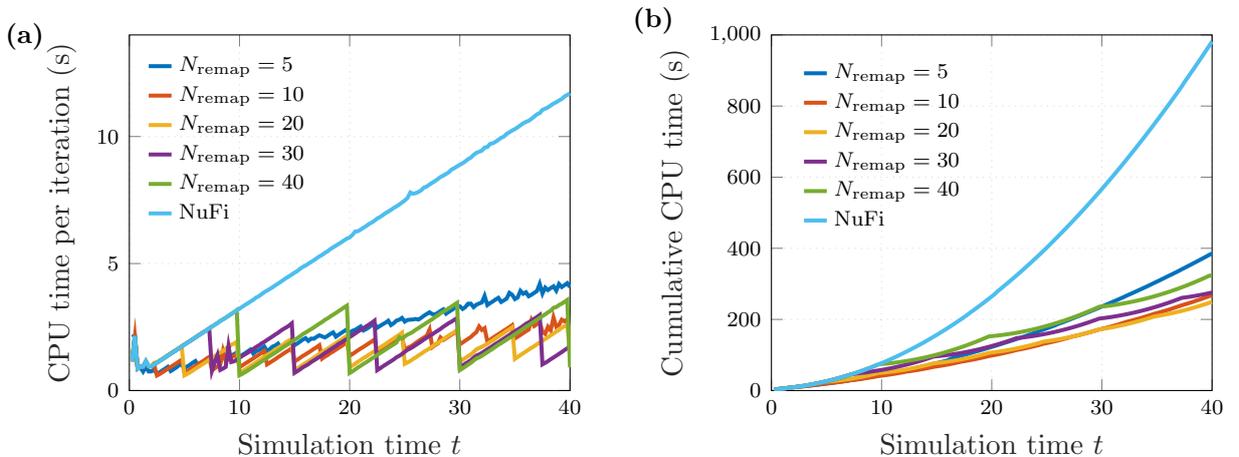

Besides improving the conservation properties of the classical CMM approach, the hybrid method outperforms NuFI concerning computational efficiency. This is graphically visualized in \cref{fig:cpu_time_analysis} for the CPU time per iteration and cumulative CPU time. We use a sample grid size of $\Nsample=256$ and a time stepsize of $\tau = 0.1$. Depending on the remapping frequency $\Nremap$, we observe different behavior of our method. 
If no remapping is performed (plain NuFI), the CPU time per iteration increases linearly with time, as for each new time step, an additional backward step is necessary. This results in a quadratic increase in the cumulative CPU time, visible in \cref{subfig:cummulative_cpu_time} for the NuFI iteration.
If we introduce remapping, the CPU time per iteration decreases drastically as the number of NuFI iterations is replaced by an interpolation. The characteristic sawtooth structure is shown for $\Nremap=10-40$. If we remap rapidly with $\Nremap=5$ we also see a linear increase in CPU time, which is due to the increasing number of compositions to construct the CMM backward map using interpolation.
Although the resulting slope for a low remapping frequency is also quadratic, it is less steep as pure NuFI. The slope of the two extremes depends on the software realization, size of the grids $(\Nsample,\Nmapgrid)$ and dimension of the problem. 
Nevertheless, for our test case, we found that $\Nremap=20$ gives the best compromise between performance and introduced incompressibility error shown in \cref{fig:incomp_err_landau_damping}. For the Landau damping presented here, we observe that the cumulative CPU time is indeed the lowest for $\Nremap=20$ at $t=40$.

\subsubsection*{Spatial Convergence}
The convergence of our scheme is studied in  \cref{fig:convergence_study} for different grid resolutions $\Nmapgrid = \Nsample=32,64,128,256$ and interpolation schemes to confirm \cref{eq:flowmap_error}. In \cref{fig:convergence_study}
we upsample the solution of the NuFI-CMM simulations to our reference solution grid $512\times 512$ that is computed with NuFI and we calculate the maximums error. To minimize the time error, we use for both the reference and the CMM-NuFI simulation a small time step size of $\tau = 0.01$ and we compare our solution at $T=10$. The remapping frequency was chosen as before $\Nremap=20$.

\subsection{Two-stream instability}
\label{subsec:two-stream}

\begin{figure}[htp!]
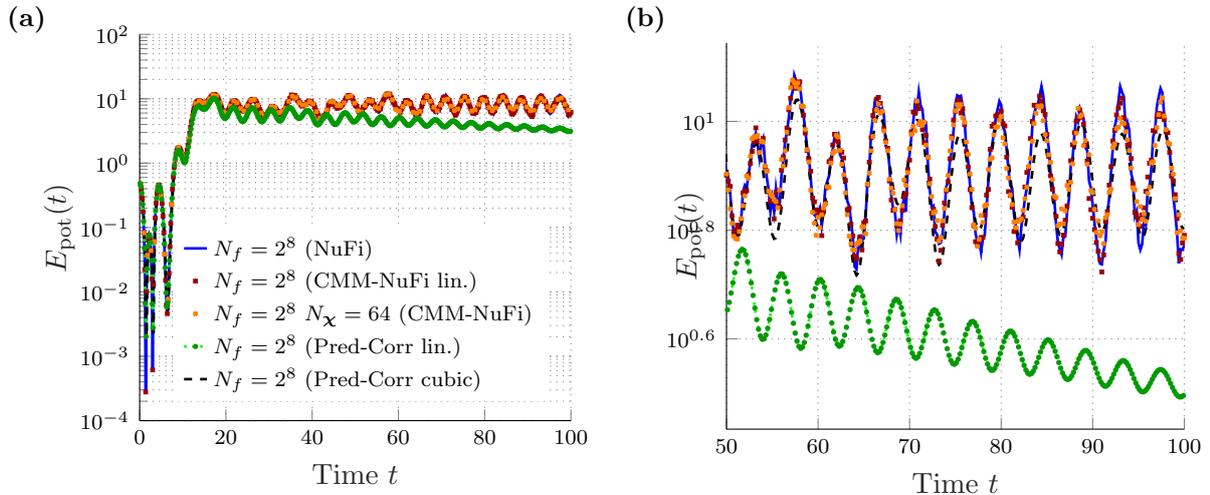

	  \centering
\begin{subfigure}[t]{0.48\textwidth}
     \caption{}
     \centering
	\setlength\figureheight{0.63\linewidth}
	\setlength\figurewidth{1\linewidth}
     \label{subfig:epot_2stream}
     \vspace{-0.3cm}
	  \input{figures/two_stream_resolution_study_Tend100}
	  \end{subfigure}%
\begin{subfigure}[t]{0.48\textwidth}
     \caption{}
     \label{subfig:epot_2streamzoom}
     \centering
	  \setlength\figureheight{0.63\linewidth}
	\setlength\figurewidth{1	\linewidth}
    \vspace{-0.3cm}
    \hspace*{-1.2cm}
    \input{figures/two_stream_resolution_study_Tend100_zoom}%
      \end{subfigure}
      \caption{Evolution of the potential energy as a function of time for NuFI, CMM-NuFI and the predictor-corrector scheme (a) and a zoom in the non-linear regime (b).}
      \label{fig:epot_2stream}
\end{figure}

The two-stream instability is a classical test case for studying fine-scale behavior in phase space.
The initial condition is given by
	\begin{equation}
		f_0(x,v) = (1+\epsilon\cos(kx))\frac{1}{2\sqrt{2\pi}} (e^{-(v-v_0)^2/2}+e^{-(v+v_0)^2/2})\quad \text{with } \epsilon=0.05,k = 0.2,v_0=3 .
	\end{equation}
We compare the evolution of the potential energy for the three different schemes in \cref{fig:epot_2stream}. To this end, we use for all the schemes a time step of $\tau=0.2$. Furthermore, remapping is performed in the CMM-NuFI method after $\Nremap=20$ NuFI steps. As reported in \cite{WilhelmKormann2025}, we observe diffusive behavior of the classical semi-Lagrangian scheme, which implies a variation of the total amplitude (maximum principle in \cref{theo:conserv_nufi}) and results in a decay of the potential energy over time. This is especially pronounced in the linear predictor-corrector visualized in \cref{fig:epot_2stream}.
As in the Landau damping case shown in \cref{fig:linear_landau_damping}, we see no change in the potential energy when varying the mapgrid resolution $\Nmapgrid$ ($\Nmapgrid < \Nsample$) or when using different interpolation orders.

\begin{figure}[htp!]
	\setlength\figureheight{0.3\linewidth}
	\setlength\figurewidth{0.45	\linewidth}
	  \centering
    \includegraphics[width=1\linewidth]{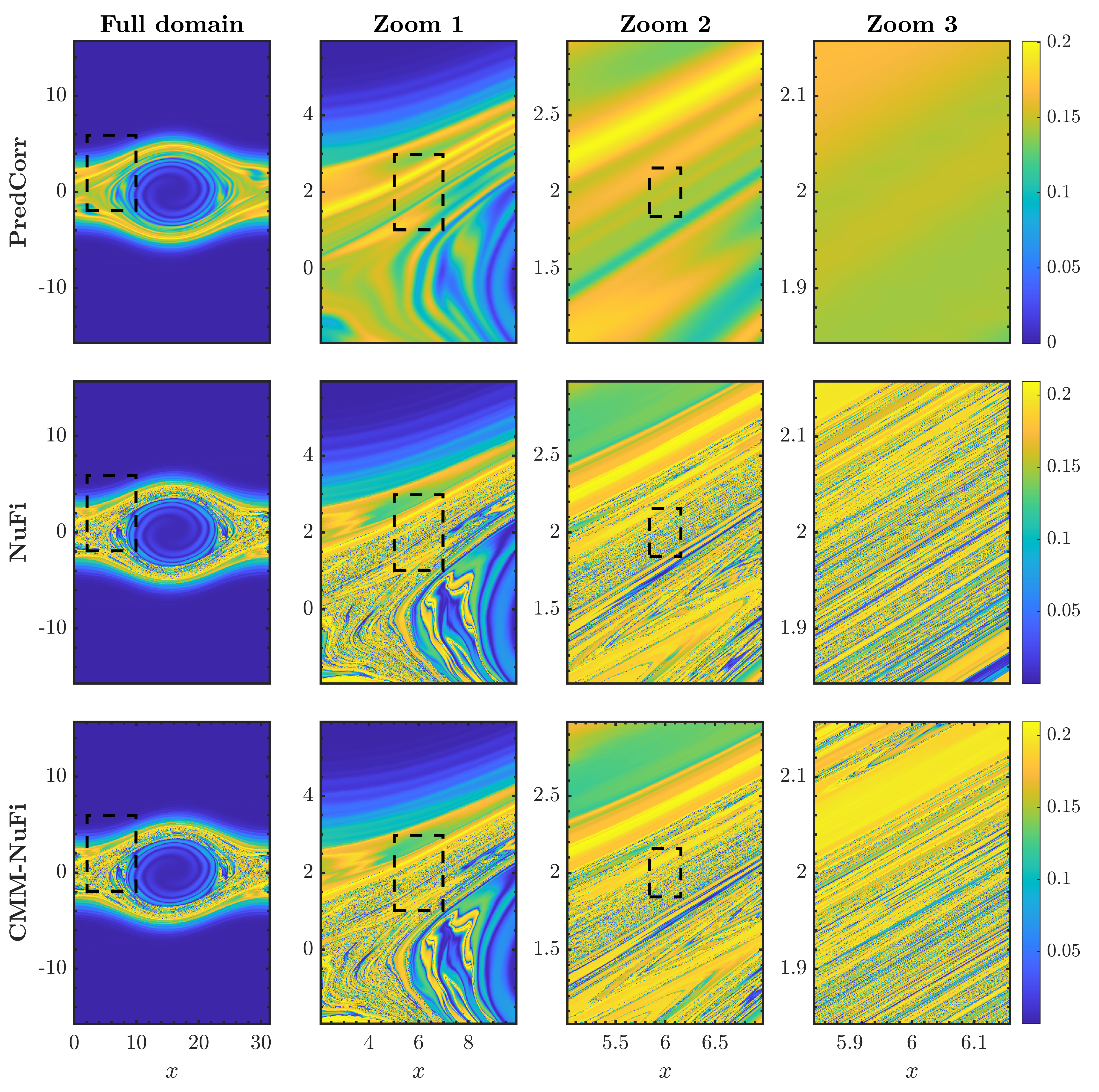}
		    \caption{Zoom properties for the two-stream instability at $t=100$. We compare three zooms: NuFI, CMM-NuFI, and the predictor-corrector scheme. Each zoom is evaluated on a $1024\times 1024$ grid. The results from the predictor-corrector simulation are interpolated inside the zoom windows using cubic splines.}
        \label{fig:zoom_2stream}
        \end{figure}

\begin{figure}[htp!]
    \centering  
  \setlength\figureheight{0.7\textwidth}
  \setlength\figurewidth{0.6\textwidth}
  \input{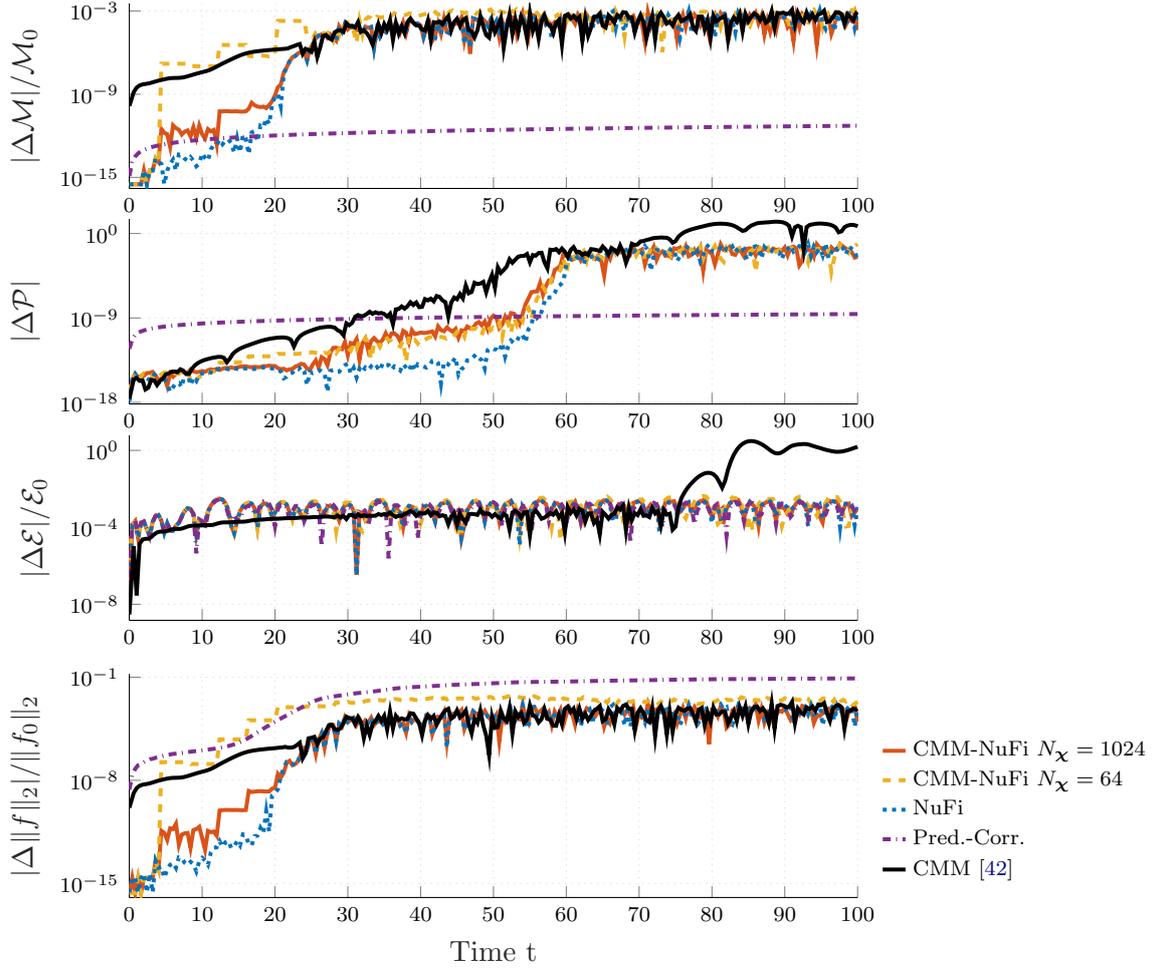}
    \caption{Conservation errors for the two-stream instability test case using a sample grid with $\Nsample=1024$. For the CMM method, we use adaptive remapping with incomp. threshold $\delta_\epsilon=10^{-2}$ (i.e., $\epsilon_\text{incomp}<\delta_\epsilon$) and CMM-NuFI is remapped after $\Nremap=20$ NuFI iterations. All simulations use the same time step size of $\tau=0.2$.}
    \label{fig:conservation_two_stream}
\end{figure}

To visualize the fine-scale properties, we have plotted the solution at $t=100$ for the different methods in \cref{fig:zoom_2stream}. Comparing the flow map methods with the classical semi-Lagrangian representation, we estimate a storage of 3.9MB for NuFI, storing 500 electric field vectors $E^k\in\mathbb{R}^{\Nsample}, \Nsample = 1024$, 0.8MB for CMM-NuFI storing 25 maps of size $64\times 64$, and the predictor-corrector scheme 8MB for the point values of $f$ on a $1024\times 1024$ sample grid. We see that less memory is needed in order to resolve the solution. 
Further, as shown in \cref{fig:zoom_2stream} the flow mapping approach allows us to zoom into the solution to access its fine-scale properties, which are lost for the classical predictor-corrector semi-Lagrangian scheme. Additionally, the initial range is not preserved by the predictor-corrector scheme, as a consequence of its diffusivity.
The conservation of mass, momentum and energy is shown for the three simulations in \cref{fig:conservation_two_stream}. We note that, in contrast to \cref{fig:conservation_landau_damping}, the simulation of the two-stream instability using the older CMM \cite{KrahYinBergmannNaveSchneider2024} shows less of a discrepancy in the conservation errors relative to the other methods than in the Landau damping case.

\section{Conclusion}
\label{sec:concl}
We have presented a hybrid flow mapping approach that merges the Characteristic Mapping Method (CMM) and Numerical Flow Iteration (NuFI) to compute the backward map for the Vlasov--Poisson equation. By integrating CMM's efficient backmapping with NuFI's accurate time-stepping, our method achieves enhanced conservation properties over CMM \cite{KrahYinBergmannNaveSchneider2024} and accelerates performance by orders of magnitude over NuFI \cite{KirchhartWilhelm2024}.
This enhanced performance enables high-resolution, long-duration simulations illustrated for different test-cases in 1D+1V.

Benchmarking against a semi-Lagrangian predictor-corrector scheme reveals improvements in the $L_2$ and maximum norms, as well as in momentum conservation. Importantly, the method uses the semigroup property of diffeomorphic flow maps to break down the flow into a sequence of submaps. This approach means memory needs grow linearly, but resolution improves exponentially compared to standard semi-Lagrangian schemes. We can thus "zoom in" to achieve extremely fine resolution, revealing small-scale details that are important for problems like plasma filamentation and for simulations in higher dimensions.
The fine-scale structures resolved in the two-stream instability distribution functions underscore this capability.

This hybrid framework opens several promising research avenues. Future work will focus on incorporating source terms, following the approach established for CMM \cite{YinKrahNaveSchneider2024}, to enable the simulation of collisional kinetic systems. Furthermore, we plan to generalize the method to handle boundary conditions, as demonstrated in NuFI, and to implement low-rank compression of the CMM maps \cite{WilhelmKormann2025} to ensure memory efficiency for high-dimensional plasma simulations.

%% file: figures/submaps.tex
\tikzset{every picture/.style={line width=0.75pt}} 

\begin{tikzpicture}[x=0.75pt,y=0.75pt,yscale=-1,xscale=1]

\draw [color={rgb, 255:red, 0; green, 38; blue, 255 }  ,draw opacity=1 ][line width=1.5]    (123.19,194.38) .. controls (204.51,171.58) and (215.05,64.99) .. (261.46,49.67) ;
\draw [shift={(193,138.93)}, rotate = 303.58] [fill={rgb, 255:red, 0; green, 38; blue, 255 }  ,fill opacity=1 ][line width=0.08]  [draw opacity=0] (15.6,-3.9) -- (0,0) -- (15.6,3.9) -- cycle    ;
\draw [shift={(118.13,195.67)}, rotate = 347.09] [fill={rgb, 255:red, 0; green, 38; blue, 255 }  ,fill opacity=1 ][line width=0.08]  [draw opacity=0] (11.61,-5.58) -- (0,0) -- (11.61,5.58) -- cycle    ;
\draw [line width=0.75]  [dash pattern={on 0.84pt off 2.51pt}]  (56,199) -- (286.33,199) ;
\draw [line width=1.5]    (46.33,199) -- (56,199) ;
\draw [line width=1.5]    (46.33,50) -- (56,50) ;
\draw [line width=0.75]  [dash pattern={on 4.5pt off 4.5pt}]  (210.5,71.17) .. controls (244.5,88.67) and (260.5,60.67) .. (285.5,75.5) ;
\draw [line width=0.75]  [dash pattern={on 0.84pt off 2.51pt}]  (56,50) -- (281,50) ;
\draw  [fill={rgb, 255:red, 155; green, 155; blue, 155 }  ,fill opacity=0.42 ][line width=1.5]  (210,13.5) -- (285,13.5) -- (285,86.5) -- (210,86.5) -- cycle ;
\draw [line width=0.75]  [dash pattern={on 4.5pt off 4.5pt}]  (209.5,51.17) .. controls (246.5,41.17) and (261.5,61.17) .. (285.5,51.5) ;
\draw [line width=0.75]  [dash pattern={on 4.5pt off 4.5pt}]  (210,28.5) .. controls (241.98,4.52) and (260,55.67) .. (286,33.5) ;
\draw [line width=0.75]  [dash pattern={on 4.5pt off 4.5pt}]  (221.5,13.67) .. controls (265,24.67) and (197.5,76.17) .. (225,86.17) ;
\draw [line width=0.75]  [dash pattern={on 4.5pt off 4.5pt}]  (241.5,13.67) .. controls (271,48.17) and (239.5,57.67) .. (245,86.17) ;
\draw [line width=0.75]  [dash pattern={on 4.5pt off 4.5pt}]  (262,13.67) .. controls (285.5,45.17) and (270,65.67) .. (265.5,86.17) ;
\draw  [color={rgb, 255:red, 189; green, 16; blue, 224 }  ,draw opacity=1 ] (258.63,45.67) -- (266.88,45.67) -- (266.88,53.67) -- (258.63,53.67) -- cycle ;
\draw  [fill={rgb, 255:red, 189; green, 16; blue, 224 }  ,fill opacity=1 ] (261.46,49.67) .. controls (261.46,48.75) and (262.2,48) .. (263.13,48) .. controls (264.05,48) and (264.79,48.75) .. (264.79,49.67) .. controls (264.79,50.59) and (264.05,51.33) .. (263.13,51.33) .. controls (262.2,51.33) and (261.46,50.59) .. (261.46,49.67) -- cycle ;
\draw [line width=1.5]    (51.3,242.9) -- (51.3,14.9) ;
\draw [shift={(51.3,10.9)}, rotate = 90] [fill={rgb, 255:red, 0; green, 0; blue, 0 }  ][line width=0.08]  [draw opacity=0] (11.61,-5.58) -- (0,0) -- (11.61,5.58) -- cycle    ;
\draw  [draw opacity=0][dash pattern={on 4.5pt off 4.5pt}] (74,161.9) -- (147.3,161.9) -- (147.3,233.9) -- (74,233.9) -- cycle ; \draw  [dash pattern={on 4.5pt off 4.5pt}] (75,161.9) -- (75,233.9)(93,161.9) -- (93,233.9)(111,161.9) -- (111,233.9)(129,161.9) -- (129,233.9)(147,161.9) -- (147,233.9) ; \draw  [dash pattern={on 4.5pt off 4.5pt}] (74,162.9) -- (147.3,162.9)(74,180.9) -- (147.3,180.9)(74,198.9) -- (147.3,198.9)(74,216.9) -- (147.3,216.9) ; \draw  [dash pattern={on 4.5pt off 4.5pt}]  ;
\draw  [fill={rgb, 255:red, 155; green, 155; blue, 155 }  ,fill opacity=0.49 ][line width=1.5]  (75,162.9) -- (147.5,162.9) -- (147.5,235.4) -- (75,235.4) -- cycle ;
\draw  [color={rgb, 255:red, 189; green, 16; blue, 224 }  ,draw opacity=1 ] (113.63,191.67) -- (121.88,191.67) -- (121.88,199.67) -- (113.63,199.67) -- cycle ;
\draw  [fill={rgb, 255:red, 189; green, 16; blue, 224 }  ,fill opacity=1 ] (116.46,195.67) .. controls (116.46,194.75) and (117.2,194) .. (118.13,194) .. controls (119.05,194) and (119.79,194.75) .. (119.79,195.67) .. controls (119.79,196.59) and (119.05,197.33) .. (118.13,197.33) .. controls (117.2,197.33) and (116.46,196.59) .. (116.46,195.67) -- cycle ;
\draw [line width=0.75]  [dash pattern={on 4.5pt off 4.5pt}]  (157.5,155.17) .. controls (191.5,172.67) and (207.5,144.67) .. (232.5,159.5) ;
\draw [line width=0.75]  [dash pattern={on 0.84pt off 2.51pt}]  (51.05,134) -- (285.55,134) ;
\draw  [fill={rgb, 255:red, 155; green, 155; blue, 155 }  ,fill opacity=0.42 ][line width=1.5]  (157,97.5) -- (232,97.5) -- (232,170.5) -- (157,170.5) -- cycle ;
\draw [line width=0.75]  [dash pattern={on 4.5pt off 4.5pt}]  (156.5,135.17) .. controls (193.5,125.17) and (208.5,145.17) .. (232.5,135.5) ;
\draw [line width=0.75]  [dash pattern={on 4.5pt off 4.5pt}]  (157,112.5) .. controls (188.98,88.52) and (207,139.67) .. (233,117.5) ;
\draw [line width=0.75]  [dash pattern={on 4.5pt off 4.5pt}]  (168.5,97.67) .. controls (191.05,124.03) and (163.55,147.03) .. (172,170.17) ;
\draw [line width=0.75]  [dash pattern={on 4.5pt off 4.5pt}]  (188.5,97.67) .. controls (195.55,122.53) and (186.5,141.67) .. (192,170.17) ;
\draw [line width=0.75]  [dash pattern={on 4.5pt off 4.5pt}]  (209,97.67) .. controls (205.05,146.53) and (217,149.67) .. (212.5,170.17) ;
\draw  [color={rgb, 255:red, 189; green, 16; blue, 224 }  ,draw opacity=1 ] (191.63,129.67) -- (199.88,129.67) -- (199.88,137.67) -- (191.63,137.67) -- cycle ;
\draw  [fill={rgb, 255:red, 189; green, 16; blue, 224 }  ,fill opacity=1 ] (194.46,133.67) .. controls (194.46,132.75) and (195.2,132) .. (196.13,132) .. controls (197.05,132) and (197.79,132.75) .. (197.79,133.67) .. controls (197.79,134.59) and (197.05,135.33) .. (196.13,135.33) .. controls (195.2,135.33) and (194.46,134.59) .. (194.46,133.67) -- cycle ;
\draw [line width=1.5]    (46.33,134) -- (56,134) ;

\draw (32,39.4) node [anchor=north west][inner sep=0.75pt]    {$t$};
\draw (24,192.4) node [anchor=north west][inner sep=0.75pt]    {$t_{0}$};
\draw (16,3) node [anchor=north west][inner sep=0.75pt]   [align=left] {time};
\draw (59,194.4) node [anchor=north west][inner sep=0.75pt]    {$v$};
\draw (94,235.4) node [anchor=north west][inner sep=0.75pt]    {$x$};
\draw (154,196.4) node [anchor=north west][inner sep=0.75pt]    {$\boldsymbol{\Phi }_{t_{0}}^{t_{0}}( x,v) =( x,v)$};
\draw (147,24.4) node [anchor=north west][inner sep=0.75pt]    {$\boldsymbol{\Phi }_{t_{1}}^{t}( x,v)$};
\draw (234,87.43) node [anchor=north west][inner sep=0.75pt]    {$x$};
\draw (191,47.4) node [anchor=north west][inner sep=0.75pt]    {$v$};
\draw (94,106.4) node [anchor=north west][inner sep=0.75pt]    {$\boldsymbol{\Phi }_{t_{0}}^{t_{1}}( x,v)$};
\draw (28,123.4) node [anchor=north west][inner sep=0.75pt]    {$t_{1}$};

\end{tikzpicture}

%% file: figures/NuFi-algo.tex
\tikzset{every picture/.style={line width=0.75pt}} 

\begin{tikzpicture}[x=0.75pt,y=0.75pt,yscale=-1,xscale=1]

\draw  [fill={rgb, 255:red, 155; green, 155; blue, 155 }  ,fill opacity=1 ] (64,43.85) -- (90.66,43.85) -- (90.66,39.8) -- (108.43,47.9) -- (90.66,56) -- (90.66,51.95) -- (64,51.95) -- cycle ;
\draw  [fill={rgb, 255:red, 155; green, 155; blue, 155 }  ,fill opacity=1 ] (144,43.85) -- (170.66,43.85) -- (170.66,39.8) -- (188.43,47.9) -- (170.66,56) -- (170.66,51.95) -- (144,51.95) -- cycle ;
\draw  [fill={rgb, 255:red, 208; green, 2; blue, 27 }  ,fill opacity=1 ] (219.94,41.8) -- (234.34,41.78) .. controls (246.27,41.77) and (255.95,51.43) .. (255.96,63.35) -- (256.02,108.98) -- (260.02,108.98) -- (250.04,128.99) -- (240.02,109) -- (244.02,109) -- (243.96,63.37) .. controls (243.96,58.07) and (239.65,53.77) .. (234.35,53.78) -- (219.95,53.8) -- cycle ;
\draw  [fill={rgb, 255:red, 155; green, 155; blue, 155 }  ,fill opacity=1 ] (254.97,166) -- (254.97,181.9) .. controls (254.97,195.07) and (244.29,205.75) .. (231.12,205.75) -- (207.87,205.75) -- (207.87,211.8) -- (184.97,200.35) -- (207.87,188.9) -- (207.87,194.96) -- (231.12,194.96) .. controls (238.33,194.96) and (244.18,189.11) .. (244.18,181.9) -- (244.18,166) -- cycle ;
\draw  [fill={rgb, 255:red, 155; green, 155; blue, 155 }  ,fill opacity=1 ] (147.66,208.33) -- (127.28,207.98) .. controls (119.03,207.83) and (112.46,201.03) .. (112.6,192.78) -- (112.77,182.82) -- (108.83,182.75) -- (119,163.29) -- (128.49,183.09) -- (124.55,183.03) -- (124.38,192.98) .. controls (124.35,194.73) and (125.74,196.17) .. (127.48,196.2) -- (147.87,196.55) -- cycle ;
\draw  [color={rgb, 255:red, 0; green, 0; blue, 0 }  ,draw opacity=1 ][fill={rgb, 255:red, 208; green, 2; blue, 27 }  ,fill opacity=1 ] (112.43,134.8) -- (112.43,103.98) .. controls (112.43,97.92) and (117.35,93) .. (123.41,93) -- (220.65,93) .. controls (226.72,93) and (231.63,97.92) .. (231.63,103.98) -- (231.63,108.8) -- (235.43,108.8) -- (224.98,123.93) -- (214.53,108.8) -- (218.33,108.8) -- (218.33,106.3) .. controls (218.33,106.3) and (218.33,106.3) .. (218.33,106.3) -- (125.73,106.3) .. controls (125.73,106.3) and (125.73,106.3) .. (125.73,106.3) -- (125.73,134.8) -- cycle ;
\draw  [fill={rgb, 255:red, 155; green, 155; blue, 155 }  ,fill opacity=1 ] (260.43,165.8) -- (260.43,188.58) .. controls (260.43,197.58) and (267.73,204.88) .. (276.73,204.88) -- (290.73,204.88) -- (290.73,209.66) -- (312.25,199.23) -- (290.73,188.8) -- (290.73,193.58) -- (276.73,193.58) .. controls (273.97,193.58) and (271.73,191.34) .. (271.73,188.58) -- (271.73,165.8) -- cycle ;

\draw (40,37.4) node [anchor=north west][inner sep=0.75pt]    {$f_{0}$};
\draw (115,37.4) node [anchor=north west][inner sep=0.75pt]    {$\rho ^{0}$};
\draw (195,39.4) node [anchor=north west][inner sep=0.75pt]    {$E^{0}$};
\draw (207,138.4) node [anchor=north west][inner sep=0.75pt]   {$f^{k+1} =f_{0}\left( \Psi _{0}^{k\tau }\right)$};
\draw (315,190.4) node [anchor=north west][inner sep=0.75pt]    {$f^{n}$};
\draw (71,139.4) node [anchor=north west][inner sep=0.75pt]    {$\left[ E^{0} ,\cdots ,E^{k}\right]$};
\draw (157,190.4) node [anchor=north west][inner sep=0.75pt]    {$\rho ^{k}$};
\draw (160,92) node [anchor=north west][inner sep=0.75pt]  [font=\small] [align=left] {NuFi};
\draw (260,82) node [anchor=north west][inner sep=0.75pt]  [font=\small] [align=left] {NuFi};
\draw (139,19) node [anchor=north west][inner sep=0.75pt]  [font=\small] [align=left] {Poisson};
\draw (50,188) node [anchor=north west][inner sep=0.75pt]   [align=left] {Poisson};
\draw (97.73,69.7) node [anchor=north west][inner sep=0.75pt] {$\Psi _{0}^{k\tau } =\Psi _{0}^{\tau } \circ \cdots \circ \Psi _{( k-1) \tau }^{k\tau }$};

\end{tikzpicture}

%% file: figures/volumepreserving_upsample.tex
%
\definecolor{mycolor1}{rgb}{1.00000,0.00000,1.00000}%
\definecolor{mycolor2}{rgb}{0.00000,1.00000,1.00000}%
\begin{tikzpicture}

\begin{axis}[%
width=0.976\figurewidth,
height=\figureheight,
at={(0\figurewidth,0\figureheight)},
scale only axis,
xmin=0.0000000000e+00,
xmax=2.0000000000e+01,
xminorticks=true,
xlabel style={font=\color{white!15!black}},
xlabel={Number of Composed Maps},
ymode=log,
ymin=2.5657254099e-13,
ymax=2.1744828160e-12,
yminorticks=true,
ylabel style={font=\color{white!15!black}},
ylabel={Max Incompressibility Error},
axis background/.style={fill=white},
title style={font=\bfseries},
xmajorgrids,
ymajorgrids,
yminorgrids,
legend style={at={(0.03,0.97)}, anchor=north west, legend cell align=left, align=left, draw=white!15!black},
scaled ticks=false, xticklabel style={/pgf/number format/fixed},yticklabel style={/pgf/number format/fixed}
]
\addplot [color=blue, line width=1.5pt, mark size=1.0pt, mark=o, mark options={solid, blue}]
  table[row sep=crcr]{%
1.0000000000e+00	2.5657254099e-13\\
2.0000000000e+00	5.2169379927e-13\\
3.0000000000e+00	6.9644290335e-13\\
4.0000000000e+00	9.4069196877e-13\\
5.0000000000e+00	1.1017853296e-12\\
6.0000000000e+00	1.1758372054e-12\\
7.0000000000e+00	1.4195311593e-12\\
8.0000000000e+00	1.2501111257e-12\\
9.0000000000e+00	1.2668754934e-12\\
1.0000000000e+01	1.1510792319e-12\\
1.1000000000e+01	1.2232437285e-12\\
1.2000000000e+01	1.2796430582e-12\\
1.3000000000e+01	1.3360423878e-12\\
1.4000000000e+01	1.3531398224e-12\\
1.5000000000e+01	1.6413537196e-12\\
1.6000000000e+01	1.6102674749e-12\\
1.7000000000e+01	1.8001156121e-12\\
1.8000000000e+01	1.7352785875e-12\\
1.9000000000e+01	2.0303758674e-12\\
2.0000000000e+01	2.1744828160e-12\\
};
\addlegendentry{1024x1024}

\addplot [color=red, line width=1.5pt, mark size=1pt, mark=square, mark options={solid, red}]
  table[row sep=crcr]{%
1.0000000000e+00	3.6703973194e-13\\
2.0000000000e+00	3.1885605267e-13\\
3.0000000000e+00	3.7381209239e-13\\
4.0000000000e+00	3.8247183198e-13\\
5.0000000000e+00	3.6859404418e-13\\
6.0000000000e+00	3.9923619966e-13\\
7.0000000000e+00	4.1955328101e-13\\
8.0000000000e+00	4.0012437807e-13\\
9.0000000000e+00	4.1033842990e-13\\
1.0000000000e+01	4.0767389464e-13\\
1.1000000000e+01	4.4098058538e-13\\
1.2000000000e+01	4.0545344859e-13\\
1.3000000000e+01	4.6607162574e-13\\
1.4000000000e+01	4.9449333517e-13\\
1.5000000000e+01	5.0182080713e-13\\
1.6000000000e+01	4.9227288912e-13\\
1.7000000000e+01	4.9316106754e-13\\
1.8000000000e+01	4.9560355819e-13\\
1.9000000000e+01	5.7032156775e-13\\
2.0000000000e+01	5.1447734961e-13\\
};
\addlegendentry{32x32}

\addplot [color=green, line width=1.5pt, mark size=1pt, mark=triangle, mark options={solid, green}]
  table[row sep=crcr]{%
1.0000000000e+00	3.6060043840e-13\\
2.0000000000e+00	3.4194869158e-13\\
3.0000000000e+00	3.3173463976e-13\\
4.0000000000e+00	3.8280489889e-13\\
5.0000000000e+00	3.5860203695e-13\\
6.0000000000e+00	3.5060843118e-13\\
7.0000000000e+00	4.3998138466e-13\\
8.0000000000e+00	4.0101255649e-13\\
9.0000000000e+00	4.4075854078e-13\\
1.0000000000e+01	3.9734882051e-13\\
1.1000000000e+01	4.3809400552e-13\\
1.2000000000e+01	5.0071058411e-13\\
1.3000000000e+01	4.7606363296e-13\\
1.4000000000e+01	5.0515147620e-13\\
1.5000000000e+01	4.9271697833e-13\\
1.6000000000e+01	5.1980642013e-13\\
1.7000000000e+01	5.0714987765e-13\\
1.8000000000e+01	5.4578563891e-13\\
1.9000000000e+01	5.4389825976e-13\\
2.0000000000e+01	6.1728400169e-13\\
};
\addlegendentry{64x64}

\addplot [color=mycolor1, line width=1.5pt, mark size=1pt, mark=diamond, mark options={solid, mycolor1}]
  table[row sep=crcr]{%
1.0000000000e+00	3.4616753908e-13\\
2.0000000000e+00	4.3054448895e-13\\
3.0000000000e+00	3.7059244562e-13\\
4.0000000000e+00	4.2499337383e-13\\
5.0000000000e+00	5.2902127123e-13\\
6.0000000000e+00	4.4186876380e-13\\
7.0000000000e+00	4.2321701699e-13\\
8.0000000000e+00	4.5541348470e-13\\
9.0000000000e+00	4.1033842990e-13\\
1.0000000000e+01	4.3853809473e-13\\
1.1000000000e+01	5.4967141949e-13\\
1.2000000000e+01	5.0426329778e-13\\
1.3000000000e+01	5.2169379927e-13\\
1.4000000000e+01	5.1181281435e-13\\
1.5000000000e+01	5.2491344604e-13\\
1.6000000000e+01	4.9649173661e-13\\
1.7000000000e+01	5.2646775828e-13\\
1.8000000000e+01	5.8641980161e-13\\
1.9000000000e+01	5.8464344477e-13\\
2.0000000000e+01	5.9829918797e-13\\
};
\addlegendentry{128x128}

\addplot [color=mycolor2, line width=1.5pt, mark size=1pt, mark=triangle, mark options={solid, rotate=180, mycolor2}]
  table[row sep=crcr]{%
1.0000000000e+00	3.4994229736e-13\\
2.0000000000e+00	4.5297099405e-13\\
3.0000000000e+00	4.1444625509e-13\\
4.0000000000e+00	4.6040948831e-13\\
5.0000000000e+00	3.8002934133e-13\\
6.0000000000e+00	4.2199577166e-13\\
7.0000000000e+00	4.4164671920e-13\\
8.0000000000e+00	4.8872017544e-13\\
9.0000000000e+00	5.6843418861e-13\\
1.0000000000e+01	6.2949645496e-13\\
1.1000000000e+01	6.0518257072e-13\\
1.2000000000e+01	6.1206595348e-13\\
1.3000000000e+01	6.8955952059e-13\\
1.4000000000e+01	7.0921046813e-13\\
1.5000000000e+01	6.9788619328e-13\\
1.6000000000e+01	7.4751316248e-13\\
1.7000000000e+01	7.5595085747e-13\\
1.8000000000e+01	7.8914652590e-13\\
1.9000000000e+01	7.7671202803e-13\\
2.0000000000e+01	7.3852035598e-13\\
};
\addlegendentry{256x256}

\addplot [color=black, line width=1.5pt, mark size=1.0pt, mark=star, mark options={solid, black}]
  table[row sep=crcr]{%
1.0000000000e+00	4.2443826231e-13\\
2.0000000000e+00	4.4908521346e-13\\
3.0000000000e+00	4.2665870836e-13\\
4.0000000000e+00	4.6962433942e-13\\
5.0000000000e+00	5.2557957986e-13\\
6.0000000000e+00	5.4478643818e-13\\
7.0000000000e+00	5.9841021027e-13\\
8.0000000000e+00	6.5325522769e-13\\
9.0000000000e+00	6.8889338678e-13\\
1.0000000000e+01	6.9322325658e-13\\
1.1000000000e+01	8.4299234260e-13\\
1.2000000000e+01	8.5387252824e-13\\
1.3000000000e+01	9.4990681987e-13\\
1.4000000000e+01	1.0287326546e-12\\
1.5000000000e+01	1.0575984533e-12\\
1.6000000000e+01	9.7721830627e-13\\
1.7000000000e+01	1.0170753129e-12\\
1.8000000000e+01	1.0080825064e-12\\
1.9000000000e+01	1.1555201240e-12\\
2.0000000000e+01	1.1947109968e-12\\
};
\addlegendentry{512x512}

\end{axis}
\end{tikzpicture}%

%% file: figures/volumepreserving_l2error.tex
%
\definecolor{mycolor1}{rgb}{0.00000,0.44700,0.74100}%
\definecolor{mycolor2}{rgb}{0.85000,0.32500,0.09800}%
\begin{tikzpicture}

\begin{axis}[%
width=0.976\figurewidth,
height=\figureheight,
at={(0\figurewidth,0\figureheight)},
scale only axis,
xmode=log,
xmin=30.0000000000,
xmax=542.0000000000,
xminorticks=true,
xlabel style={font=\color{white!15!black}},
xlabel={$1/\Nmapgrid$},
ymode=log,
ymin=0.0001000000,
ymax=10000.0000000000,
yminorticks=true,
ylabel style={font=\color{white!15!black}},
ylabel={$\norm{(\flowmap-\mapcmm)}_\infty$},
axis background/.style={fill=white},
xmajorgrids,
xminorgrids,
ymajorgrids,
yminorgrids,
legend style={legend cell align=left, align=left, draw=white!15!black},
scaled ticks=false, xticklabel style={/pgf/number format/fixed},yticklabel style={/pgf/number format/fixed}
]
\addplot [color=mycolor1, line width=1.5pt, only marks, mark size=0.6pt, mark=o, mark options={solid, mycolor1}]
  table[row sep=crcr]{%
32.0000000000	23.5620411474\\
64.0000000000	2.0915273229\\
128.0000000000	0.1426506038\\
256.0000000000	0.0091050704\\
512.0000000000	0.0005633276\\
};
\addlegendentry{1 maps}

\addplot [color=mycolor2, line width=1.5pt, only marks, mark size=0.6, mark=square, mark options={solid, mycolor2}]
  table[row sep=crcr]{%
32.0000000000	471.2408229475\\
64.0000000000	41.8305464587\\
128.0000000000	2.8530120757\\
256.0000000000	0.1821014084\\
512.0000000000	0.0112665511\\
};
\addlegendentry{20 maps}

\addplot [color=black, dashed, line width=1.5pt]
  table[row sep=crcr]{%
32.0000000000	95.3674316406\\
64.0000000000	5.9604644775\\
128.0000000000	0.3725290298\\
256.0000000000	0.0232830644\\
512.0000000000	0.0014551915\\
1024.0000000000	0.0000909495\\
};
\addlegendentry{$\mathcal{O}(\Nmapgrid^{-4})$}

\end{axis}
\end{tikzpicture}%

%% file: figures/hybrid-scheme.tex
\tikzset{every picture/.style={line width=0.75pt}} 

\begin{tikzpicture}[x=0.75pt,y=0.75pt,yscale=-1,xscale=1]

\draw    (60.1,89.7) -- (429.1,89.7) ;
\draw [shift={(431.1,89.7)}, rotate = 180] [color={rgb, 255:red, 0; green, 0; blue, 0 }  ][line width=0.75]    (10.93,-3.29) .. controls (6.95,-1.4) and (3.31,-0.3) .. (0,0) .. controls (3.31,0.3) and (6.95,1.4) .. (10.93,3.29)   ;
\draw [color={rgb, 255:red, 0; green, 0; blue, 255 }  ,draw opacity=0.22 ]   (60.21,84.41) .. controls (62.25,58.8) and (92.58,50.51) .. (115.36,56.72) .. controls (129.34,60.53) and (140.48,69.81) .. (140.48,83.9) ;
\draw [shift={(60.1,87.7)}, rotate = 269.21] [fill={rgb, 255:red, 0; green, 0; blue, 255 }  ,fill opacity=0.22 ][line width=0.08]  [draw opacity=0] (8.93,-4.29) -- (0,0) -- (8.93,4.29) -- cycle    ;
\draw    (60.1,89.7) -- (411.77,89.7) (80.1,85.7) -- (80.1,93.7)(100.1,85.7) -- (100.1,93.7)(120.1,85.7) -- (120.1,93.7)(140.1,85.7) -- (140.1,93.7)(160.1,85.7) -- (160.1,93.7)(180.1,85.7) -- (180.1,93.7)(200.1,85.7) -- (200.1,93.7)(220.1,85.7) -- (220.1,93.7)(240.1,85.7) -- (240.1,93.7)(260.1,85.7) -- (260.1,93.7)(280.1,85.7) -- (280.1,93.7)(300.1,85.7) -- (300.1,93.7)(320.1,85.7) -- (320.1,93.7)(340.1,85.7) -- (340.1,93.7)(360.1,85.7) -- (360.1,93.7)(380.1,85.7) -- (380.1,93.7)(400.1,85.7) -- (400.1,93.7) ;
\draw [shift={(60.1,89.7)}, rotate = 0] [color={rgb, 255:red, 0; green, 0; blue, 0 }  ][fill={rgb, 255:red, 0; green, 0; blue, 0 }  ][line width=0.75]      (0, 0) circle [x radius= 1.34, y radius= 1.34]   ;
\draw [color={rgb, 255:red, 0; green, 0; blue, 255 }  ,draw opacity=0.55 ]   (140.19,84.6) .. controls (141.61,65.44) and (158.8,55.92) .. (176.91,54.85) .. controls (198.06,53.6) and (220.48,63.91) .. (220.48,83.9) ;
\draw [shift={(140.1,87.7)}, rotate = 269.21] [fill={rgb, 255:red, 0; green, 0; blue, 255 }  ,fill opacity=0.55 ][line width=0.08]  [draw opacity=0] (8.93,-4.29) -- (0,0) -- (8.93,4.29) -- cycle    ;
\draw [color={rgb, 255:red, 0; green, 0; blue, 255 }  ,draw opacity=1 ]   (220.2,84.43) .. controls (223.48,43.07) and (300.48,46.85) .. (300.48,83.9) ;
\draw [shift={(220.1,87.7)}, rotate = 269.21] [fill={rgb, 255:red, 0; green, 0; blue, 255 }  ,fill opacity=1 ][line width=0.08]  [draw opacity=0] (8.93,-4.29) -- (0,0) -- (8.93,4.29) -- cycle    ;
\draw [color={rgb, 255:red, 0; green, 100; blue, 0 }  ,draw opacity=1 ]   (300.57,98.75) .. controls (302.59,117.39) and (320.38,115.53) .. (320.38,94.63) ;
\draw [shift={(300.38,95.63)}, rotate = 89.26] [fill={rgb, 255:red, 0; green, 100; blue, 0 }  ,fill opacity=1 ][line width=0.08]  [draw opacity=0] (8.93,-4.29) -- (0,0) -- (8.93,4.29) -- cycle    ;
\draw [color={rgb, 255:red, 0; green, 100; blue, 0 }  ,draw opacity=1 ]   (320.57,98.75) .. controls (322.59,117.39) and (340.38,115.53) .. (340.38,94.63) ;
\draw [shift={(320.38,95.63)}, rotate = 89.26] [fill={rgb, 255:red, 0; green, 100; blue, 0 }  ,fill opacity=1 ][line width=0.08]  [draw opacity=0] (8.93,-4.29) -- (0,0) -- (8.93,4.29) -- cycle    ;
\draw [color={rgb, 255:red, 0; green, 100; blue, 0 }  ,draw opacity=1 ]   (340.57,98.75) .. controls (342.59,117.39) and (360.38,115.53) .. (360.38,94.63) ;
\draw [shift={(340.38,95.63)}, rotate = 89.26] [fill={rgb, 255:red, 0; green, 100; blue, 0 }  ,fill opacity=1 ][line width=0.08]  [draw opacity=0] (8.93,-4.29) -- (0,0) -- (8.93,4.29) -- cycle    ;
\draw [color={rgb, 255:red, 0; green, 100; blue, 0 }  ,draw opacity=1 ]   (360.57,98.75) .. controls (362.59,117.39) and (380.38,115.53) .. (380.38,94.63) ;
\draw [shift={(360.38,95.63)}, rotate = 89.26] [fill={rgb, 255:red, 0; green, 100; blue, 0 }  ,fill opacity=1 ][line width=0.08]  [draw opacity=0] (8.93,-4.29) -- (0,0) -- (8.93,4.29) -- cycle    ;
\draw    (60.1,169.7) -- (429.1,169.7) ;
\draw [shift={(431.1,169.7)}, rotate = 180] [color={rgb, 255:red, 0; green, 0; blue, 0 }  ][line width=0.75]    (10.93,-3.29) .. controls (6.95,-1.4) and (3.31,-0.3) .. (0,0) .. controls (3.31,0.3) and (6.95,1.4) .. (10.93,3.29)   ;
\draw [color={rgb, 255:red, 0; green, 0; blue, 255 }  ,draw opacity=0.22 ]   (60.21,164.41) .. controls (62.25,138.8) and (92.58,130.51) .. (115.36,136.72) .. controls (129.34,140.53) and (140.48,149.81) .. (140.48,163.9) ;
\draw [shift={(60.1,167.7)}, rotate = 269.21] [fill={rgb, 255:red, 0; green, 0; blue, 255 }  ,fill opacity=0.22 ][line width=0.08]  [draw opacity=0] (8.93,-4.29) -- (0,0) -- (8.93,4.29) -- cycle    ;
\draw    (60.1,169.7) -- (411.77,169.7) (80.1,165.7) -- (80.1,173.7)(100.1,165.7) -- (100.1,173.7)(120.1,165.7) -- (120.1,173.7)(140.1,165.7) -- (140.1,173.7)(160.1,165.7) -- (160.1,173.7)(180.1,165.7) -- (180.1,173.7)(200.1,165.7) -- (200.1,173.7)(220.1,165.7) -- (220.1,173.7)(240.1,165.7) -- (240.1,173.7)(260.1,165.7) -- (260.1,173.7)(280.1,165.7) -- (280.1,173.7)(300.1,165.7) -- (300.1,173.7)(320.1,165.7) -- (320.1,173.7)(340.1,165.7) -- (340.1,173.7)(360.1,165.7) -- (360.1,173.7)(380.1,165.7) -- (380.1,173.7)(400.1,165.7) -- (400.1,173.7) ;
\draw [shift={(60.1,169.7)}, rotate = 0] [color={rgb, 255:red, 0; green, 0; blue, 0 }  ][fill={rgb, 255:red, 0; green, 0; blue, 0 }  ][line width=0.75]      (0, 0) circle [x radius= 1.34, y radius= 1.34]   ;
\draw [color={rgb, 255:red, 0; green, 0; blue, 255 }  ,draw opacity=0.55 ]   (140.19,164.6) .. controls (141.61,145.44) and (158.8,135.92) .. (176.91,134.85) .. controls (198.06,133.6) and (220.48,143.91) .. (220.48,163.9) ;
\draw [shift={(140.1,167.7)}, rotate = 269.21] [fill={rgb, 255:red, 0; green, 0; blue, 255 }  ,fill opacity=0.55 ][line width=0.08]  [draw opacity=0] (8.93,-4.29) -- (0,0) -- (8.93,4.29) -- cycle    ;
\draw [color={rgb, 255:red, 0; green, 0; blue, 255 }  ,draw opacity=1 ]   (220.2,164.43) .. controls (223.48,123.07) and (300.48,126.85) .. (300.48,163.9) ;
\draw [shift={(220.1,167.7)}, rotate = 269.21] [fill={rgb, 255:red, 0; green, 0; blue, 255 }  ,fill opacity=1 ][line width=0.08]  [draw opacity=0] (8.93,-4.29) -- (0,0) -- (8.93,4.29) -- cycle    ;
\draw [color={rgb, 255:red, 0; green, 0; blue, 161 }  ,draw opacity=1 ]   (300.21,164.4) .. controls (301.84,144.04) and (321.34,134.63) .. (340.71,134.75) .. controls (360.67,134.87) and (380.48,145.11) .. (380.48,163.9) ;
\draw [shift={(300.1,167.7)}, rotate = 269.21] [fill={rgb, 255:red, 0; green, 0; blue, 161 }  ,fill opacity=1 ][line width=0.08]  [draw opacity=0] (8.93,-4.29) -- (0,0) -- (8.93,4.29) -- cycle    ;
\draw [color={rgb, 255:red, 0; green, 100; blue, 0 }  ,draw opacity=1 ]   (380.57,178.75) .. controls (382.59,197.39) and (400.38,195.53) .. (400.38,174.63) ;
\draw [shift={(380.38,175.63)}, rotate = 89.26] [fill={rgb, 255:red, 0; green, 100; blue, 0 }  ,fill opacity=1 ][line width=0.08]  [draw opacity=0] (8.93,-4.29) -- (0,0) -- (8.93,4.29) -- cycle    ;
\draw  [color={rgb, 255:red, 208; green, 2; blue, 27 }  ,draw opacity=1 ][fill={rgb, 255:red, 181; green, 0; blue, 0 }  ,fill opacity=1 ] (377.34,89.69) .. controls (377.34,88.22) and (378.54,87.02) .. (380.01,87.02) .. controls (381.48,87.02) and (382.68,88.22) .. (382.68,89.69) .. controls (382.68,91.16) and (381.48,92.36) .. (380.01,92.36) .. controls (378.54,92.36) and (377.34,91.16) .. (377.34,89.69) -- cycle ;
\draw  [color={rgb, 255:red, 208; green, 2; blue, 27 }  ,draw opacity=1 ][fill={rgb, 255:red, 181; green, 0; blue, 0 }  ,fill opacity=1 ] (57.43,89.7) .. controls (57.43,88.23) and (58.63,87.03) .. (60.1,87.03) .. controls (61.57,87.03) and (62.77,88.23) .. (62.77,89.7) .. controls (62.77,91.17) and (61.57,92.37) .. (60.1,92.37) .. controls (58.63,92.37) and (57.43,91.17) .. (57.43,89.7) -- cycle ;
\draw  [color={rgb, 255:red, 208; green, 2; blue, 27 }  ,draw opacity=1 ][fill={rgb, 255:red, 181; green, 0; blue, 0 }  ,fill opacity=1 ] (397.34,169.49) .. controls (397.34,168.02) and (398.54,166.82) .. (400.01,166.82) .. controls (401.48,166.82) and (402.68,168.02) .. (402.68,169.49) .. controls (402.68,170.96) and (401.48,172.16) .. (400.01,172.16) .. controls (398.54,172.16) and (397.34,170.96) .. (397.34,169.49) -- cycle ;
\draw  [color={rgb, 255:red, 208; green, 2; blue, 27 }  ,draw opacity=1 ][fill={rgb, 255:red, 181; green, 0; blue, 0 }  ,fill opacity=1 ] (57.43,170.37) .. controls (57.43,168.89) and (58.63,167.7) .. (60.1,167.7) .. controls (61.57,167.7) and (62.77,168.89) .. (62.77,170.37) .. controls (62.77,171.84) and (61.57,173.03) .. (60.1,173.03) .. controls (58.63,173.03) and (57.43,171.84) .. (57.43,170.37) -- cycle ;
\draw [color={rgb, 255:red, 208; green, 2; blue, 27 }  ,draw opacity=1 ] [dash pattern={on 4.5pt off 4.5pt}]  (60.23,84.5) .. controls (67.63,-19.97) and (380.61,-2.06) .. (380.01,89.69) ;
\draw [shift={(60.1,87.7)}, rotate = 270.62] [fill={rgb, 255:red, 208; green, 2; blue, 27 }  ,fill opacity=1 ][line width=0.08]  [draw opacity=0] (8.93,-4.29) -- (0,0) -- (8.93,4.29) -- cycle    ;

\draw (16,79.4) node [anchor=north west][inner sep=0.75pt]  [font=\small]  {$N$};
\draw (5,52) node [anchor=north west][inner sep=0.75pt]   [align=left] {iteration};
\draw (433.1,93.1) node [anchor=north west][inner sep=0.75pt]    {$t$};
\draw (305,70.4) node [anchor=north west][inner sep=0.75pt]    {$\tau $};
\draw (86,59) node [anchor=north west][inner sep=0.75pt]   [align=left] {$\displaystyle \mapcmm^{t_{1}}_{0}$};
\draw (169,59) node [anchor=north west][inner sep=0.75pt]   [align=left] {$\displaystyle \mapcmm^{t_{2}}_{t_{1}}$};
\draw (249,59) node [anchor=north west][inner sep=0.75pt]   [align=left] {$\displaystyle \mapcmm^{t_{3}}_{t_{2}}$};
\draw (299.63,113.33) node [anchor=north west][inner sep=0.75pt]  [font=\scriptsize] [align=left] {$\displaystyle \NuFiMap^{t_{3} +\tau }_{t_{3}} \ \cdots \ \ \ \NuFiMap^{t_{3} +N\tau }_{t_{3} +( N-1) \tau }$};
\draw (13,159.4) node [anchor=north west][inner sep=0.75pt]  [font=\small]  {$N+1$};
\draw (5,132) node [anchor=north west][inner sep=0.75pt]   [align=left] {iteration};
\draw (433.1,173.1) node [anchor=north west][inner sep=0.75pt]    {$t$};
\draw (385,149.4) node [anchor=north west][inner sep=0.75pt]    {$\tau $};
\draw (86,139) node [anchor=north west][inner sep=0.75pt]   [align=left] {$\displaystyle \mapcmm^{t_{1}}_{0}$};
\draw (169,139) node [anchor=north west][inner sep=0.75pt]   [align=left] {$\displaystyle \mapcmm^{t_{2}}_{t_{1}}$};
\draw (249,139) node [anchor=north west][inner sep=0.75pt]   [align=left] {$\displaystyle \mapcmm^{t_{3}}_{t_{2}}$};
\draw (374.63,193.33) node [anchor=north west][inner sep=0.75pt]  [font=\scriptsize] [align=left] {$\displaystyle \NuFiMap^{t_{4} +\tau }_{t_{4}}$};
\draw (329,139) node [anchor=north west][inner sep=0.75pt]   [align=left] {$\displaystyle \mapcmm^{t_{4}}_{t_{3}}$};
\draw (381,70.4) node [anchor=north west][inner sep=0.75pt]  [font=\footnotesize,color={rgb, 255:red, 208; green, 2; blue, 27 }  ,opacity=1 ]  {$( x,v)$};
\draw (400,147.4) node [anchor=north west][inner sep=0.75pt]  [font=\footnotesize,color={rgb, 255:red, 208; green, 2; blue, 27 }  ,opacity=1 ]  {$( x,v)$};
\draw (38,96.73) node [anchor=north west][inner sep=0.75pt]  [font=\footnotesize,color={rgb, 255:red, 208; green, 2; blue, 27 }  ,opacity=1 ]  {$\left( x^{*} ,v^{*}\right)$};
\draw (38,176.07) node [anchor=north west][inner sep=0.75pt]  [font=\footnotesize,color={rgb, 255:red, 208; green, 2; blue, 27 }  ,opacity=1 ]  {$\left( x^{*} ,v^{*}\right)$};
\draw (197.33,14.33) node [anchor=north west][inner sep=0.75pt]  [color={rgb, 255:red, 208; green, 2; blue, 27 }  ,opacity=1 ] [align=left] {$\displaystyle \flowmap^{t}_{0}$};

\end{tikzpicture}

%% file: figures/convergence_comparison_landau_damping.tikz
%
\definecolor{mycolor1}{rgb}{0.00000,0.44700,0.74100}%
\definecolor{mycolor2}{rgb}{0.85000,0.32500,0.09800}%
\definecolor{mycolor3}{rgb}{0.92900,0.69400,0.12500}%
\begin{tikzpicture}

\begin{axis}[%
width=0.951\figurewidth,
height=\figureheight,
at={(0\figurewidth,0\figureheight)},
scale only axis,
xmode=log,
xmin=0.002734375,
xmax=0.046875,
xtick={0.00390625,0.0078125,0.015625,0.03125,0.0625},
xticklabels={{1/256},{1/128},{1/64},{1/32},{1/16}},
xminorticks=true,
xlabel style={font=\color{white!15!black}},
xlabel={$1/\Nmapgrid$},
ymode=log,
ymin=2.00410188444522e-11,
ymax=3.91988539327759e-06,
yminorticks=true,
ylabel style={font=\color{white!15!black}},
ylabel={$||f - f^{\mathrm{NuFi}}||_\infty$},
axis background/.style={fill=white},
xmajorgrids=true,
xminorgrids=true,
ymajorgrids=true,
yminorgrids=true,
grid style={dashed, draw=gray!40},
legend style={fill=none,
    draw=none,
    at={(0.01,0.99)},
    anchor=north west,
    legend cell align=left,
    align=left,
    /tikz/every mark/.append style={scale=0.5} 
}
]

\node[anchor=south east, align=left] at (rel axis cs:0.95, 0.01) {\small$||f - f^{\mathrm{NuFi}}||_\infty = \mathcal{O}(1/\Nmapgrid^\alpha)$};

\addplot [color=mycolor1, line width=1pt, only marks, mark size=1.0pt, mark=o, mark options={solid, draw=mycolor1}]
  table[row sep=crcr]{%
0.03125	1.28841879576402e-06\\
0.015625	8.24509422958641e-08\\
0.0078125	1.37060700566849e-08\\
0.00390625	1.68369501674093e-09\\
};
\addlegendentry{Lagrange-bary degree 2, $\alpha = 3.13$}

\addplot [color=mycolor1, dashed, line width=1pt, forget plot]
  table[row sep=crcr]{%
0.03125	1.02791403294352e-06\\
0.015625	1.17189969125886e-07\\
0.0078125	1.3360542247292e-08\\
0.00390625	1.52320280031752e-09\\
};
\addplot [color=mycolor2, line width=1pt, only marks, mark size=1pt, mark=square, mark options={solid, draw=mycolor2}]
  table[row sep=crcr]{%
0.03125	1.27837135666109e-06\\
0.015625	4.41329117695233e-08\\
0.0078125	2.08801861911745e-09\\
0.00390625	7.68460572952989e-11\\
};
\addlegendentry{Lagrange-bary degree 3, $\alpha = 4.65$}

\addplot [color=mycolor2, dashed, line width=1pt, forget plot]
  table[row sep=crcr]{%
0.03125	1.2229651271062e-06\\
0.015625	4.88203197491138e-08\\
0.0078125	1.9488892754002e-09\\
0.00390625	7.77989457522726e-11\\
};
\addplot [color=mycolor3, line width=1pt, only marks, mark size=1.0pt, mark=triangle, mark options={solid, draw=mycolor3}]
  table[row sep=crcr]{%
0.03125	1.30662846442586e-06\\
0.015625	4.25519337543534e-08\\
0.0078125	1.96387581197222e-09\\
0.00390625	6.6803396148174e-11\\
};
\addlegendentry{B-spline degree 3, $\alpha = 4.72$}

\addplot [color=mycolor3, dashed, line width=1pt, forget plot]
  table[row sep=crcr]{%
0.03125	1.25090951129996e-06\\
0.015625	4.74503362922385e-08\\
0.0078125	1.79991789486571e-09\\
0.00390625	6.82756895189328e-11\\
};
\end{axis}
\end{tikzpicture}%

%% file: figures/incomp_err_landau_damping.tex
\definecolor{mycolor1}{rgb}{0.85000,0.32500,0.09800}%
\begin{tikzpicture}

\begin{axis}[%
width=0.953\figurewidth,
height=0.163\figureheight,
at={(0\figurewidth,0\figureheight)},
scale only axis,
xmin=0.00000000000000,
xmax=50.00000000000000,
xminorticks=true,
xlabel style={font=\color{white!15!black}},
xlabel={Time $t$},
ymode=log,
ymin=0.00000000000001,
ymax=0.01996425963199,
yminorticks=true,
ylabel style={font=\color{white!15!black}},
ylabel={\small Incomp. Err.},
axis background/.style={fill=white},
axis x line*=bottom,
axis y line*=left,
xmajorgrids,
ymajorgrids,
yminorgrids,
grid style={opacity=0.3},
minor grid style={opacity=0.1},
legend style={at={(1.02,0.03)}, anchor=south west, legend cell align=left, align=left, fill=none, draw=none},
scaled ticks=false, xticklabel style={/pgf/number format/fixed},yticklabel style={/pgf/number format/fixed}
]
\addplot [color=mycolor1, line width=1.5pt]
  table[row sep=crcr]{%
0.10000000000000	0.00000000230092\\
0.20000000000000	0.00000011377721\\
0.30000000000000	0.00000038198202\\
0.40000000000000	0.00000080795456\\
0.50000000000000	0.00000134534210\\
0.60000000000000	0.00000186374881\\
0.70000000000000	0.00000234953130\\
0.80000000000000	0.00000269307905\\
0.90000000000000	0.00000283892361\\
1.00000000000000	0.00000283050429\\
1.10000000000000	0.00000267867999\\
1.20000000000000	0.00000247947878\\
1.30000000000000	0.00000228127156\\
1.40000000000000	0.00000212008546\\
1.50000000000000	0.00000206394078\\
1.60000000000000	0.00000218334395\\
1.70000000000000	0.00000240129373\\
1.80000000000000	0.00000282379130\\
1.90000000000000	0.00000312804997\\
2.00000000000000	0.00000314470638\\
2.10000000000000	0.00000000171987\\
2.20000000000000	0.00000004456072\\
2.30000000000000	0.00000014822736\\
2.40000000000000	0.00000032065839\\
2.50000000000000	0.00000056082659\\
2.60000000000000	0.00000082115922\\
2.70000000000000	0.00000110742964\\
2.80000000000000	0.00000133774394\\
2.90000000000000	0.00000149209373\\
3.00000000000000	0.00000161163201\\
3.10000000000000	0.00000160807442\\
3.20000000000000	0.00000156019197\\
3.30000000000000	0.00000146405415\\
3.40000000000000	0.00000130873977\\
3.50000000000000	0.00000110924319\\
3.60000000000000	0.00000096040928\\
3.70000000000000	0.00000095136537\\
3.80000000000000	0.00000104001531\\
3.90000000000000	0.00000113967823\\
4.00000000000000	0.00000120865994\\
4.10000000000000	0.00000000074680\\
4.20000000000000	0.00000002073338\\
4.30000000000000	0.00000007130985\\
4.40000000000000	0.00000016056742\\
4.50000000000000	0.00000028923484\\
4.60000000000000	0.00000044200039\\
4.70000000000000	0.00000062644987\\
4.80000000000000	0.00000078827234\\
4.90000000000000	0.00000094791261\\
5.00000000000000	0.00000107212921\\
5.10000000000000	0.00000114676747\\
5.20000000000000	0.00000114669828\\
5.30000000000000	0.00000110271804\\
5.40000000000000	0.00000101898862\\
5.50000000000000	0.00000087972181\\
5.60000000000000	0.00000070950453\\
5.70000000000000	0.00000071940068\\
5.80000000000000	0.00000076511157\\
5.90000000000000	0.00000081141688\\
6.00000000000000	0.00000085702665\\
6.10000000000000	0.00000000014441\\
6.20000000000000	0.00000000676280\\
6.30000000000000	0.00000002715226\\
6.40000000000000	0.00000006593937\\
6.50000000000000	0.00000012721545\\
6.60000000000000	0.00000021004663\\
6.70000000000000	0.00000031140573\\
6.80000000000000	0.00000042811619\\
6.90000000000000	0.00000054134630\\
7.00000000000000	0.00000065904339\\
7.10000000000000	0.00000076739654\\
7.20000000000000	0.00000084444422\\
7.30000000000000	0.00000087399020\\
7.40000000000000	0.00000085564699\\
7.50000000000000	0.00000076176204\\
7.60000000000000	0.00000066917684\\
7.70000000000000	0.00000055327966\\
7.80000000000000	0.00000055121247\\
7.90000000000000	0.00000059481573\\
8.00000000000000	0.00000063760718\\
8.10000000000000	0.00000000019778\\
8.20000000000000	0.00000000147548\\
8.30000000000000	0.00000000149727\\
8.40000000000000	0.00000001255533\\
8.50000000000000	0.00000003649547\\
8.60000000000000	0.00000007418204\\
8.70000000000000	0.00000012823348\\
8.80000000000000	0.00000019759018\\
8.90000000000000	0.00000028090559\\
9.00000000000000	0.00000036526964\\
9.10000000000000	0.00000046011766\\
9.20000000000000	0.00000054159772\\
9.30000000000000	0.00000062537163\\
9.40000000000000	0.00000068719121\\
9.50000000000000	0.00000070426609\\
9.60000000000000	0.00000068464184\\
9.70000000000000	0.00000061834373\\
9.80000000000000	0.00000053321312\\
9.90000000000000	0.00000043739416\\
10.00000000000000	0.00000043541294\\
10.10000000000000	0.00000000035355\\
10.20000000000000	0.00000000551093\\
10.30000000000000	0.00000001250478\\
10.40000000000000	0.00000001707784\\
10.50000000000000	0.00000001689002\\
10.60000000000000	0.00000001600954\\
10.70000000000000	0.00000003509610\\
10.80000000000000	0.00000007498449\\
10.90000000000000	0.00000013165213\\
11.00000000000000	0.00000019834036\\
11.10000000000000	0.00000027572526\\
11.20000000000000	0.00000035091077\\
11.30000000000000	0.00000042369658\\
11.40000000000000	0.00000048370131\\
11.50000000000000	0.00000052922369\\
11.60000000000000	0.00000057241331\\
11.70000000000000	0.00000059461355\\
11.80000000000000	0.00000058006942\\
11.90000000000000	0.00000053141178\\
12.00000000000000	0.00000042910564\\
12.10000000000000	0.00000000038575\\
12.20000000000000	0.00000000694592\\
12.30000000000000	0.00000001815493\\
12.40000000000000	0.00000003074210\\
12.50000000000000	0.00000004089348\\
12.60000000000000	0.00000004831622\\
12.70000000000000	0.00000005467076\\
12.80000000000000	0.00000006620473\\
12.90000000000000	0.00000008657227\\
13.00000000000000	0.00000012345046\\
13.10000000000000	0.00000017571923\\
13.20000000000000	0.00000023618952\\
13.30000000000000	0.00000029877409\\
13.40000000000000	0.00000035492165\\
13.50000000000000	0.00000040357162\\
13.60000000000000	0.00000043933290\\
13.70000000000000	0.00000047585248\\
13.80000000000000	0.00000049433242\\
13.90000000000000	0.00000047494974\\
14.00000000000000	0.00000041891063\\
14.10000000000000	0.00000000034628\\
14.20000000000000	0.00000000669023\\
14.30000000000000	0.00000001851615\\
14.40000000000000	0.00000003351429\\
14.50000000000000	0.00000004797935\\
14.60000000000000	0.00000006059301\\
14.70000000000000	0.00000007210319\\
14.80000000000000	0.00000008189418\\
14.90000000000000	0.00000009343169\\
15.00000000000000	0.00000010858126\\
15.10000000000000	0.00000013494636\\
15.20000000000000	0.00000016949041\\
15.30000000000000	0.00000020670468\\
15.40000000000000	0.00000024977580\\
15.50000000000000	0.00000028639275\\
15.60000000000000	0.00000031380014\\
15.70000000000000	0.00000035054370\\
15.80000000000000	0.00000037481582\\
15.90000000000000	0.00000037061293\\
16.00000000000000	0.00000033570398\\
16.10000000000000	0.00000000027240\\
16.20000000000000	0.00000000553033\\
16.30000000000000	0.00000001594491\\
16.40000000000000	0.00000002976543\\
16.50000000000000	0.00000004405062\\
16.60000000000000	0.00000005946317\\
16.70000000000000	0.00000007268509\\
16.80000000000000	0.00000008332942\\
16.90000000000000	0.00000009289445\\
17.00000000000000	0.00000010064925\\
17.10000000000000	0.00000010996809\\
17.20000000000000	0.00000011876062\\
17.30000000000000	0.00000013092278\\
17.40000000000000	0.00000014874847\\
17.50000000000000	0.00000016720324\\
17.60000000000000	0.00000019993081\\
17.70000000000000	0.00000023094104\\
17.80000000000000	0.00000025301973\\
17.90000000000000	0.00000025872510\\
18.00000000000000	0.00000024338403\\
18.10000000000000	0.00000000019021\\
18.20000000000000	0.00000000405084\\
18.30000000000000	0.00000001205749\\
18.40000000000000	0.00000002321193\\
18.50000000000000	0.00000003538898\\
18.60000000000000	0.00000004837996\\
18.70000000000000	0.00000006152889\\
18.80000000000000	0.00000007209582\\
18.90000000000000	0.00000008064957\\
19.00000000000000	0.00000008729320\\
19.10000000000000	0.00000009086412\\
19.20000000000000	0.00000009456813\\
19.30000000000000	0.00000009385280\\
19.40000000000000	0.00000009510539\\
19.50000000000000	0.00000010108533\\
19.60000000000000	0.00000011435634\\
19.70000000000000	0.00000013404164\\
19.80000000000000	0.00000015010283\\
19.90000000000000	0.00000015768109\\
20.00000000000000	0.00000015320612\\
20.10000000000000	0.00000000011492\\
20.20000000000000	0.00000000258965\\
20.30000000000000	0.00000000809193\\
20.40000000000000	0.00000001575945\\
20.50000000000000	0.00000002485351\\
20.60000000000000	0.00000003544549\\
20.70000000000000	0.00000004529882\\
20.80000000000000	0.00000005354535\\
20.90000000000000	0.00000006272715\\
21.00000000000000	0.00000006936162\\
21.10000000000000	0.00000007291986\\
21.20000000000000	0.00000007485362\\
21.30000000000000	0.00000007471672\\
21.40000000000000	0.00000007505149\\
21.50000000000000	0.00000006705147\\
21.60000000000000	0.00000006659776\\
21.70000000000000	0.00000007489703\\
21.80000000000000	0.00000007655643\\
21.90000000000000	0.00000008898854\\
22.00000000000000	0.00000009416822\\
22.10000000000000	0.00000000005457\\
22.20000000000000	0.00000000135573\\
22.30000000000000	0.00000000436326\\
22.40000000000000	0.00000000907166\\
22.50000000000000	0.00000001493828\\
22.60000000000000	0.00000002232869\\
22.70000000000000	0.00000002951148\\
22.80000000000000	0.00000003663693\\
22.90000000000000	0.00000004333813\\
23.00000000000000	0.00000005189441\\
23.10000000000000	0.00000006654615\\
23.20000000000000	0.00000008007523\\
23.30000000000000	0.00000006952474\\
23.40000000000000	0.00000012520328\\
23.50000000000000	0.00000007304292\\
23.60000000000000	0.00000005959727\\
23.70000000000000	0.00000011781309\\
23.80000000000000	0.00000008569514\\
23.90000000000000	0.00000007757309\\
24.00000000000000	0.00000012760186\\
24.10000000000000	0.00000000001209\\
24.20000000000000	0.00000000042958\\
24.30000000000000	0.00000000162228\\
24.40000000000000	0.00000000382363\\
24.50000000000000	0.00000000686389\\
24.60000000000000	0.00000001103441\\
24.70000000000000	0.00000001597113\\
24.80000000000000	0.00000002158477\\
24.90000000000000	0.00000002729996\\
25.00000000000000	0.00000003309067\\
25.10000000000000	0.00000005369284\\
25.20000000000000	0.00000009282494\\
25.30000000000000	0.00000010361961\\
25.40000000000000	0.00000027538519\\
25.50000000000000	0.00000024578429\\
25.60000000000000	0.00000012208843\\
25.70000000000000	0.00000038252192\\
25.80000000000000	0.00000032327396\\
25.90000000000000	0.00000016622161\\
26.00000000000000	0.00000029360584\\
26.10000000000000	0.00000000001769\\
26.20000000000000	0.00000000018222\\
26.30000000000000	0.00000000023113\\
26.40000000000000	0.00000000018803\\
26.50000000000000	0.00000000112930\\
26.60000000000000	0.00000000288090\\
26.70000000000000	0.00000000554037\\
26.80000000000000	0.00000000879757\\
26.90000000000000	0.00000001344016\\
27.00000000000000	0.00000001875181\\
27.10000000000000	0.00000003664869\\
27.20000000000000	0.00000008235995\\
27.30000000000000	0.00000017594339\\
27.40000000000000	0.00000045193476\\
27.50000000000000	0.00000042013493\\
27.60000000000000	0.00000025300247\\
27.70000000000000	0.00000090966646\\
27.80000000000000	0.00000118540370\\
27.90000000000000	0.00000056590398\\
28.00000000000000	0.00000084233382\\
28.10000000000000	0.00000000003125\\
28.20000000000000	0.00000000051256\\
28.30000000000000	0.00000000127926\\
28.40000000000000	0.00000000208167\\
28.50000000000000	0.00000000272421\\
28.60000000000000	0.00000000308812\\
28.70000000000000	0.00000000349399\\
28.80000000000000	0.00000000442886\\
28.90000000000000	0.00000000638970\\
29.00000000000000	0.00000000966001\\
29.10000000000000	0.00000001943478\\
29.20000000000000	0.00000006467202\\
29.30000000000000	0.00000019973328\\
29.40000000000000	0.00000045031451\\
29.50000000000000	0.00000059450752\\
29.60000000000000	0.00000122311114\\
29.70000000000000	0.00000218331008\\
29.80000000000000	0.00000194836984\\
29.90000000000000	0.00000103063225\\
30.00000000000000	0.00000288004028\\
30.10000000000000	0.00000000004810\\
30.20000000000000	0.00000000062410\\
30.30000000000000	0.00000000169466\\
30.40000000000000	0.00000000302765\\
30.50000000000000	0.00000000437240\\
30.60000000000000	0.00000000568476\\
30.70000000000000	0.00000000664380\\
30.80000000000000	0.00000000768927\\
30.90000000000000	0.00000000863355\\
31.00000000000000	0.00000001043232\\
31.10000000000000	0.00000001465209\\
31.20000000000000	0.00000004198799\\
31.30000000000000	0.00000015382717\\
31.40000000000000	0.00000042990776\\
31.50000000000000	0.00000097770915\\
31.60000000000000	0.00000232564464\\
31.70000000000000	0.00000291869010\\
31.80000000000000	0.00000227146696\\
31.90000000000000	0.00000543331242\\
32.00000000000000	0.00000750347482\\
32.10000000000000	0.00000000006342\\
32.20000000000000	0.00000000061016\\
32.30000000000000	0.00000000161629\\
32.40000000000000	0.00000000310442\\
32.50000000000000	0.00000000475285\\
32.60000000000000	0.00000000633280\\
32.70000000000000	0.00000000781505\\
32.80000000000000	0.00000000923263\\
32.90000000000000	0.00000001149222\\
33.00000000000000	0.00000001271827\\
33.10000000000000	0.00000001602920\\
33.20000000000000	0.00000002429539\\
33.30000000000000	0.00000008677227\\
33.40000000000000	0.00000028677030\\
33.50000000000000	0.00000082344468\\
33.60000000000000	0.00000182335538\\
33.70000000000000	0.00000240051194\\
33.80000000000000	0.00000476341139\\
33.90000000000000	0.00000980894574\\
34.00000000000000	0.00000920464521\\
34.10000000000000	0.00000000003399\\
34.20000000000000	0.00000000045251\\
34.30000000000000	0.00000000141110\\
34.40000000000000	0.00000000259437\\
34.50000000000000	0.00000000406066\\
34.60000000000000	0.00000000572521\\
34.70000000000000	0.00000000738758\\
34.80000000000000	0.00000000880689\\
34.90000000000000	0.00000001034641\\
35.00000000000000	0.00000001616454\\
35.10000000000000	0.00000001655860\\
35.20000000000000	0.00000002109775\\
35.30000000000000	0.00000003747298\\
35.40000000000000	0.00000011097755\\
35.50000000000000	0.00000031827268\\
35.60000000000000	0.00000057461674\\
35.70000000000000	0.00000053654651\\
35.80000000000000	0.00000194872221\\
35.90000000000000	0.00000214375040\\
36.00000000000000	0.00000571297507\\
36.10000000000000	0.00000000006394\\
36.20000000000000	0.00000000034412\\
36.30000000000000	0.00000000093960\\
36.40000000000000	0.00000000197423\\
36.50000000000000	0.00000000311410\\
36.60000000000000	0.00000000439373\\
36.70000000000000	0.00000000612751\\
36.80000000000000	0.00000000751844\\
36.90000000000000	0.00000001073008\\
37.00000000000000	0.00000001348505\\
37.10000000000000	0.00000003200741\\
37.20000000000000	0.00000002847331\\
37.30000000000000	0.00000003677116\\
37.40000000000000	0.00000004129404\\
37.50000000000000	0.00000004810218\\
37.60000000000000	0.00000031349101\\
37.70000000000000	0.00000145707419\\
37.80000000000000	0.00000459015898\\
37.90000000000000	0.00000983196509\\
38.00000000000000	0.00001495899201\\
38.10000000000000	0.00000000004452\\
38.20000000000000	0.00000000023904\\
38.30000000000000	0.00000000067653\\
38.40000000000000	0.00000000114770\\
38.50000000000000	0.00000000163117\\
38.60000000000000	0.00000000264284\\
38.70000000000000	0.00000000390460\\
38.80000000000000	0.00000000553030\\
38.90000000000000	0.00000000812870\\
39.00000000000000	0.00000001441107\\
39.10000000000000	0.00000002431824\\
39.20000000000000	0.00000009155593\\
39.30000000000000	0.00000003921071\\
39.40000000000000	0.00000009419084\\
39.50000000000000	0.00000015652048\\
39.60000000000000	0.00000039561320\\
39.70000000000000	0.00000130302704\\
39.80000000000000	0.00000340490250\\
39.90000000000000	0.00000668590662\\
40.00000000000000	0.00000767744765\\
40.10000000000000	0.00000000011370\\
40.20000000000000	0.00000000033437\\
40.30000000000000	0.00000000043924\\
40.40000000000000	0.00000000069976\\
40.50000000000000	0.00000000117236\\
40.60000000000000	0.00000000095016\\
40.70000000000000	0.00000000233010\\
40.80000000000000	0.00000000340778\\
40.90000000000000	0.00000000468576\\
41.00000000000000	0.00000001883058\\
41.10000000000000	0.00000005730929\\
41.20000000000000	0.00000016329635\\
41.30000000000000	0.00000012532847\\
41.40000000000000	0.00000012484627\\
41.50000000000000	0.00000020584282\\
41.60000000000000	0.00000029738051\\
41.70000000000000	0.00000032372007\\
41.80000000000000	0.00000036525521\\
41.90000000000000	0.00000067329269\\
42.00000000000000	0.00000378839033\\
42.10000000000000	0.00000000003564\\
42.20000000000000	0.00000000023151\\
42.30000000000000	0.00000000065794\\
42.40000000000000	0.00000000111554\\
42.50000000000000	0.00000000132150\\
42.60000000000000	0.00000000111380\\
42.70000000000000	0.00000000139102\\
42.80000000000000	0.00000000253473\\
42.90000000000000	0.00000000508576\\
43.00000000000000	0.00000001537364\\
43.10000000000000	0.00000006307314\\
43.20000000000000	0.00000016637571\\
43.30000000000000	0.00000025089576\\
43.40000000000000	0.00000027699355\\
43.50000000000000	0.00000042529296\\
43.60000000000000	0.00000032533371\\
43.70000000000000	0.00000053428225\\
43.80000000000000	0.00000157648751\\
43.90000000000000	0.00000295369862\\
44.00000000000000	0.00000424145922\\
44.10000000000000	0.00000000010875\\
44.20000000000000	0.00000000055182\\
44.30000000000000	0.00000000124230\\
44.40000000000000	0.00000000172532\\
44.50000000000000	0.00000000221248\\
44.60000000000000	0.00000000195394\\
44.70000000000000	0.00000000297327\\
44.80000000000000	0.00000000445971\\
44.90000000000000	0.00000000624398\\
45.00000000000000	0.00000001376289\\
45.10000000000000	0.00000003899178\\
45.20000000000000	0.00000013174913\\
45.30000000000000	0.00000029110817\\
45.40000000000000	0.00000049261584\\
45.50000000000000	0.00000090328965\\
45.60000000000000	0.00000056827416\\
45.70000000000000	0.00000091275756\\
45.80000000000000	0.00000144899830\\
45.90000000000000	0.00000187942164\\
46.00000000000000	0.00000285796702\\
46.10000000000000	0.00000000006741\\
46.20000000000000	0.00000000047725\\
46.30000000000000	0.00000000157564\\
46.40000000000000	0.00000000284294\\
46.50000000000000	0.00000000293041\\
46.60000000000000	0.00000000296289\\
46.70000000000000	0.00000000366788\\
46.80000000000000	0.00000000511368\\
46.90000000000000	0.00000001181559\\
47.00000000000000	0.00000001569807\\
47.10000000000000	0.00000002156623\\
47.20000000000000	0.00000009098384\\
47.30000000000000	0.00000033213916\\
47.40000000000000	0.00000075496436\\
47.50000000000000	0.00000113405628\\
47.60000000000000	0.00000066605645\\
47.70000000000000	0.00000157575465\\
47.80000000000000	0.00000201242380\\
47.90000000000000	0.00000112900469\\
48.00000000000000	0.00000099360694\\
48.10000000000000	0.00000000026188\\
48.20000000000000	0.00000000124191\\
48.30000000000000	0.00000000237532\\
48.40000000000000	0.00000000301489\\
48.50000000000000	0.00000000412027\\
48.60000000000000	0.00000000641808\\
48.70000000000000	0.00000000638619\\
48.80000000000000	0.00000000649392\\
48.90000000000000	0.00000001095173\\
49.00000000000000	0.00000002884585\\
49.10000000000000	0.00000004331807\\
49.20000000000000	0.00000009010394\\
49.30000000000000	0.00000034668328\\
49.40000000000000	0.00000101220115\\
49.50000000000000	0.00000166915701\\
49.60000000000000	0.00000178181853\\
49.70000000000000	0.00000192678707\\
49.80000000000000	0.00000249840727\\
49.90000000000000	0.00000340960380\\
50.00000000000000	0.00000292088907\\
};
\addlegendentry{CMM-Nufi $\Nmapgrid=256$}
\addplot [color=black, line width=1.5pt]
  table[row sep=crcr]{%
0.10000000000000	0.00000097743551\\
0.20000000000000	0.00000454647588\\
0.30000000000000	0.00001103262831\\
0.40000000000000	0.00002895270134\\
0.50000000000000	0.00006443995415\\
0.60000000000000	0.00012378573715\\
0.70000000000000	0.00020289408683\\
0.80000000000000	0.00031153505481\\
0.90000000000000	0.00046604325360\\
1.00000000000000	0.00067745839776\\
1.10000000000000	0.00092277882623\\
1.20000000000000	0.00117883934988\\
1.30000000000000	0.00142943021778\\
1.40000000000000	0.00160973239402\\
1.50000000000000	0.00166705145435\\
1.60000000000000	0.00170504352644\\
1.70000000000000	0.00181378847862\\
1.80000000000000	0.00215834160198\\
1.90000000000000	0.00344715508047\\
2.00000000000000	0.00620070070131\\
2.20000000000000	0.00000018709181\\
2.30000000000000	0.00000095670322\\
2.40000000000000	0.00000262751852\\
2.50000000000000	0.00000710350148\\
2.60000000000000	0.00001774647876\\
2.70000000000000	0.00003617774072\\
2.80000000000000	0.00006395681478\\
2.90000000000000	0.00010575759946\\
3.00000000000000	0.00016518083555\\
3.10000000000000	0.00024095358973\\
3.20000000000000	0.00032475909160\\
3.30000000000000	0.00043324905225\\
3.40000000000000	0.00055611032512\\
3.50000000000000	0.00068466537682\\
3.60000000000000	0.00078839879226\\
3.70000000000000	0.00083545855123\\
3.80000000000000	0.00081357223514\\
3.90000000000000	0.00082318007017\\
4.00000000000000	0.00102230236144\\
4.10000000000000	0.00141431450680\\
4.20000000000000	0.00227807450813\\
4.30000000000000	0.00374942135929\\
4.40000000000000	0.00568999345294\\
4.50000000000000	0.00813246613651\\
4.60000000000000	0.00000000000167\\
4.70000000000000	0.00000015735958\\
4.80000000000000	0.00000064061154\\
4.90000000000000	0.00000186281035\\
5.00000000000000	0.00000471491383\\
5.10000000000000	0.00001015213737\\
5.20000000000000	0.00001962250988\\
5.30000000000000	0.00003231645346\\
5.40000000000000	0.00004852145729\\
5.50000000000000	0.00006847500898\\
5.60000000000000	0.00009370746977\\
5.70000000000000	0.00012094836309\\
5.80000000000000	0.00014391013936\\
5.90000000000000	0.00015683544182\\
6.00000000000000	0.00016024379517\\
6.10000000000000	0.00018307201652\\
6.20000000000000	0.00024796022142\\
6.30000000000000	0.00035524934105\\
6.40000000000000	0.00055543028355\\
6.50000000000000	0.00095203037566\\
6.60000000000000	0.00149404851958\\
6.70000000000000	0.00219995855132\\
6.80000000000000	0.00308114179468\\
6.90000000000000	0.00413878407191\\
7.00000000000000	0.00536120751396\\
7.10000000000000	0.00672220358694\\
7.20000000000000	0.00818090045003\\
7.30000000000000	0.00968353490918\\
7.40000000000000	0.00000001069805\\
7.50000000000000	0.00000010778034\\
7.60000000000000	0.00000031474437\\
7.70000000000000	0.00000073481470\\
7.80000000000000	0.00000140280285\\
7.90000000000000	0.00000235958577\\
8.00000000000000	0.00000335600500\\
8.10000000000000	0.00000403433448\\
8.20000000000000	0.00000432793005\\
8.30000000000000	0.00000628007555\\
8.40000000000000	0.00001222355105\\
8.50000000000000	0.00002246036552\\
8.60000000000000	0.00003877725035\\
8.70000000000000	0.00007014753595\\
8.80000000000000	0.00011726143598\\
8.90000000000000	0.00018266745570\\
9.00000000000000	0.00026906526838\\
9.10000000000000	0.00037815665561\\
9.20000000000000	0.00051016567768\\
9.30000000000000	0.00067018610558\\
9.40000000000000	0.00085419540198\\
9.50000000000000	0.00105610733219\\
9.60000000000000	0.00126877563910\\
9.70000000000000	0.00148312882241\\
9.80000000000000	0.00168886066254\\
9.90000000000000	0.00187541045043\\
10.00000000000000	0.00203314495519\\
10.10000000000000	0.00217555174222\\
10.20000000000000	0.00227948183756\\
10.30000000000000	0.00233605145165\\
10.40000000000000	0.00235069910832\\
10.50000000000000	0.00242955646567\\
10.60000000000000	0.00262815610110\\
10.70000000000000	0.00297001968595\\
10.80000000000000	0.00347303916277\\
10.90000000000000	0.00414693767436\\
11.00000000000000	0.00499149368988\\
11.10000000000000	0.00599570560019\\
11.20000000000000	0.00713807889476\\
11.30000000000000	0.00838792207041\\
11.40000000000000	0.00974650181443\\
11.50000000000000	0.00000000371683\\
11.60000000000000	0.00000004038512\\
11.70000000000000	0.00000013732806\\
11.80000000000000	0.00000034599571\\
11.90000000000000	0.00000072761451\\
12.00000000000000	0.00000139968494\\
12.10000000000000	0.00000235070823\\
12.20000000000000	0.00000363886112\\
12.30000000000000	0.00000501621419\\
12.40000000000000	0.00000631442466\\
12.50000000000000	0.00000766361505\\
12.60000000000000	0.00000822806988\\
12.70000000000000	0.00000883345681\\
12.80000000000000	0.00001223052950\\
12.90000000000000	0.00001868833222\\
13.00000000000000	0.00002893042240\\
13.10000000000000	0.00004664883914\\
13.20000000000000	0.00007889873855\\
13.30000000000000	0.00012330339573\\
13.40000000000000	0.00018165666514\\
13.50000000000000	0.00025527228794\\
13.60000000000000	0.00034474768767\\
13.70000000000000	0.00044974664274\\
13.80000000000000	0.00056883124696\\
13.90000000000000	0.00069937455736\\
14.00000000000000	0.00083758042430\\
14.10000000000000	0.00097862827999\\
14.20000000000000	0.00111694709024\\
14.30000000000000	0.00124660681366\\
14.40000000000000	0.00136179995608\\
14.50000000000000	0.00146846024978\\
14.60000000000000	0.00155385774860\\
14.70000000000000	0.00163894179267\\
14.80000000000000	0.00167564035671\\
14.90000000000000	0.00168420693745\\
15.00000000000000	0.00174311261080\\
15.10000000000000	0.00188317347193\\
15.20000000000000	0.00212783630674\\
15.30000000000000	0.00247056112081\\
15.40000000000000	0.00291526638394\\
15.50000000000000	0.00346158421482\\
15.60000000000000	0.00411344650125\\
15.70000000000000	0.00489194059896\\
15.80000000000000	0.00575822039475\\
15.90000000000000	0.00669343593030\\
16.00000000000000	0.00767470834848\\
16.10000000000000	0.00867628812528\\
16.20000000000000	0.00967087157121\\
16.30000000000000	0.00000000076310\\
16.40000000000000	0.00000000899995\\
16.50000000000000	0.00000002801173\\
16.60000000000000	0.00000006298202\\
16.70000000000000	0.00000011657166\\
16.80000000000000	0.00000018391957\\
16.90000000000000	0.00000025573510\\
17.00000000000000	0.00000032515802\\
17.10000000000000	0.00000036688620\\
17.20000000000000	0.00000038412471\\
17.30000000000000	0.00000069562026\\
17.40000000000000	0.00000129132515\\
17.50000000000000	0.00000224243473\\
17.60000000000000	0.00000380551848\\
17.70000000000000	0.00000650263124\\
17.80000000000000	0.00001029745545\\
17.90000000000000	0.00001537831466\\
18.00000000000000	0.00002188945162\\
18.10000000000000	0.00002996335948\\
18.20000000000000	0.00003983142028\\
18.30000000000000	0.00005128493473\\
18.40000000000000	0.00006412300079\\
18.50000000000000	0.00007801225011\\
18.60000000000000	0.00009249290961\\
18.70000000000000	0.00010699941840\\
18.80000000000000	0.00012089571144\\
18.90000000000000	0.00013352379133\\
19.00000000000000	0.00014426250477\\
19.10000000000000	0.00015259167699\\
19.20000000000000	0.00015815589334\\
19.30000000000000	0.00016082144195\\
19.40000000000000	0.00016182808999\\
19.50000000000000	0.00016761917420\\
19.60000000000000	0.00018035736800\\
19.70000000000000	0.00020135003371\\
19.80000000000000	0.00023367004824\\
19.90000000000000	0.00027698757578\\
20.00000000000000	0.00033108268563\\
20.10000000000000	0.00039573999109\\
20.20000000000000	0.00047017408430\\
20.30000000000000	0.00055303147947\\
20.40000000000000	0.00064243599329\\
20.50000000000000	0.00073607601097\\
20.60000000000000	0.00083132756181\\
20.70000000000000	0.00092540651604\\
20.80000000000000	0.00101553806126\\
20.90000000000000	0.00110231617462\\
21.00000000000000	0.00118693271852\\
21.10000000000000	0.00126079490869\\
21.20000000000000	0.00132222013863\\
21.30000000000000	0.00137025917617\\
21.40000000000000	0.00140477489807\\
21.50000000000000	0.00142647168612\\
21.60000000000000	0.00143687284619\\
21.70000000000000	0.00144092832916\\
21.80000000000000	0.00150351029828\\
21.90000000000000	0.00159172703077\\
22.00000000000000	0.00173203038668\\
22.10000000000000	0.00190327843673\\
22.20000000000000	0.00210206802573\\
22.30000000000000	0.00234791439969\\
22.40000000000000	0.00263502402359\\
22.50000000000000	0.00291533866148\\
22.60000000000000	0.00324692514767\\
22.70000000000000	0.00355037806384\\
22.80000000000000	0.00386247550808\\
22.90000000000000	0.00413087391152\\
23.00000000000000	0.00438312311483\\
23.10000000000000	0.00455218891926\\
23.20000000000000	0.00487179831031\\
23.30000000000000	0.00529866496849\\
23.40000000000000	0.00569989325405\\
23.50000000000000	0.00616307900036\\
23.60000000000000	0.00657174179509\\
23.70000000000000	0.00711006000853\\
23.80000000000000	0.00760947042921\\
23.90000000000000	0.00821065081226\\
24.00000000000000	0.00884984877668\\
24.10000000000000	0.00954260714528\\
24.20000000000000	0.00000000000002\\
24.30000000000000	0.00000000222502\\
24.40000000000000	0.00000000989218\\
24.50000000000000	0.00000002793288\\
24.60000000000000	0.00000006559882\\
24.70000000000000	0.00000013537227\\
24.80000000000000	0.00000027891504\\
24.90000000000000	0.00000053441002\\
25.00000000000000	0.00000096710712\\
25.10000000000000	0.00000167688926\\
25.20000000000000	0.00000273294178\\
25.30000000000000	0.00000434780581\\
25.40000000000000	0.00000682773720\\
25.50000000000000	0.00001034821195\\
25.60000000000000	0.00001519236376\\
25.70000000000000	0.00002249241648\\
25.80000000000000	0.00003341762940\\
25.90000000000000	0.00004870057864\\
26.00000000000000	0.00006971948956\\
26.10000000000000	0.00009869366104\\
26.20000000000000	0.00013861353879\\
26.30000000000000	0.00019184248067\\
26.40000000000000	0.00026184712864\\
26.50000000000000	0.00035270393432\\
26.60000000000000	0.00046912467629\\
26.70000000000000	0.00061645674810\\
26.80000000000000	0.00080542371895\\
26.90000000000000	0.00104770825841\\
27.00000000000000	0.00134830474642\\
27.10000000000000	0.00171717473524\\
27.20000000000000	0.00216491659701\\
27.30000000000000	0.00270248388518\\
27.40000000000000	0.00334077273193\\
27.50000000000000	0.00417843626682\\
27.60000000000000	0.00532509048227\\
27.70000000000000	0.00674248422783\\
27.80000000000000	0.00846858236837\\
27.90000000000000	0.00000000000190\\
28.00000000000000	0.00000017781295\\
28.10000000000000	0.00000076695750\\
28.20000000000000	0.00000237494252\\
28.30000000000000	0.00000703957689\\
28.40000000000000	0.00001721685625\\
28.50000000000000	0.00003720434801\\
28.60000000000000	0.00007040723326\\
28.70000000000000	0.00012666886565\\
28.80000000000000	0.00021621117106\\
28.90000000000000	0.00035316394559\\
29.00000000000000	0.00056265209323\\
29.10000000000000	0.00085954567578\\
29.20000000000000	0.00132358740163\\
29.30000000000000	0.00195911238770\\
29.40000000000000	0.00284521260610\\
29.50000000000000	0.00402986154404\\
29.60000000000000	0.00560346597188\\
29.70000000000000	0.00765671395580\\
29.90000000000000	0.00000180826194\\
30.00000000000000	0.00000897588101\\
30.10000000000000	0.00002332225531\\
30.20000000000000	0.00006855948852\\
30.30000000000000	0.00017441141734\\
30.40000000000000	0.00038403514228\\
30.50000000000000	0.00078975469100\\
30.60000000000000	0.00143530659222\\
30.70000000000000	0.00244888579034\\
30.80000000000000	0.00391347357093\\
30.90000000000000	0.00581060425731\\
31.00000000000000	0.00835753301559\\
31.10000000000000	0.00000000017817\\
31.20000000000000	0.00001350127252\\
31.30000000000000	0.00005815276000\\
31.40000000000000	0.00015680451090\\
31.50000000000000	0.00047754976483\\
31.60000000000000	0.00122297751935\\
31.70000000000000	0.00250480642863\\
31.80000000000000	0.00476555576564\\
31.90000000000000	0.00791050259069\\
32.00000000000000	0.00000000074515\\
32.10000000000000	0.00004571296915\\
32.20000000000000	0.00018932765148\\
32.30000000000000	0.00057184943586\\
32.40000000000000	0.00166331421705\\
32.50000000000000	0.00408278576225\\
32.60000000000000	0.00774147358493\\
32.70000000000000	0.00000000257888\\
32.80000000000000	0.00011542469760\\
32.90000000000000	0.00046135401868\\
33.00000000000000	0.00155072186188\\
33.10000000000000	0.00414723621774\\
33.20000000000000	0.00932523132520\\
33.30000000000000	0.00004880336345\\
33.40000000000000	0.00044194861875\\
33.50000000000000	0.00141139280657\\
33.60000000000000	0.00486772180505\\
33.80000000000000	0.00028469422488\\
33.90000000000000	0.00104133681891\\
34.00000000000000	0.00403450375306\\
34.10000000000000	0.00988952434813\\
34.20000000000000	0.00012591959217\\
34.30000000000000	0.00096722818675\\
34.40000000000000	0.00328015982888\\
34.50000000000000	0.00949101592249\\
34.60000000000000	0.00015449428283\\
34.70000000000000	0.00112135054958\\
34.80000000000000	0.00354850487886\\
34.90000000000000	0.00957045793491\\
35.00000000000000	0.00016907747054\\
35.10000000000000	0.00113313293336\\
35.20000000000000	0.00311487420208\\
35.30000000000000	0.00822801288140\\
35.40000000000000	0.00000048132157\\
35.50000000000000	0.00064676123050\\
35.60000000000000	0.00151586587434\\
35.70000000000000	0.00597496199103\\
35.90000000000000	0.00043205347156\\
36.00000000000000	0.00128040325904\\
36.10000000000000	0.00835225842435\\
36.20000000000000	0.00022151953164\\
36.30000000000000	0.00084192524828\\
36.40000000000000	0.00595735908398\\
36.50000000000000	0.00000016056149\\
36.60000000000000	0.00053569682235\\
36.70000000000000	0.00375162704536\\
36.90000000000000	0.00065587494259\\
37.00000000000000	0.00199737905644\\
37.10000000000000	0.00829247216106\\
37.20000000000000	0.00028833763756\\
37.30000000000000	0.00079768283521\\
37.40000000000000	0.00411456832146\\
37.50000000000000	0.00000006714472\\
37.60000000000000	0.00060796547717\\
37.70000000000000	0.00142507598598\\
37.80000000000000	0.00674414983443\\
37.90000000000000	0.00028164365427\\
38.00000000000000	0.00073301083912\\
38.10000000000000	0.00180111950985\\
38.20000000000000	0.00488902292492\\
38.30000000000000	0.00983900108362\\
38.40000000000001	0.00030769255047\\
38.50000000000000	0.00063848538748\\
38.60000000000000	0.00146540525283\\
38.70000000000000	0.00623725330443\\
38.80000000000000	0.00029587586503\\
38.90000000000001	0.00057975645429\\
39.00000000000000	0.00132535412429\\
39.10000000000000	0.00572214974851\\
39.20000000000000	0.00000006621842\\
39.30000000000000	0.00051172823555\\
39.40000000000001	0.00099892411925\\
39.50000000000000	0.00322603415711\\
39.60000000000000	0.00998212981599\\
39.70000000000000	0.00020439647868\\
39.80000000000000	0.00071946042377\\
39.90000000000001	0.00147595903223\\
40.00000000000000	0.00586940578271\\
40.20000000000000	0.00052462042227\\
40.30000000000000	0.00262114332052\\
40.40000000000001	0.00593265348398\\
40.50000000000000	0.00025238817904\\
40.60000000000000	0.00125795227088\\
40.80000000000000	0.00043255270521\\
40.90000000000001	0.00249093633134\\
41.00000000000000	0.00629513763702\\
41.10000000000000	0.00033972525964\\
41.20000000000000	0.00266751751324\\
41.30000000000000	0.00833264235053\\
41.50000000000000	0.00070366971867\\
41.60000000000000	0.00320540102987\\
41.70000000000000	0.00594403325200\\
41.80000000000000	0.00000051568151\\
41.90000000000001	0.00183288104487\\
42.00000000000000	0.00694934320361\\
42.10000000000000	0.00000021965875\\
42.20000000000000	0.00093357587036\\
42.30000000000000	0.00361741739302\\
42.40000000000001	0.00000023827146\\
42.50000000000000	0.00140212674173\\
42.60000000000000	0.00376211190288\\
42.70000000000000	0.00689399599303\\
42.80000000000000	0.00077134482735\\
42.90000000000001	0.00268148537196\\
43.00000000000000	0.00537129098016\\
43.10000000000000	0.00000007878046\\
43.20000000000000	0.00116150563872\\
43.30000000000000	0.00163560718801\\
43.40000000000001	0.00710625539082\\
43.50000000000000	0.00055088199759\\
43.60000000000000	0.00097107307992\\
43.70000000000000	0.00283469788519\\
43.80000000000000	0.00781961658906\\
43.90000000000001	0.00000006183122\\
44.00000000000000	0.00091427107255\\
44.10000000000000	0.00140960214332\\
44.20000000000000	0.00696809434780\\
44.30000000000000	0.00059848218531\\
44.40000000000001	0.00095146530238\\
44.50000000000000	0.00562951146697\\
44.70000000000000	0.00071069727722\\
44.80000000000000	0.00183710666473\\
45.00000000000000	0.00058039177303\\
45.10000000000000	0.00133030350005\\
45.20000000000000	0.00558592329770\\
45.40000000000001	0.00073084663911\\
45.50000000000000	0.00200251082365\\
45.60000000000000	0.00964313507707\\
45.70000000000000	0.00072520779052\\
45.80000000000000	0.00101640477410\\
45.90000000000001	0.00582041379087\\
46.00000000000000	0.00000040662991\\
46.10000000000000	0.00063333278188\\
46.20000000000000	0.00225011847303\\
46.30000000000000	0.00706464537697\\
46.50000000000000	0.00060760813054\\
46.60000000000000	0.00142509359376\\
46.70000000000000	0.00768246348156\\
46.80000000000000	0.00067835207532\\
46.90000000000001	0.00147581356991\\
47.00000000000000	0.00763441400204\\
47.10000000000000	0.00937832800078\\
47.20000000000000	0.00000106839922\\
47.30000000000000	0.00138681242287\\
47.40000000000001	0.00637553563558\\
47.50000000000000	0.00649400507024\\
47.60000000000000	0.00104921992163\\
47.70000000000000	0.00297152206144\\
47.90000000000001	0.00190365019559\\
48.00000000000000	0.00268265336447\\
48.20000000000000	0.00247993990578\\
48.30000000000000	0.00399821199326\\
48.40000000000001	0.00622869253863\\
48.50000000000000	0.00000016699844\\
48.60000000000000	0.00116995446313\\
48.70000000000000	0.00276127082408\\
48.80000000000000	0.00866586692756\\
48.90000000000001	0.00078112328961\\
49.00000000000000	0.00157461066169\\
49.10000000000000	0.00538607863572\\
49.20000000000000	0.00000018673773\\
49.30000000000000	0.00115132464693\\
49.40000000000001	0.00133181839321\\
49.50000000000000	0.00498372306666\\
49.60000000000000	0.00000010806760\\
49.70000000000000	0.00118672607108\\
49.80000000000000	0.00114567509334\\
49.90000000000001	0.00490865146360\\
50.00000000000000	0.00000023278409\\
};
\addlegendentry{CMM \cite{KrahYinBergmannNaveSchneider2024}}
\end{axis}

\end{tikzpicture}

%% file: figures/cpu_timing_evolution.tex
%
\definecolor{mycolor1}{rgb}{0.00000,0.44700,0.74100}%
\definecolor{mycolor2}{rgb}{0.85000,0.32500,0.09800}%
\definecolor{mycolor3}{rgb}{0.92900,0.69400,0.12500}%
\definecolor{mycolor4}{rgb}{0.49400,0.18400,0.55600}%
\definecolor{mycolor5}{rgb}{0.46600,0.67400,0.18800}%
\definecolor{mycolor6}{rgb}{0.30100,0.74500,0.93300}%
\begin{tikzpicture}

\begin{axis}[%
width=0.951\figurewidth,
height=\figureheight,
at={(0\figurewidth,0\figureheight)},
scale only axis,
xmin=0.00000000000000,
xmax=40.00000000000000,
xminorticks=true,
xlabel style={font=\color{white!15!black}},
xlabel={Simulation time $t$},
ymin=0.00000000000000,
ymax=14.00000000000000,
yminorticks=true,
ylabel style={font=\color{white!15!black}},
ylabel={CPU time per iteration (s)},
axis background/.style={fill=white},
title style={font=\bfseries},
title={},
xmajorgrids,
ymajorgrids,
grid style={opacity=0.3},
legend style={at={(0.02,0.97)}, anchor=north west, legend cell align=left, align=left, fill=none, draw=none},
scaled ticks=false, xticklabel style={/pgf/number format/fixed},yticklabel style={/pgf/number format/fixed}
]
\addplot [color=mycolor1, line width=1.5pt]
  table[row sep=crcr]{%
0.25000000000000	2.11054737500000\\
0.50000000000000	2.16423041700000\\
0.75000000000000	0.86295479200000\\
1.00000000000000	0.90617891700000\\
1.25000000000000	0.74434787500000\\
1.50000000000000	0.93570316700000\\
1.75000000000000	0.74980312500000\\
2.00000000000000	0.75155350000000\\
2.25000000000000	0.82786866600000\\
2.50000000000000	0.65829479200000\\
2.75000000000000	0.72040745800000\\
3.00000000000000	0.79338762500000\\
3.25000000000000	0.86248620800000\\
3.50000000000000	0.93581883400000\\
3.75000000000000	0.76419670900000\\
4.00000000000000	0.83393195900000\\
4.25000000000000	0.90374216600000\\
4.50000000000000	0.98765595800000\\
4.75000000000000	1.04659854200000\\
5.00000000000000	0.88179245800000\\
5.25000000000000	0.94619245800000\\
5.50000000000000	1.01197054200000\\
5.75000000000000	1.09301666700000\\
6.00000000000000	1.16327845800000\\
6.25000000000000	0.99231333400000\\
6.50000000000000	1.06020116700000\\
6.75000000000000	1.13259937500000\\
7.00000000000000	1.20827633300000\\
7.25000000000000	1.27412975000000\\
7.50000000000000	1.10286250000000\\
7.75000000000000	1.17957425000000\\
8.00000000000000	1.24832254200000\\
8.25000000000000	1.31908345900000\\
8.50000000000000	1.38604745800000\\
8.75000000000000	1.22053945800000\\
9.00000000000000	1.28561595800000\\
9.25000000000000	1.36461400000000\\
9.50000000000000	1.43261187500000\\
9.75000000000000	1.49879704100000\\
10.00000000000000	1.32953941700000\\
10.25000000000000	1.39982170900000\\
10.50000000000000	1.46771083400000\\
10.75000000000000	1.54151991600000\\
11.00000000000000	1.61379270800000\\
11.25000000000000	1.45212145800000\\
11.50000000000000	1.54655704200000\\
11.75000000000000	1.58914279200000\\
12.00000000000000	1.65772329200000\\
12.25000000000000	1.73521975000000\\
12.50000000000000	1.56028454100000\\
12.75000000000000	1.62537854200000\\
13.00000000000000	1.72748195800000\\
13.25000000000000	1.77287062500000\\
13.50000000000000	1.84967133400000\\
13.75000000000000	1.67402670800000\\
14.00000000000000	1.73981691700000\\
14.25000000000000	1.82005783400000\\
14.50000000000000	1.88388962500000\\
14.75000000000000	1.94855516700000\\
15.00000000000000	1.78437345800000\\
15.25000000000000	1.85227691700000\\
15.50000000000000	1.92663720800000\\
15.75000000000000	2.00617175000000\\
16.00000000000000	2.06454691700000\\
16.25000000000000	1.95458987500000\\
16.50000000000000	2.04511887500000\\
16.75000000000000	2.03796358300000\\
17.00000000000000	2.11346912500000\\
17.25000000000000	2.17553641700000\\
17.50000000000000	2.13182691700000\\
17.75000000000000	2.08991083300000\\
18.00000000000000	2.15408516600000\\
18.25000000000000	2.22759870900000\\
18.50000000000000	2.29041216600000\\
18.75000000000000	2.11814575000000\\
19.00000000000000	2.19687295800000\\
19.25000000000000	2.27083575000000\\
19.50000000000000	2.34024508400000\\
19.75000000000000	2.40861204100000\\
20.00000000000000	2.24338358300000\\
20.25000000000000	2.35397979200000\\
20.50000000000000	2.38819575000000\\
20.75000000000000	2.45903379200000\\
21.00000000000000	2.52549633300000\\
21.25000000000000	2.36546629200000\\
21.50000000000000	2.41928708300000\\
21.75000000000000	2.49066854100000\\
22.00000000000000	2.56524025000000\\
22.25000000000000	2.63748104200000\\
22.50000000000000	2.47324333400000\\
22.75000000000000	2.54071816700000\\
23.00000000000000	2.60972191700000\\
23.25000000000000	2.67035975000000\\
23.50000000000000	2.74936483300000\\
23.75000000000000	2.57775025000000\\
24.00000000000000	2.65253604200000\\
24.25000000000000	2.75523654100000\\
24.50000000000000	2.78724229200000\\
24.75000000000000	2.85887479200000\\
25.00000000000000	2.71240912500000\\
25.25000000000000	2.77836462500000\\
25.50000000000000	2.78996870900000\\
25.75000000000000	2.79697870900000\\
26.00000000000000	2.81847658400000\\
26.25000000000000	2.79439491700000\\
26.50000000000000	2.86589908400000\\
26.75000000000000	2.94752808400000\\
27.00000000000000	3.01526791700000\\
27.25000000000000	3.08161033300000\\
27.50000000000000	3.00789516700000\\
27.75000000000000	2.99701004200000\\
28.00000000000000	3.05828466700000\\
28.25000000000000	3.13263529200000\\
28.50000000000000	3.19608325000000\\
28.75000000000000	3.02472608400000\\
29.00000000000000	3.09879745900000\\
29.25000000000000	3.16627587500000\\
29.50000000000000	3.24058954200000\\
29.75000000000000	3.30925229100000\\
30.00000000000000	3.31695058400000\\
30.25000000000000	3.21136054200000\\
30.50000000000000	3.28437175000000\\
30.75000000000000	3.41938800000000\\
31.00000000000000	3.42150133300000\\
31.25000000000000	3.28224950000000\\
31.50000000000000	3.37831187500000\\
31.75000000000000	3.42583029200000\\
32.00000000000000	3.46549270800000\\
32.25000000000000	3.72258275000000\\
32.50000000000000	3.36484100000000\\
32.75000000000000	3.44417245800000\\
33.00000000000000	3.50808016600000\\
33.25000000000000	3.58593108300000\\
33.50000000000000	3.64965875000000\\
33.75000000000000	3.48671816700000\\
34.00000000000000	3.54931412500000\\
34.25000000000000	3.61965666700000\\
34.50000000000000	3.84590295800000\\
34.75000000000000	3.77079425000000\\
35.00000000000000	3.60013037500000\\
35.25000000000000	3.66364137500000\\
35.50000000000000	3.92408687500000\\
35.75000000000000	3.80968729200000\\
36.00000000000000	3.87103183300000\\
36.25000000000000	3.71144370900000\\
36.50000000000000	3.77156887500000\\
36.75000000000000	3.84605900000000\\
37.00000000000000	3.91682241700000\\
37.25000000000000	4.19553975000000\\
37.50000000000000	3.82137895800000\\
37.75000000000000	4.05393504100000\\
38.00000000000000	3.95675604200000\\
38.25000000000000	4.03386362500000\\
38.50000000000000	4.10406883400000\\
38.75000000000000	3.93503775000000\\
39.00000000000000	4.19552270900000\\
39.25000000000000	4.07498670900000\\
39.50000000000000	4.14884541600000\\
39.75000000000000	4.22774537500000\\
40.00000000000000	4.07317183300000\\
};
\addlegendentry{$N_{\mathrm{remap}} = 5$}

\addplot [color=mycolor2, line width=1.5pt]
  table[row sep=crcr]{%
0.25000000000000	1.25206737500000\\
0.50000000000000	2.13747625000000\\
0.75000000000000	0.86336620900000\\
1.00000000000000	0.93745583300000\\
1.25000000000000	0.94265033300000\\
1.50000000000000	1.13562783300000\\
1.75000000000000	0.95942979200000\\
2.00000000000000	0.95169133400000\\
2.25000000000000	1.02557720800000\\
2.50000000000000	0.59141870800000\\
2.75000000000000	0.64925454200000\\
3.00000000000000	0.71739254200000\\
3.25000000000000	0.78610233400000\\
3.50000000000000	0.86130937500000\\
3.75000000000000	0.92744875000000\\
4.00000000000000	1.00019537500000\\
4.25000000000000	1.06560987500000\\
4.50000000000000	1.14223870800000\\
4.75000000000000	1.21124787500000\\
5.00000000000000	0.72073729200000\\
5.25000000000000	0.75702183300000\\
5.50000000000000	0.82709741700000\\
5.75000000000000	0.89805516700000\\
6.00000000000000	0.96654879100000\\
6.25000000000000	1.04099362500000\\
6.50000000000000	1.11196050000000\\
6.75000000000000	1.17651420800000\\
7.00000000000000	1.25582375000000\\
7.25000000000000	1.32447762500000\\
7.50000000000000	0.79442079100000\\
7.75000000000000	0.87143091600000\\
8.00000000000000	0.95230420900000\\
8.25000000000000	1.01855908300000\\
8.50000000000000	1.55838270800000\\
8.75000000000000	1.15474262500000\\
9.00000000000000	1.22940570800000\\
9.25000000000000	1.29515437500000\\
9.50000000000000	1.36406016700000\\
9.75000000000000	1.43123625000000\\
10.00000000000000	0.91466762500000\\
10.25000000000000	0.98632170900000\\
10.50000000000000	1.05874900000000\\
10.75000000000000	1.12608916600000\\
11.00000000000000	1.19723512500000\\
11.25000000000000	1.26661725000000\\
11.50000000000000	1.33708937500000\\
11.75000000000000	1.40975995800000\\
12.00000000000000	1.47924987500000\\
12.25000000000000	1.55396279200000\\
12.50000000000000	1.02809529200000\\
12.75000000000000	1.09763579200000\\
13.00000000000000	1.16677220800000\\
13.25000000000000	1.23581233400000\\
13.50000000000000	1.31251479200000\\
13.75000000000000	1.38231079200000\\
14.00000000000000	1.45798983300000\\
14.25000000000000	1.52761895800000\\
14.50000000000000	1.59183166600000\\
14.75000000000000	1.66646795900000\\
15.00000000000000	1.14084529200000\\
15.25000000000000	1.21143700000000\\
15.50000000000000	1.28566945800000\\
15.75000000000000	1.36267187500000\\
16.00000000000000	1.42581854100000\\
16.25000000000000	1.49379679100000\\
16.50000000000000	1.56564795900000\\
16.75000000000000	1.63796587500000\\
17.00000000000000	1.70617800000000\\
17.25000000000000	1.78336241600000\\
17.50000000000000	1.25503929200000\\
17.75000000000000	1.33540720800000\\
18.00000000000000	1.39429800000000\\
18.25000000000000	1.47512679200000\\
18.50000000000000	1.53518820800000\\
18.75000000000000	1.61083708300000\\
19.00000000000000	1.67787133400000\\
19.25000000000000	1.74875195800000\\
19.50000000000000	1.82603295800000\\
19.75000000000000	1.89162062500000\\
20.00000000000000	1.36969895900000\\
20.25000000000000	1.43397012500000\\
20.50000000000000	1.50933354200000\\
20.75000000000000	1.58490562500000\\
21.00000000000000	1.65391425000000\\
21.25000000000000	1.72001666600000\\
21.50000000000000	1.78975050000000\\
21.75000000000000	1.86272008300000\\
22.00000000000000	1.93161691600000\\
22.25000000000000	2.00792320800000\\
22.50000000000000	1.48449304200000\\
22.75000000000000	1.55259008300000\\
23.00000000000000	1.62189687500000\\
23.25000000000000	1.69418450000000\\
23.50000000000000	1.76391487500000\\
23.75000000000000	1.83201287500000\\
24.00000000000000	1.89714175000000\\
24.25000000000000	1.98309375000000\\
24.50000000000000	2.03719150000000\\
24.75000000000000	2.10939420800000\\
25.00000000000000	1.59824616600000\\
25.25000000000000	1.68103179200000\\
25.50000000000000	1.72665658300000\\
25.75000000000000	1.79232629200000\\
26.00000000000000	1.87668350000000\\
26.25000000000000	1.94985633400000\\
26.50000000000000	2.02097008400000\\
26.75000000000000	2.08887187500000\\
27.00000000000000	2.15770337500000\\
27.25000000000000	2.39695045800000\\
27.50000000000000	1.70596858300000\\
27.75000000000000	1.79012070900000\\
28.00000000000000	1.84995095800000\\
28.25000000000000	2.09118100000000\\
28.50000000000000	1.98848262500000\\
28.75000000000000	2.05759354200000\\
29.00000000000000	2.13176604200000\\
29.25000000000000	2.20324495800000\\
29.50000000000000	2.26844766700000\\
29.75000000000000	2.34270145800000\\
30.00000000000000	1.81950666700000\\
30.25000000000000	1.88864212500000\\
30.50000000000000	1.96442175000000\\
30.75000000000000	2.03708762500000\\
31.00000000000000	2.10772333300000\\
31.25000000000000	2.20963483300000\\
31.50000000000000	2.29781004100000\\
31.75000000000000	2.52273137500000\\
32.00000000000000	2.38507366700000\\
32.25000000000000	2.45664970800000\\
32.50000000000000	1.92818475000000\\
32.75000000000000	1.99979579100000\\
33.00000000000000	2.07733729200000\\
33.25000000000000	2.14983812500000\\
33.50000000000000	2.21698804200000\\
33.75000000000000	2.28729304200000\\
34.00000000000000	2.35964879200000\\
34.25000000000000	2.42582587500000\\
34.50000000000000	2.49748233300000\\
34.75000000000000	2.57270266700000\\
35.00000000000000	2.05210650000000\\
35.25000000000000	2.11713391700000\\
35.50000000000000	2.18941575000000\\
35.75000000000000	2.25914616700000\\
36.00000000000000	2.33430408300000\\
36.25000000000000	2.39878216700000\\
36.50000000000000	2.76977433300000\\
36.75000000000000	2.53977683400000\\
37.00000000000000	2.61578641700000\\
37.25000000000000	2.68092275000000\\
37.50000000000000	2.32387241700000\\
37.75000000000000	2.23120483300000\\
38.00000000000000	2.50865100000000\\
38.25000000000000	2.36499604200000\\
38.50000000000000	2.87799316600000\\
38.75000000000000	2.51407333300000\\
39.00000000000000	2.58454175000000\\
39.25000000000000	2.79068404100000\\
39.50000000000000	2.72492583300000\\
39.75000000000000	2.79524183300000\\
40.00000000000000	2.28634866700000\\
};
\addlegendentry{$N_{\mathrm{remap}} = 10$}

\addplot [color=mycolor3, line width=1.5pt]
  table[row sep=crcr]{%
0.25000000000000	1.13973766700000\\
0.50000000000000	2.14098875000000\\
0.75000000000000	0.88548766600000\\
1.00000000000000	0.93939758400000\\
1.25000000000000	0.94784420800000\\
1.50000000000000	1.14680704200000\\
1.75000000000000	0.95103483400000\\
2.00000000000000	0.95066287500000\\
2.25000000000000	1.09344900000000\\
2.50000000000000	1.09899287500000\\
2.75000000000000	1.16992858300000\\
3.00000000000000	1.23935816700000\\
3.25000000000000	1.30898725000000\\
3.50000000000000	1.37857795800000\\
3.75000000000000	1.45570741600000\\
4.00000000000000	1.52205204200000\\
4.25000000000000	1.59051883300000\\
4.50000000000000	1.66043541700000\\
4.75000000000000	1.73368412500000\\
5.00000000000000	0.59573233400000\\
5.25000000000000	0.65294833300000\\
5.50000000000000	0.71536508300000\\
5.75000000000000	0.78563070900000\\
6.00000000000000	0.85933920900000\\
6.25000000000000	0.92589991600000\\
6.50000000000000	0.99501462500000\\
6.75000000000000	1.06630712500000\\
7.00000000000000	1.13760412500000\\
7.25000000000000	1.21005208300000\\
7.50000000000000	1.28270720800000\\
7.75000000000000	1.35692025000000\\
8.00000000000000	1.42219575000000\\
8.25000000000000	1.49537766700000\\
8.50000000000000	1.56028400000000\\
8.75000000000000	1.65001583300000\\
9.00000000000000	1.70804775000000\\
9.25000000000000	1.80051466600000\\
9.50000000000000	1.84715020900000\\
9.75000000000000	1.91814008300000\\
10.00000000000000	0.69281070800000\\
10.25000000000000	0.75577845800000\\
10.50000000000000	0.82443958400000\\
10.75000000000000	0.89537733300000\\
11.00000000000000	0.96792066700000\\
11.25000000000000	1.03640183300000\\
11.50000000000000	1.10676883300000\\
11.75000000000000	1.18473004200000\\
12.00000000000000	1.25411612500000\\
12.25000000000000	1.33364220800000\\
12.50000000000000	1.39942829200000\\
12.75000000000000	1.46526987500000\\
13.00000000000000	1.53286212500000\\
13.25000000000000	1.60803350000000\\
13.50000000000000	1.67499766700000\\
13.75000000000000	1.74876620800000\\
14.00000000000000	1.81792575000000\\
14.25000000000000	1.91100841700000\\
14.50000000000000	1.96806000000000\\
14.75000000000000	2.03693537500000\\
15.00000000000000	0.80351970900000\\
15.25000000000000	0.86949791700000\\
15.50000000000000	0.94389662500000\\
15.75000000000000	1.01684025000000\\
16.00000000000000	1.08739300000000\\
16.25000000000000	1.15055837500000\\
16.50000000000000	1.22201037500000\\
16.75000000000000	1.29735775000000\\
17.00000000000000	1.36652587500000\\
17.25000000000000	1.44403833300000\\
17.50000000000000	1.50923100000000\\
17.75000000000000	1.58450012500000\\
18.00000000000000	1.64800141700000\\
18.25000000000000	1.72038683300000\\
18.50000000000000	1.78883245800000\\
18.75000000000000	1.86043350000000\\
19.00000000000000	1.93987329200000\\
19.25000000000000	2.01079000000000\\
19.50000000000000	2.08370104100000\\
19.75000000000000	2.14566883300000\\
20.00000000000000	0.91534850000000\\
20.25000000000000	0.98516275000000\\
20.50000000000000	1.05431929200000\\
20.75000000000000	1.13557083300000\\
21.00000000000000	1.20078758300000\\
21.25000000000000	1.28039900000000\\
21.50000000000000	1.33757854100000\\
21.75000000000000	1.41307150000000\\
22.00000000000000	1.47907662500000\\
22.25000000000000	1.55375054100000\\
22.50000000000000	1.62498133300000\\
22.75000000000000	1.69492050000000\\
23.00000000000000	1.76135045900000\\
23.25000000000000	1.82941012500000\\
23.50000000000000	1.90663945900000\\
23.75000000000000	1.98063808300000\\
24.00000000000000	2.05089108400000\\
24.25000000000000	2.12532483300000\\
24.50000000000000	2.19764325000000\\
24.75000000000000	2.27660654200000\\
25.00000000000000	1.05212212500000\\
25.25000000000000	1.10522400000000\\
25.50000000000000	1.16409920800000\\
25.75000000000000	1.24002183300000\\
26.00000000000000	1.30597641700000\\
26.25000000000000	1.38107237500000\\
26.50000000000000	1.44882575000000\\
26.75000000000000	1.52868525000000\\
27.00000000000000	1.59154504200000\\
27.25000000000000	1.66008795900000\\
27.50000000000000	1.73513787500000\\
27.75000000000000	1.81548495800000\\
28.00000000000000	1.87465979200000\\
28.25000000000000	1.94878720900000\\
28.50000000000000	2.01792312500000\\
28.75000000000000	2.09656637500000\\
29.00000000000000	2.17021137500000\\
29.25000000000000	2.23808608400000\\
29.50000000000000	2.30039116600000\\
29.75000000000000	2.37736908300000\\
30.00000000000000	1.14174691600000\\
30.25000000000000	1.20972350000000\\
30.50000000000000	1.28388341700000\\
30.75000000000000	1.35516212500000\\
31.00000000000000	1.42182916600000\\
31.25000000000000	1.51645966700000\\
31.50000000000000	1.61931633400000\\
31.75000000000000	1.67285716700000\\
32.00000000000000	1.72512083300000\\
32.25000000000000	1.78046891700000\\
32.50000000000000	1.83844458400000\\
32.75000000000000	1.91375641700000\\
33.00000000000000	1.98047116700000\\
33.25000000000000	2.05728979200000\\
33.50000000000000	2.12921608400000\\
33.75000000000000	2.19901483300000\\
34.00000000000000	2.26903537500000\\
34.25000000000000	2.34657058300000\\
34.50000000000000	2.41222704100000\\
34.75000000000000	2.48703066600000\\
35.00000000000000	1.24607104200000\\
35.25000000000000	1.32240254200000\\
35.50000000000000	1.39615204100000\\
35.75000000000000	1.46797262500000\\
36.00000000000000	1.53960395800000\\
36.25000000000000	1.61290537500000\\
36.50000000000000	1.68045654200000\\
36.75000000000000	1.74446016600000\\
37.00000000000000	1.81992383300000\\
37.25000000000000	1.89037733300000\\
37.50000000000000	1.95905937500000\\
37.75000000000000	2.03048400000000\\
38.00000000000000	2.10269650000000\\
38.25000000000000	2.17427733400000\\
38.50000000000000	2.24344025000000\\
38.75000000000000	2.32399408300000\\
39.00000000000000	2.38270391700000\\
39.25000000000000	2.47387700000000\\
39.50000000000000	2.53056266700000\\
39.75000000000000	2.60676495800000\\
40.00000000000000	1.38067950000000\\
};
\addlegendentry{$N_{\mathrm{remap}} = 20$}

\addplot [color=mycolor4, line width=1.5pt]
  table[row sep=crcr]{%
0.25000000000000	1.14315841700000\\
0.50000000000000	2.14811745800000\\
0.75000000000000	0.86617712500000\\
1.00000000000000	0.93530162500000\\
1.25000000000000	0.94351562500000\\
1.50000000000000	1.14351550000000\\
1.75000000000000	0.95405804200000\\
2.00000000000000	0.95479666700000\\
2.25000000000000	1.02537216700000\\
2.50000000000000	1.09886820800000\\
2.75000000000000	1.16915016700000\\
3.00000000000000	1.23511558400000\\
3.25000000000000	1.30554175000000\\
3.50000000000000	1.38457433300000\\
3.75000000000000	1.45411070800000\\
4.00000000000000	1.51867195900000\\
4.25000000000000	1.59144287500000\\
4.50000000000000	1.66015537500000\\
4.75000000000000	1.73356254100000\\
5.00000000000000	1.81011608300000\\
5.25000000000000	1.88205212500000\\
5.50000000000000	1.94589804100000\\
5.75000000000000	2.02129045900000\\
6.00000000000000	2.08899825000000\\
6.25000000000000	2.16554975000000\\
6.50000000000000	2.23522629200000\\
6.75000000000000	2.30751229200000\\
7.00000000000000	2.36718412500000\\
7.25000000000000	2.44363837500000\\
7.50000000000000	0.77850745800000\\
7.75000000000000	1.10160462500000\\
8.00000000000000	1.39841037500000\\
8.25000000000000	0.80991820800000\\
8.50000000000000	0.86325954200000\\
8.75000000000000	0.93017175000000\\
9.00000000000000	1.75518479200000\\
9.25000000000000	1.12033854200000\\
9.50000000000000	1.15022370800000\\
9.75000000000000	1.21482970800000\\
10.00000000000000	1.28938679100000\\
10.25000000000000	1.36087770900000\\
10.50000000000000	1.42457208300000\\
10.75000000000000	1.48929041600000\\
11.00000000000000	1.56452175000000\\
11.25000000000000	1.63642554200000\\
11.50000000000000	1.71292550000000\\
11.75000000000000	1.77385445900000\\
12.00000000000000	1.85619300000000\\
12.25000000000000	1.91561262500000\\
12.50000000000000	1.99023879200000\\
12.75000000000000	2.05821983300000\\
13.00000000000000	2.13103516700000\\
13.25000000000000	2.20554508300000\\
13.50000000000000	2.27241079100000\\
13.75000000000000	2.35356241600000\\
14.00000000000000	2.42638400000000\\
14.25000000000000	2.48674833300000\\
14.50000000000000	2.55729429100000\\
14.75000000000000	2.64603025000000\\
15.00000000000000	0.69431708300000\\
15.25000000000000	0.76172479100000\\
15.50000000000000	0.83315375000000\\
15.75000000000000	0.91429379100000\\
16.00000000000000	0.97402708300000\\
16.25000000000000	1.03542741700000\\
16.50000000000000	1.11937533300000\\
16.75000000000000	1.18461433300000\\
17.00000000000000	1.25276054100000\\
17.25000000000000	1.32207329200000\\
17.50000000000000	1.39298375000000\\
17.75000000000000	1.47528591700000\\
18.00000000000000	1.53764166700000\\
18.25000000000000	1.61481220800000\\
18.50000000000000	1.68764504200000\\
18.75000000000000	1.75181800000000\\
19.00000000000000	1.83857279100000\\
19.25000000000000	1.89525354200000\\
19.50000000000000	1.96624062500000\\
19.75000000000000	2.03572062500000\\
20.00000000000000	2.10296541700000\\
20.25000000000000	2.17248045900000\\
20.50000000000000	2.24805466600000\\
20.75000000000000	2.32294637500000\\
21.00000000000000	2.38661975000000\\
21.25000000000000	2.45598483300000\\
21.50000000000000	2.52860883300000\\
21.75000000000000	2.59676008400000\\
22.00000000000000	2.66900345800000\\
22.25000000000000	2.75015754200000\\
22.50000000000000	0.80016891700000\\
22.75000000000000	0.87091970900000\\
23.00000000000000	0.94787483300000\\
23.25000000000000	1.01787395800000\\
23.50000000000000	1.08657137500000\\
23.75000000000000	1.15548587500000\\
24.00000000000000	1.22545187500000\\
24.25000000000000	1.30963020800000\\
24.50000000000000	1.36714304200000\\
24.75000000000000	1.43863687500000\\
25.00000000000000	1.52405037500000\\
25.25000000000000	1.60356812500000\\
25.50000000000000	1.65487150000000\\
25.75000000000000	1.71986175000000\\
26.00000000000000	1.79202104100000\\
26.25000000000000	1.86125808400000\\
26.50000000000000	1.93318941600000\\
26.75000000000000	2.01312158300000\\
27.00000000000000	2.08026145800000\\
27.25000000000000	2.15018466600000\\
27.50000000000000	2.21460050000000\\
27.75000000000000	2.29037620800000\\
28.00000000000000	2.35630191600000\\
28.25000000000000	2.42833108400000\\
28.50000000000000	2.51279858300000\\
28.75000000000000	2.58125987500000\\
29.00000000000000	2.65625283300000\\
29.25000000000000	2.71307387500000\\
29.50000000000000	2.78726745800000\\
29.75000000000000	2.85224125000000\\
30.00000000000000	0.91743625000000\\
30.25000000000000	0.98584029200000\\
30.50000000000000	1.06127229100000\\
30.75000000000000	1.12600437500000\\
31.00000000000000	1.19761316600000\\
31.25000000000000	1.29480433300000\\
31.50000000000000	1.37649145900000\\
31.75000000000000	1.42238375000000\\
32.00000000000000	1.48111270900000\\
32.25000000000000	1.55148954200000\\
32.50000000000000	1.62552850000000\\
32.75000000000000	1.69101933300000\\
33.00000000000000	1.76358604200000\\
33.25000000000000	1.83799804100000\\
33.50000000000000	1.90212695900000\\
33.75000000000000	1.97701304200000\\
34.00000000000000	2.04594025000000\\
34.25000000000000	2.11807437500000\\
34.50000000000000	2.18928262500000\\
34.75000000000000	2.26957950000000\\
35.00000000000000	2.33009204100000\\
35.25000000000000	2.40270675000000\\
35.50000000000000	2.46876166700000\\
35.75000000000000	2.54697037500000\\
36.00000000000000	2.62751966600000\\
36.25000000000000	2.68446666700000\\
36.50000000000000	2.75446733300000\\
36.75000000000000	2.82424995800000\\
37.00000000000000	2.89569933300000\\
37.25000000000000	2.97321262500000\\
37.50000000000000	1.02994637500000\\
37.75000000000000	1.08835404200000\\
38.00000000000000	1.15705125000000\\
38.25000000000000	1.23292041700000\\
38.50000000000000	1.30094929100000\\
38.75000000000000	1.37045091600000\\
39.00000000000000	1.44440741600000\\
39.25000000000000	1.51434200000000\\
39.50000000000000	1.59220245900000\\
39.75000000000000	1.66152629200000\\
40.00000000000000	1.73719475000000\\
};
\addlegendentry{$N_{\mathrm{remap}} = 30$}

\addplot [color=mycolor5, line width=1.5pt]
  table[row sep=crcr]{%
0.25000000000000	1.15166720800000\\
0.50000000000000	2.14414616700000\\
0.75000000000000	0.86704129100000\\
1.00000000000000	0.93514383400000\\
1.25000000000000	0.94350683400000\\
1.50000000000000	1.14908641600000\\
1.75000000000000	0.95163158400000\\
2.00000000000000	0.95454612500000\\
2.25000000000000	1.02645595800000\\
2.50000000000000	1.10124916700000\\
2.75000000000000	1.17205808300000\\
3.00000000000000	1.24189420800000\\
3.25000000000000	1.30867754200000\\
3.50000000000000	1.38534066700000\\
3.75000000000000	1.44916679200000\\
4.00000000000000	1.52814520900000\\
4.25000000000000	1.59485654100000\\
4.50000000000000	1.66345204200000\\
4.75000000000000	1.74020666700000\\
5.00000000000000	1.80767295800000\\
5.25000000000000	1.88009333300000\\
5.50000000000000	1.94995462500000\\
5.75000000000000	2.01776383400000\\
6.00000000000000	2.09258779200000\\
6.25000000000000	2.17173308300000\\
6.50000000000000	2.22989470800000\\
6.75000000000000	2.30404700000000\\
7.00000000000000	2.37635320800000\\
7.25000000000000	2.44695825000000\\
7.50000000000000	2.52112933300000\\
7.75000000000000	2.59219933300000\\
8.00000000000000	2.65998795800000\\
8.25000000000000	2.74154712500000\\
8.50000000000000	2.80247783300000\\
8.75000000000000	2.86850637500000\\
9.00000000000000	2.94104895900000\\
9.25000000000000	3.01081816600000\\
9.50000000000000	3.08277254100000\\
9.75000000000000	3.16410987500000\\
10.00000000000000	0.60100504200000\\
10.25000000000000	0.65527729200000\\
10.50000000000000	0.72073962500000\\
10.75000000000000	0.78549445800000\\
11.00000000000000	0.85734595900000\\
11.25000000000000	0.92700529200000\\
11.50000000000000	1.00130975000000\\
11.75000000000000	1.06955116600000\\
12.00000000000000	1.14188737500000\\
12.25000000000000	1.20859966700000\\
12.50000000000000	1.28773875000000\\
12.75000000000000	1.35677750000000\\
13.00000000000000	1.42343841700000\\
13.25000000000000	1.49137033300000\\
13.50000000000000	1.56469225000000\\
13.75000000000000	1.63698450000000\\
14.00000000000000	1.71790200000000\\
14.25000000000000	1.77900525000000\\
14.50000000000000	1.84933458300000\\
14.75000000000000	1.92422983300000\\
15.00000000000000	2.00064370800000\\
15.25000000000000	2.06056929200000\\
15.50000000000000	2.13868912500000\\
15.75000000000000	2.22159670800000\\
16.00000000000000	2.28373033400000\\
16.25000000000000	2.35024016700000\\
16.50000000000000	2.41690108400000\\
16.75000000000000	2.48637412500000\\
17.00000000000000	2.55923354100000\\
17.25000000000000	2.63157175000000\\
17.50000000000000	2.69915095800000\\
17.75000000000000	2.78041291600000\\
18.00000000000000	2.84369604100000\\
18.25000000000000	2.91593500000000\\
18.50000000000000	2.98789145800000\\
18.75000000000000	3.05372725000000\\
19.00000000000000	3.12927695800000\\
19.25000000000000	3.19700737500000\\
19.50000000000000	3.27551108300000\\
19.75000000000000	3.34239645800000\\
20.00000000000000	0.69014200000000\\
20.25000000000000	0.75540666700000\\
20.50000000000000	0.82934129200000\\
20.75000000000000	0.90820391700000\\
21.00000000000000	0.97353950000000\\
21.25000000000000	1.04621841700000\\
21.50000000000000	1.11264062500000\\
21.75000000000000	1.18216645900000\\
22.00000000000000	1.25687312500000\\
22.25000000000000	1.32202020800000\\
22.50000000000000	1.40141912500000\\
22.75000000000000	1.47109025000000\\
23.00000000000000	1.54098904100000\\
23.25000000000000	1.61098870900000\\
23.50000000000000	1.67994983300000\\
23.75000000000000	1.75233437500000\\
24.00000000000000	1.81767770800000\\
24.25000000000000	1.89246337500000\\
24.50000000000000	2.09773829200000\\
24.75000000000000	2.03931933300000\\
25.00000000000000	2.12962620800000\\
25.25000000000000	2.20239054200000\\
25.50000000000000	2.25110895800000\\
25.75000000000000	2.31493170800000\\
26.00000000000000	2.39531550000000\\
26.25000000000000	2.46091600000000\\
26.50000000000000	2.53167687500000\\
26.75000000000000	2.59915308300000\\
27.00000000000000	2.67176391700000\\
27.25000000000000	2.74157862500000\\
27.50000000000000	2.81021266600000\\
27.75000000000000	2.89213154200000\\
28.00000000000000	2.96069812500000\\
28.25000000000000	3.04134412500000\\
28.50000000000000	3.10819983400000\\
28.75000000000000	3.16881312500000\\
29.00000000000000	3.23988987500000\\
29.25000000000000	3.30809737500000\\
29.50000000000000	3.38324175000000\\
29.75000000000000	3.45197712500000\\
30.00000000000000	0.80095691600000\\
30.25000000000000	0.87105662500000\\
30.50000000000000	0.94277900000000\\
30.75000000000000	1.01930658300000\\
31.00000000000000	1.08817425000000\\
31.25000000000000	1.18079554200000\\
31.50000000000000	1.26159179100000\\
31.75000000000000	1.30394470800000\\
32.00000000000000	1.36400937500000\\
32.25000000000000	1.43807554100000\\
32.50000000000000	1.50843487500000\\
32.75000000000000	1.58062808400000\\
33.00000000000000	1.64940533400000\\
33.25000000000000	1.72629762500000\\
33.50000000000000	1.80438295900000\\
33.75000000000000	1.86462295900000\\
34.00000000000000	1.93283695900000\\
34.25000000000000	2.00554112500000\\
34.50000000000000	2.07522066700000\\
34.75000000000000	2.14039000000000\\
35.00000000000000	2.21805445900000\\
35.25000000000000	2.29168541600000\\
35.50000000000000	2.36050854200000\\
35.75000000000000	2.43219883300000\\
36.00000000000000	2.50920833300000\\
36.25000000000000	2.57170145800000\\
36.50000000000000	2.64499545900000\\
36.75000000000000	2.71488000000000\\
37.00000000000000	2.78780033400000\\
37.25000000000000	2.85679125000000\\
37.50000000000000	2.93193845800000\\
37.75000000000000	3.00029695800000\\
38.00000000000000	3.07805787500000\\
38.25000000000000	3.14504079100000\\
38.50000000000000	3.25793729100000\\
38.75000000000000	3.28531266700000\\
39.00000000000000	3.36630791700000\\
39.25000000000000	3.42059295800000\\
39.50000000000000	3.49222641700000\\
39.75000000000000	3.56643725000000\\
40.00000000000000	0.91321670800000\\
};
\addlegendentry{$N_{\mathrm{remap}} = 40$}

\addplot [color=mycolor6, line width=1.5pt]
  table[row sep=crcr]{%
0.25000000000000	1.15651558300000\\
0.50000000000000	2.14579687500000\\
0.75000000000000	0.86191400000000\\
1.00000000000000	0.93360920800000\\
1.25000000000000	0.93841745800000\\
1.50000000000000	1.16206295900000\\
1.75000000000000	0.95638170900000\\
2.00000000000000	0.95808170800000\\
2.25000000000000	1.02182729200000\\
2.50000000000000	1.09421116700000\\
2.75000000000000	1.16313879200000\\
3.00000000000000	1.23438925000000\\
3.25000000000000	1.31059041700000\\
3.50000000000000	1.37088608300000\\
3.75000000000000	1.45073283300000\\
4.00000000000000	1.51405570800000\\
4.25000000000000	1.58541995900000\\
4.50000000000000	1.66142241700000\\
4.75000000000000	1.73488595800000\\
5.00000000000000	1.79688462500000\\
5.25000000000000	1.87854625000000\\
5.50000000000000	1.94809970900000\\
5.75000000000000	2.02295129200000\\
6.00000000000000	2.08672012500000\\
6.25000000000000	2.15928079200000\\
6.50000000000000	2.22754204200000\\
6.75000000000000	2.29968545900000\\
7.00000000000000	2.37311462500000\\
7.25000000000000	2.45425408300000\\
7.50000000000000	2.51492812500000\\
7.75000000000000	2.58256925000000\\
8.00000000000000	2.65435541700000\\
8.25000000000000	2.73019408300000\\
8.50000000000000	2.80123354100000\\
8.75000000000000	2.86590337500000\\
9.00000000000000	2.94023816700000\\
9.25000000000000	3.01758141700000\\
9.50000000000000	3.08387554200000\\
9.75000000000000	3.15487320900000\\
10.00000000000000	3.21978641700000\\
10.25000000000000	3.29540595900000\\
10.50000000000000	3.36758150000000\\
10.75000000000000	3.43321491700000\\
11.00000000000000	3.50629975000000\\
11.25000000000000	3.57290787500000\\
11.50000000000000	3.65859820800000\\
11.75000000000000	3.71193079200000\\
12.00000000000000	3.79101020900000\\
12.25000000000000	3.87088658300000\\
12.50000000000000	3.93328258400000\\
12.75000000000000	3.99938516600000\\
13.00000000000000	4.06643116600000\\
13.25000000000000	4.14010370800000\\
13.50000000000000	4.21075837500000\\
13.75000000000000	4.28078470900000\\
14.00000000000000	4.35791691700000\\
14.25000000000000	4.43377125000000\\
14.50000000000000	4.49658595800000\\
14.75000000000000	4.56883812500000\\
15.00000000000000	4.64143287500000\\
15.25000000000000	4.70609937500000\\
15.50000000000000	4.79471800000000\\
15.75000000000000	4.86766058300000\\
16.00000000000000	4.92290383400000\\
16.25000000000000	4.99532091600000\\
16.50000000000000	5.06222475000000\\
16.75000000000000	5.13006204100000\\
17.00000000000000	5.20517025000000\\
17.25000000000000	5.27980033400000\\
17.50000000000000	5.34030025000000\\
17.75000000000000	5.42098450000000\\
18.00000000000000	5.48131245900000\\
18.25000000000000	5.57049333300000\\
18.50000000000000	5.63955787500000\\
18.75000000000000	5.70531854100000\\
19.00000000000000	5.77223358300000\\
19.25000000000000	5.83939512500000\\
19.50000000000000	5.91273050000000\\
19.75000000000000	5.98328533400000\\
20.00000000000000	6.02843729200000\\
20.25000000000000	6.12634995800000\\
20.50000000000000	6.25945258300000\\
20.75000000000000	6.27279050000000\\
21.00000000000000	6.33876679200000\\
21.25000000000000	6.41029191600000\\
21.50000000000000	6.47541304200000\\
21.75000000000000	6.55839220800000\\
22.00000000000000	6.61347583300000\\
22.25000000000000	6.68734316700000\\
22.50000000000000	6.77314625000000\\
22.75000000000000	6.83024412500000\\
23.00000000000000	6.90538883400000\\
23.25000000000000	6.98268370900000\\
23.50000000000000	7.04618695800000\\
23.75000000000000	7.11274325000000\\
24.00000000000000	7.19087875000000\\
24.25000000000000	7.25164716600000\\
24.50000000000000	7.32273825000000\\
24.75000000000000	7.40532970800000\\
25.00000000000000	7.48210375000000\\
25.25000000000000	7.58270416700000\\
25.50000000000000	7.79937783300000\\
25.75000000000000	7.74831837500000\\
26.00000000000000	7.77232137500000\\
26.25000000000000	7.81999262500000\\
26.50000000000000	7.88936112500000\\
26.75000000000000	7.99196587500000\\
27.00000000000000	8.03415866700000\\
27.25000000000000	8.11116587500000\\
27.50000000000000	8.18176083300000\\
27.75000000000000	8.24878762500000\\
28.00000000000000	8.32526445900000\\
28.25000000000000	8.40789341700000\\
28.50000000000000	8.46076366700000\\
28.75000000000000	8.52851762500000\\
29.00000000000000	8.59415579200000\\
29.25000000000000	8.69913870800000\\
29.50000000000000	8.75689675000000\\
29.75000000000000	8.82663904200000\\
30.00000000000000	8.89196487500000\\
30.25000000000000	8.94825212500000\\
30.50000000000000	9.02249412500000\\
30.75000000000000	9.10841704200000\\
31.00000000000000	9.17230029200000\\
31.25000000000000	9.27621995900000\\
31.50000000000000	9.35835137500000\\
31.75000000000000	9.40160870800000\\
32.00000000000000	9.44619441700000\\
32.25000000000000	9.56631562500000\\
32.50000000000000	9.58739166700000\\
32.75000000000000	9.71128029200000\\
33.00000000000000	9.73378700000000\\
33.25000000000000	9.79659695900000\\
33.50000000000000	9.87074245800000\\
33.75000000000000	9.94969412500000\\
34.00000000000000	9.99789525000000\\
34.25000000000000	10.10470762500000\\
34.50000000000000	10.17008137500000\\
34.75000000000000	10.23759245800000\\
35.00000000000000	10.29501866700000\\
35.25000000000000	10.36533054200000\\
35.50000000000000	10.44271825000000\\
35.75000000000000	10.52067408400000\\
36.00000000000000	10.57764216600000\\
36.25000000000000	10.61907037500000\\
36.50000000000000	10.69657370900000\\
36.75000000000000	10.78433029200000\\
37.00000000000000	10.87349583300000\\
37.25000000000000	10.92108029100000\\
37.50000000000000	11.05198391700000\\
37.75000000000000	11.09804583400000\\
38.00000000000000	11.14827029200000\\
38.25000000000000	11.22423158300000\\
38.50000000000000	11.29667129200000\\
38.75000000000000	11.35065283300000\\
39.00000000000000	11.42022200000000\\
39.25000000000000	11.47662795800000\\
39.50000000000000	11.57727429100000\\
39.75000000000000	11.63806504200000\\
40.00000000000000	11.70022445900000\\
};
\addlegendentry{NuFi}

\end{axis}
\end{tikzpicture}%

%% file: figures/cpu_timing_cumulative.tex
%
\definecolor{mycolor1}{rgb}{0.00000,0.44700,0.74100}%
\definecolor{mycolor2}{rgb}{0.85000,0.32500,0.09800}%
\definecolor{mycolor3}{rgb}{0.92900,0.69400,0.12500}%
\definecolor{mycolor4}{rgb}{0.49400,0.18400,0.55600}%
\definecolor{mycolor5}{rgb}{0.46600,0.67400,0.18800}%
\definecolor{mycolor6}{rgb}{0.30100,0.74500,0.93300}%
\begin{tikzpicture}

\begin{axis}[%
width=0.951\figurewidth,
height=\figureheight,
at={(0\figurewidth,0\figureheight)},
scale only axis,
xmin=0.00000000000000,
xmax=40.00000000000000,
xminorticks=true,
xlabel style={font=\color{white!15!black}},
xlabel={Simulation time $t$},
ymin=0.00000000000000,
ymax=1000.00000000000000,
yminorticks=true,
ylabel style={font=\color{white!15!black}},
ylabel={Cumulative CPU time (s)},
axis background/.style={fill=white},
title style={font=\bfseries},
title={},
xmajorgrids,
ymajorgrids,
grid style={opacity=0.3},
legend style={at={(0.05,0.95)}, anchor=north west, legend cell align=left, align=left, fill=none, draw=none},
scaled ticks=false, xticklabel style={/pgf/number format/fixed},yticklabel style={/pgf/number format/fixed}
]
\addplot [color=mycolor1, line width=1.5pt]
  table[row sep=crcr]{%
0.25000000000000	2.11054737500000\\
0.50000000000000	4.27477779200000\\
0.75000000000000	5.13773258400000\\
1.00000000000000	6.04391150100000\\
1.25000000000000	6.78825937600000\\
1.50000000000000	7.72396254300000\\
1.75000000000000	8.47376566800000\\
2.00000000000000	9.22531916800000\\
2.25000000000000	10.05318783400000\\
2.50000000000000	10.71148262600000\\
2.75000000000000	11.43189008400000\\
3.00000000000000	12.22527770900000\\
3.25000000000000	13.08776391700000\\
3.50000000000000	14.02358275100000\\
3.75000000000000	14.78777946000000\\
4.00000000000000	15.62171141900000\\
4.25000000000000	16.52545358500000\\
4.50000000000000	17.51310954300000\\
4.75000000000000	18.55970808500000\\
5.00000000000000	19.44150054300000\\
5.25000000000000	20.38769300100000\\
5.50000000000000	21.39966354300000\\
5.75000000000000	22.49268021000000\\
6.00000000000000	23.65595866800000\\
6.25000000000000	24.64827200200000\\
6.50000000000000	25.70847316899999\\
6.75000000000000	26.84107254400000\\
7.00000000000000	28.04934887700000\\
7.25000000000000	29.32347862700000\\
7.50000000000000	30.42634112700000\\
7.75000000000000	31.60591537700000\\
8.00000000000000	32.85423791899999\\
8.25000000000000	34.17332137799999\\
8.50000000000000	35.55936883599999\\
8.75000000000000	36.77990829399999\\
9.00000000000000	38.06552425199999\\
9.25000000000000	39.43013825199999\\
9.50000000000000	40.86275012699999\\
9.75000000000000	42.36154716799999\\
10.00000000000000	43.69108658499999\\
10.25000000000000	45.09090829399999\\
10.50000000000000	46.55861912800000\\
10.75000000000000	48.10013904400000\\
11.00000000000000	49.71393175199999\\
11.25000000000000	51.16605320999999\\
11.50000000000000	52.71261025200000\\
11.75000000000000	54.30175304399999\\
12.00000000000000	55.95947633599999\\
12.25000000000000	57.69469608599999\\
12.50000000000000	59.25498062699999\\
12.75000000000000	60.88035916899999\\
13.00000000000000	62.60784112699999\\
13.25000000000000	64.38071175200000\\
13.50000000000000	66.23038308599999\\
13.75000000000000	67.90440979399999\\
14.00000000000000	69.64422671099999\\
14.25000000000000	71.46428454499998\\
14.50000000000000	73.34817416999998\\
14.75000000000000	75.29672933699997\\
15.00000000000000	77.08110279499998\\
15.25000000000000	78.93337971199998\\
15.50000000000000	80.86001691999998\\
15.75000000000000	82.86618866999999\\
16.00000000000000	84.93073558699999\\
16.25000000000000	88.88532546199998\\
16.50000000000000	91.93044433699998\\
16.75000000000000	93.96840791999998\\
17.00000000000000	96.08187704499997\\
17.25000000000000	98.25741346199997\\
17.50000000000000	100.38924037899997\\
17.75000000000000	102.47915121199998\\
18.00000000000000	104.63323637799998\\
18.25000000000000	106.86083508699998\\
18.50000000000000	109.15124725299998\\
18.75000000000000	111.26939300299998\\
19.00000000000000	113.46626596099998\\
19.25000000000000	115.73710171099998\\
19.50000000000000	118.07734679499998\\
19.75000000000000	120.48595883599998\\
20.00000000000000	122.72934241899998\\
20.25000000000000	125.08332221099998\\
20.50000000000000	127.47151796099998\\
20.75000000000000	129.93055175299997\\
21.00000000000000	132.45604808599998\\
21.25000000000000	134.82151437799999\\
21.50000000000000	137.24080146099999\\
21.75000000000000	139.73147000199998\\
22.00000000000000	142.29671025199997\\
22.25000000000000	144.93419129399996\\
22.50000000000000	147.40743462799995\\
22.75000000000000	149.94815279499994\\
23.00000000000000	152.55787471199994\\
23.25000000000000	155.22823446199993\\
23.50000000000000	157.97759929499992\\
23.75000000000000	160.55534954499993\\
24.00000000000000	163.20788558699994\\
24.25000000000000	166.49312212799992\\
24.50000000000000	169.28036441999993\\
24.75000000000000	172.13923921199992\\
25.00000000000000	174.85164833699992\\
25.25000000000000	177.63001296199991\\
25.50000000000000	180.46998167099991\\
25.75000000000000	183.36696037999991\\
26.00000000000000	186.88543696399992\\
26.25000000000000	189.67983188099993\\
26.50000000000000	192.54573096499993\\
26.75000000000000	195.49325904899993\\
27.00000000000000	198.50852696599992\\
27.25000000000000	201.59013729899991\\
27.50000000000000	204.79803246599991\\
27.75000000000000	207.79504250799991\\
28.00000000000000	210.85332717499992\\
28.25000000000000	213.98596246699992\\
28.50000000000000	217.18204571699991\\
28.75000000000000	220.20677180099992\\
29.00000000000000	223.30556925999991\\
29.25000000000000	226.47184513499991\\
29.50000000000000	229.71243467699992\\
29.75000000000000	233.02168696799993\\
30.00000000000000	236.33863755199994\\
30.25000000000000	239.54999809399993\\
30.50000000000000	242.83436984399992\\
30.75000000000000	246.25375784399992\\
31.00000000000000	249.67525917699993\\
31.25000000000000	252.95750867699994\\
31.50000000000000	256.33582055199992\\
31.75000000000000	259.76165084399992\\
32.00000000000000	263.22714355199992\\
32.25000000000000	266.94972630199993\\
32.50000000000000	270.31456730199994\\
32.75000000000000	273.75873975999997\\
33.00000000000000	277.26681992599998\\
33.25000000000000	280.85275100899997\\
33.50000000000000	284.50240975899999\\
33.75000000000000	287.98912792599998\\
34.00000000000000	291.53844205100000\\
34.25000000000000	295.15809871800002\\
34.50000000000000	299.00400167600003\\
34.75000000000000	302.77479592600002\\
35.00000000000000	306.37492630100002\\
35.25000000000000	310.03856767600001\\
35.50000000000000	313.96265455100001\\
35.75000000000000	317.77234184299999\\
36.00000000000000	321.64337367600001\\
36.25000000000000	325.35481738499999\\
36.50000000000000	329.12638626000000\\
36.75000000000000	332.97244526000003\\
37.00000000000000	336.88926767700002\\
37.25000000000000	341.08480742700004\\
37.50000000000000	344.90618638500007\\
37.75000000000000	348.96012142600006\\
38.00000000000000	352.91687746800005\\
38.25000000000000	356.95074109300003\\
38.50000000000000	361.05480992700001\\
38.75000000000000	364.98984767700000\\
39.00000000000000	369.18537038599999\\
39.25000000000000	373.26035709499996\\
39.50000000000000	377.40920251099999\\
39.75000000000000	381.63694788599997\\
40.00000000000000	385.71011971899998\\
};
\addlegendentry{$N_{\mathrm{remap}} = 5$}

\addplot [color=mycolor2, line width=1.5pt]
  table[row sep=crcr]{%
0.25000000000000	1.25206737500000\\
0.50000000000000	3.38954362500000\\
0.75000000000000	4.25290983400000\\
1.00000000000000	5.19036566700000\\
1.25000000000000	6.13301600000000\\
1.50000000000000	7.26864383300000\\
1.75000000000000	8.22807362500000\\
2.00000000000000	9.17976495900000\\
2.25000000000000	10.20534216700000\\
2.50000000000000	10.79676087500000\\
2.75000000000000	11.44601541700000\\
3.00000000000000	12.16340795900000\\
3.25000000000000	12.94951029300000\\
3.50000000000000	13.81081966800000\\
3.75000000000000	14.73826841800000\\
4.00000000000000	15.73846379300000\\
4.25000000000000	16.80407366800000\\
4.50000000000000	17.94631237600000\\
4.75000000000000	19.15756025100000\\
5.00000000000000	19.87829754300000\\
5.25000000000000	20.63531937600000\\
5.50000000000000	21.46241679300000\\
5.75000000000000	22.36047196000000\\
6.00000000000000	23.32702075100000\\
6.25000000000000	24.36801437600000\\
6.50000000000000	25.47997487600000\\
6.75000000000000	26.65648908400000\\
7.00000000000000	27.91231283400000\\
7.25000000000000	29.23679045900000\\
7.50000000000000	30.03121125000000\\
7.75000000000000	30.90264216600000\\
8.00000000000000	31.85494637500000\\
8.25000000000000	32.87350545800000\\
8.50000000000000	34.43188816600001\\
8.75000000000000	35.58663079100000\\
9.00000000000000	36.81603649900001\\
9.25000000000000	38.11119087400001\\
9.50000000000000	39.47525104100001\\
9.75000000000000	40.90648729100000\\
10.00000000000000	41.82115491600000\\
10.25000000000000	42.80747662500001\\
10.50000000000000	43.86622562500001\\
10.75000000000000	44.99231479100001\\
11.00000000000000	46.18954991600000\\
11.25000000000000	47.45616716600001\\
11.50000000000000	48.79325654100001\\
11.75000000000000	50.20301649900001\\
12.00000000000000	51.68226637400001\\
12.25000000000000	53.23622916600001\\
12.50000000000000	54.26432445800001\\
12.75000000000000	55.36196025000001\\
13.00000000000000	56.52873245800001\\
13.25000000000000	57.76454479200001\\
13.50000000000000	59.07705958400001\\
13.75000000000000	60.45937037600001\\
14.00000000000000	61.91736020900001\\
14.25000000000000	63.44497916700001\\
14.50000000000000	65.03681083300000\\
14.75000000000000	66.70327879200001\\
15.00000000000000	67.84412408400000\\
15.25000000000000	69.05556108400000\\
15.50000000000000	70.34123054200001\\
15.75000000000000	71.70390241700001\\
16.00000000000000	73.12972095800001\\
16.25000000000000	74.62351774900000\\
16.50000000000000	76.18916570800000\\
16.75000000000000	77.82713158300001\\
17.00000000000000	79.53330958300000\\
17.25000000000000	81.31667199900001\\
17.50000000000000	82.57171129100001\\
17.75000000000000	83.90711849900002\\
18.00000000000000	85.30141649900003\\
18.25000000000000	86.77654329100002\\
18.50000000000000	88.31173149900002\\
18.75000000000000	89.92256858200003\\
19.00000000000000	91.60043991600003\\
19.25000000000000	93.34919187400003\\
19.50000000000000	95.17522483200003\\
19.75000000000000	97.06684545700003\\
20.00000000000000	98.43654441600003\\
20.25000000000000	99.87051454100003\\
20.50000000000000	101.37984808300003\\
20.75000000000000	102.96475370800003\\
21.00000000000000	104.61866795800003\\
21.25000000000000	106.33868462400004\\
21.50000000000000	108.12843512400003\\
21.75000000000000	109.99115520700003\\
22.00000000000000	111.92277212300003\\
22.25000000000000	113.93069533100002\\
22.50000000000000	115.41518837300002\\
22.75000000000000	116.96777845600002\\
23.00000000000000	118.58967533100002\\
23.25000000000000	120.28385983100003\\
23.50000000000000	122.04777470600003\\
23.75000000000000	123.87978758100003\\
24.00000000000000	125.77692933100003\\
24.25000000000000	127.76002308100003\\
24.50000000000000	129.79721458100002\\
24.75000000000000	131.90660878900002\\
25.00000000000000	133.50485495500001\\
25.25000000000000	135.18588674700001\\
25.50000000000000	136.91254333000001\\
25.75000000000000	138.70486962200002\\
26.00000000000000	140.58155312200003\\
26.25000000000000	142.53140945600003\\
26.50000000000000	144.55237954000003\\
26.75000000000000	146.64125141500003\\
27.00000000000000	148.79895479000004\\
27.25000000000000	151.19590524800003\\
27.50000000000000	152.90187383100002\\
27.75000000000000	154.69199454000002\\
28.00000000000000	156.54194549800002\\
28.25000000000000	158.63312649800002\\
28.50000000000000	160.62160912300001\\
28.75000000000000	162.67920266500002\\
29.00000000000000	164.81096870700003\\
29.25000000000000	167.01421366500003\\
29.50000000000000	169.28266133200003\\
29.75000000000000	171.62536279000003\\
30.00000000000000	173.44486945700004\\
30.25000000000000	175.33351158200003\\
30.50000000000000	177.29793333200004\\
30.75000000000000	179.33502095700004\\
31.00000000000000	181.44274429000004\\
31.25000000000000	183.65237912300003\\
31.50000000000000	185.95018916400002\\
31.75000000000000	188.47292053900003\\
32.00000000000000	190.85799420600003\\
32.25000000000000	193.31464391400002\\
32.50000000000000	195.24282866400003\\
32.75000000000000	197.24262445500003\\
33.00000000000000	199.31996174700004\\
33.25000000000000	201.46979987200004\\
33.50000000000000	203.68678791400004\\
33.75000000000000	205.97408095600002\\
34.00000000000000	208.33372974800002\\
34.25000000000000	210.75955562300001\\
34.50000000000000	213.25703795600000\\
34.75000000000000	215.82974062299999\\
35.00000000000000	217.88184712300000\\
35.25000000000000	219.99898103999999\\
35.50000000000000	222.18839678999998\\
35.75000000000000	224.44754295700000\\
36.00000000000000	226.78184704000000\\
36.25000000000000	229.18062920700001\\
36.50000000000000	231.95040354000000\\
36.75000000000000	234.49018037400000\\
37.00000000000000	237.10596679100001\\
37.25000000000000	239.78688954100002\\
37.50000000000000	242.11076195800001\\
37.75000000000000	244.34196679100000\\
38.00000000000000	246.85061779099999\\
38.25000000000000	249.21561383299999\\
38.50000000000000	252.09360699900000\\
38.75000000000000	254.60768033200000\\
39.00000000000000	257.19222208200000\\
39.25000000000000	259.98290612300002\\
39.50000000000000	262.70783195600001\\
39.75000000000000	265.50307378899998\\
40.00000000000000	267.78942245600001\\
};
\addlegendentry{$N_{\mathrm{remap}} = 10$}

\addplot [color=mycolor3, line width=1.5pt]
  table[row sep=crcr]{%
0.25000000000000	1.13973766700000\\
0.50000000000000	3.28072641700000\\
0.75000000000000	4.16621408300000\\
1.00000000000000	5.10561166700000\\
1.25000000000000	6.05345587500000\\
1.50000000000000	7.20026291700000\\
1.75000000000000	8.15129775100000\\
2.00000000000000	9.10196062600000\\
2.25000000000000	10.19540962600000\\
2.50000000000000	11.29440250100000\\
2.75000000000000	12.46433108400000\\
3.00000000000000	13.70368925100000\\
3.25000000000000	15.01267650100000\\
3.50000000000000	16.39125445900000\\
3.75000000000000	17.84696187500000\\
4.00000000000000	19.36901391700000\\
4.25000000000000	20.95953275000000\\
4.50000000000000	22.61996816700000\\
4.75000000000000	24.35365229200000\\
5.00000000000000	24.94938462600000\\
5.25000000000000	25.60233295900000\\
5.50000000000000	26.31769804200000\\
5.75000000000000	27.10332875100000\\
6.00000000000000	27.96266796000000\\
6.25000000000000	28.88856787600000\\
6.50000000000000	29.88358250100000\\
6.75000000000000	30.94988962600000\\
7.00000000000000	32.08749375100000\\
7.25000000000000	33.29754583400000\\
7.50000000000000	34.58025304200000\\
7.75000000000000	35.93717329200000\\
8.00000000000000	37.35936904200000\\
8.25000000000000	38.85474670900000\\
8.50000000000000	40.41503070900001\\
8.75000000000000	42.06504654200000\\
9.00000000000000	43.77309429200000\\
9.25000000000000	45.57360895800000\\
9.50000000000000	47.42075916700000\\
9.75000000000000	49.33889925000000\\
10.00000000000000	50.03170995800000\\
10.25000000000000	50.78748841600000\\
10.50000000000000	51.61192800000000\\
10.75000000000000	52.50730533300000\\
11.00000000000000	53.47522600000000\\
11.25000000000000	54.51162783300000\\
11.50000000000000	55.61839666600000\\
11.75000000000000	56.80312670799999\\
12.00000000000000	58.05724283300000\\
12.25000000000000	59.39088504100000\\
12.50000000000000	60.79031333300000\\
12.75000000000000	62.25558320800000\\
13.00000000000000	63.78844533300000\\
13.25000000000000	65.39647883300000\\
13.50000000000000	67.07147650000000\\
13.75000000000000	68.82024270800001\\
14.00000000000000	70.63816845800001\\
14.25000000000000	72.54917687500001\\
14.50000000000000	74.51723687500001\\
14.75000000000000	76.55417225000001\\
15.00000000000000	77.35769195900001\\
15.25000000000000	78.22718987600001\\
15.50000000000000	79.17108650100000\\
15.75000000000000	80.18792675100001\\
16.00000000000000	81.27531975100001\\
16.25000000000000	82.42587812600001\\
16.50000000000000	83.64788850100001\\
16.75000000000000	84.94524625100001\\
17.00000000000000	86.31177212600001\\
17.25000000000000	87.75581045900000\\
17.50000000000000	89.26504145900000\\
17.75000000000000	90.84954158400001\\
18.00000000000000	92.49754300100001\\
18.25000000000000	94.21792983400002\\
18.50000000000000	96.00676229200002\\
18.75000000000000	97.86719579200002\\
19.00000000000000	99.80706908400002\\
19.25000000000000	101.81785908400002\\
19.50000000000000	103.90156012500002\\
19.75000000000000	106.04722895800002\\
20.00000000000000	106.96257745800001\\
20.25000000000000	107.94774020800001\\
20.50000000000000	109.00205950000002\\
20.75000000000000	110.13763033300002\\
21.00000000000000	111.33841791600001\\
21.25000000000000	112.61881691600001\\
21.50000000000000	113.95639545700001\\
21.75000000000000	115.36946695700001\\
22.00000000000000	116.84854358200002\\
22.25000000000000	118.40229412300002\\
22.50000000000000	120.02727545600001\\
22.75000000000000	121.72219595600001\\
23.00000000000000	123.48354641500001\\
23.25000000000000	125.31295654000000\\
23.50000000000000	127.21959599900001\\
23.75000000000000	129.20023408200001\\
24.00000000000000	131.25112516600001\\
24.25000000000000	133.37644999900002\\
24.50000000000000	135.57409324900001\\
24.75000000000000	137.85069979100001\\
25.00000000000000	138.90282191600002\\
25.25000000000000	140.00804591600001\\
25.50000000000000	141.17214512400002\\
25.75000000000000	142.41216695700001\\
26.00000000000000	143.71814337400002\\
26.25000000000000	145.09921574900002\\
26.50000000000000	146.54804149900002\\
26.75000000000000	148.07672674900002\\
27.00000000000000	149.66827179100002\\
27.25000000000000	151.32835975000003\\
27.50000000000000	153.06349762500002\\
27.75000000000000	154.87898258300004\\
28.00000000000000	156.75364237500003\\
28.25000000000000	158.70242958400001\\
28.50000000000000	160.72035270900003\\
28.75000000000000	162.81691908400003\\
29.00000000000000	164.98713045900004\\
29.25000000000000	167.22521654300004\\
29.50000000000000	169.52560770900004\\
29.75000000000000	171.90297679200003\\
30.00000000000000	173.04472370800002\\
30.25000000000000	174.25444720800002\\
30.50000000000000	175.53833062500001\\
30.75000000000000	176.89349275000001\\
31.00000000000000	178.31532191600002\\
31.25000000000000	179.83178158300001\\
31.50000000000000	181.45109791700000\\
31.75000000000000	183.12395508399999\\
32.00000000000000	184.84907591699999\\
32.25000000000000	186.62954483400000\\
32.50000000000000	188.46798941800000\\
32.75000000000000	190.38174583500000\\
33.00000000000000	192.36221700199999\\
33.25000000000000	194.41950679400000\\
33.50000000000000	196.54872287800001\\
33.75000000000000	198.74773771100001\\
34.00000000000000	201.01677308600000\\
34.25000000000000	203.36334366899999\\
34.50000000000000	205.77557070999998\\
34.75000000000000	208.26260137599999\\
35.00000000000000	209.50867241800000\\
35.25000000000000	210.83107496000000\\
35.50000000000000	212.22722700099999\\
35.75000000000000	213.69519962599998\\
36.00000000000000	215.23480358399996\\
36.25000000000000	216.84770895899996\\
36.50000000000000	218.52816550099996\\
36.75000000000000	220.27262566699997\\
37.00000000000000	222.09254949999996\\
37.25000000000000	223.98292683299996\\
37.50000000000000	225.94198620799997\\
37.75000000000000	227.97247020799998\\
38.00000000000000	230.07516670799998\\
38.25000000000000	232.24944404199999\\
38.50000000000000	234.49288429199999\\
38.75000000000000	236.81687837499999\\
39.00000000000000	239.19958229199997\\
39.25000000000000	241.67345929199996\\
39.50000000000000	244.20402195899996\\
39.75000000000000	246.81078691699997\\
40.00000000000000	248.19146641699999\\
};
\addlegendentry{$N_{\mathrm{remap}} = 20$}

\addplot [color=mycolor4, line width=1.5pt]
  table[row sep=crcr]{%
0.25000000000000	1.14315841700000\\
0.50000000000000	3.29127587500000\\
0.75000000000000	4.15745300000000\\
1.00000000000000	5.09275462500000\\
1.25000000000000	6.03627025000000\\
1.50000000000000	7.17978575000000\\
1.75000000000000	8.13384379200000\\
2.00000000000000	9.08864045900000\\
2.25000000000000	10.11401262600000\\
2.50000000000000	11.21288083400000\\
2.75000000000000	12.38203100100000\\
3.00000000000000	13.61714658500000\\
3.25000000000000	14.92268833500000\\
3.50000000000000	16.30726266800000\\
3.75000000000000	17.76137337600000\\
4.00000000000000	19.28004533500000\\
4.25000000000000	20.87148821000000\\
4.50000000000000	22.53164358500000\\
4.75000000000000	24.26520612600000\\
5.00000000000000	26.07532220900000\\
5.25000000000000	27.95737433400000\\
5.50000000000000	29.90327237500000\\
5.75000000000000	31.92456283400000\\
6.00000000000000	34.01356108400000\\
6.25000000000000	36.17911083400000\\
6.50000000000000	38.41433712600000\\
6.75000000000000	40.72184941800000\\
7.00000000000000	43.08903354300000\\
7.25000000000000	45.53267191800000\\
7.50000000000000	46.31117937600000\\
7.75000000000000	47.41278400100000\\
8.00000000000000	48.81119437600000\\
8.25000000000000	49.62111258400000\\
8.50000000000000	50.48437212600000\\
8.75000000000000	51.41454387600000\\
9.00000000000000	53.16972866800000\\
9.25000000000000	54.29006721000000\\
9.50000000000000	55.44029091800000\\
9.75000000000000	56.65512062600000\\
10.00000000000000	57.94450741700000\\
10.25000000000000	59.30538512600000\\
10.50000000000000	60.72995720900000\\
10.75000000000000	62.21924762500000\\
11.00000000000000	63.78376937500000\\
11.25000000000000	65.42019491700000\\
11.50000000000000	67.13312041700000\\
11.75000000000000	68.90697487600001\\
12.00000000000000	70.76316787600001\\
12.25000000000000	72.67878050100001\\
12.50000000000000	74.66901929300001\\
12.75000000000000	76.72723912600000\\
13.00000000000000	78.85827429299999\\
13.25000000000000	81.06381937600000\\
13.50000000000000	83.33623016700000\\
13.75000000000000	85.68979258300000\\
14.00000000000000	88.11617658300000\\
14.25000000000000	90.60292491599999\\
14.50000000000000	93.16021920700000\\
14.75000000000000	95.80624945699999\\
15.00000000000000	96.50056653999999\\
15.25000000000000	97.26229133100000\\
15.50000000000000	98.09544508099999\\
15.75000000000000	99.00973887200000\\
16.00000000000000	99.98376595500000\\
16.25000000000000	101.01919337199999\\
16.50000000000000	102.13856870499998\\
16.75000000000000	103.32318303799998\\
17.00000000000000	104.57594357899998\\
17.25000000000000	105.89801687099998\\
17.50000000000000	107.29100062099998\\
17.75000000000000	108.76628653799997\\
18.00000000000000	110.30392820499998\\
18.25000000000000	111.91874041299998\\
18.50000000000000	113.60638545499998\\
18.75000000000000	115.35820345499998\\
19.00000000000000	117.19677624599998\\
19.25000000000000	119.09202978799999\\
19.50000000000000	121.05827041299999\\
19.75000000000000	123.09399103799998\\
20.00000000000000	125.19695645499999\\
20.25000000000000	127.36943691399999\\
20.50000000000000	129.61749157999998\\
20.75000000000000	131.94043795499996\\
21.00000000000000	134.32705770499996\\
21.25000000000000	136.78304253799996\\
21.50000000000000	139.31165137099995\\
21.75000000000000	141.90841145499996\\
22.00000000000000	144.57741491299996\\
22.25000000000000	147.32757245499997\\
22.50000000000000	148.12774137199997\\
22.75000000000000	148.99866108099997\\
23.00000000000000	149.94653591399995\\
23.25000000000000	150.96440987199995\\
23.50000000000000	152.05098124699995\\
23.75000000000000	153.20646712199996\\
24.00000000000000	154.43191899699997\\
24.25000000000000	155.74154920499996\\
24.50000000000000	157.10869224699996\\
24.75000000000000	158.54732912199995\\
25.00000000000000	160.07137949699995\\
25.25000000000000	161.67494762199996\\
25.50000000000000	163.32981912199998\\
25.75000000000000	165.04968087199998\\
26.00000000000000	166.84170191299998\\
26.25000000000000	168.70295999699999\\
26.50000000000000	170.63614941300000\\
26.75000000000000	172.64927099599998\\
27.00000000000000	174.72953245399998\\
27.25000000000000	176.87971711999998\\
27.50000000000000	179.09431761999997\\
27.75000000000000	181.38469382799997\\
28.00000000000000	183.74099574399997\\
28.25000000000000	186.16932682799998\\
28.50000000000000	188.68212541099999\\
28.75000000000000	191.26338528599999\\
29.00000000000000	193.91963811899998\\
29.25000000000000	196.63271199399998\\
29.50000000000000	199.41997945199998\\
29.75000000000000	202.27222070199997\\
30.00000000000000	203.18965695199998\\
30.25000000000000	204.17549724399998\\
30.50000000000000	205.23676953499998\\
30.75000000000000	206.36277390999999\\
31.00000000000000	207.56038707599998\\
31.25000000000000	208.85519140899999\\
31.50000000000000	210.23168286799998\\
31.75000000000000	211.65406661799997\\
32.00000000000000	213.13517932699997\\
32.25000000000000	214.68666886899999\\
32.50000000000000	216.31219736899999\\
32.75000000000000	218.00321670200000\\
33.00000000000000	219.76680274400002\\
33.25000000000000	221.60480078500001\\
33.50000000000000	223.50692774400000\\
33.75000000000000	225.48394078600001\\
34.00000000000000	227.52988103600001\\
34.25000000000000	229.64795541100000\\
34.50000000000000	231.83723803600000\\
34.75000000000000	234.10681753599999\\
35.00000000000000	236.43690957699999\\
35.25000000000000	238.83961632699999\\
35.50000000000000	241.30837799399998\\
35.75000000000000	243.85534836899998\\
36.00000000000000	246.48286803500000\\
36.25000000000000	249.16733470200001\\
36.50000000000000	251.92180203500001\\
36.75000000000000	254.74605199300001\\
37.00000000000000	257.64175132600002\\
37.25000000000000	260.61496395099999\\
37.50000000000000	261.64491032600000\\
37.75000000000000	262.73326436799999\\
38.00000000000000	263.89031561799999\\
38.25000000000000	265.12323603499999\\
38.50000000000000	266.42418532599999\\
38.75000000000000	267.79463624199997\\
39.00000000000000	269.23904365799996\\
39.25000000000000	270.75338565799996\\
39.50000000000000	272.34558811699998\\
39.75000000000000	274.00711440900000\\
40.00000000000000	275.74430915900001\\
};
\addlegendentry{$N_{\mathrm{remap}} = 30$}

\addplot [color=mycolor5, line width=1.5pt]
  table[row sep=crcr]{%
0.25000000000000	1.15166720800000\\
0.50000000000000	3.29581337500000\\
0.75000000000000	4.16285466600000\\
1.00000000000000	5.09799850000000\\
1.25000000000000	6.04150533400000\\
1.50000000000000	7.19059175000000\\
1.75000000000000	8.14222333400000\\
2.00000000000000	9.09676945900000\\
2.25000000000000	10.12322541700000\\
2.50000000000000	11.22447458400000\\
2.75000000000000	12.39653266700000\\
3.00000000000000	13.63842687500000\\
3.25000000000000	14.94710441700000\\
3.50000000000000	16.33244508400000\\
3.75000000000000	17.78161187600000\\
4.00000000000000	19.30975708500000\\
4.25000000000000	20.90461362600000\\
4.50000000000000	22.56806566800000\\
4.75000000000000	24.30827233500000\\
5.00000000000000	26.11594529300000\\
5.25000000000000	27.99603862600000\\
5.50000000000000	29.94599325100000\\
5.75000000000000	31.96375708500000\\
6.00000000000000	34.05634487700000\\
6.25000000000000	36.22807796000000\\
6.50000000000000	38.45797266800000\\
6.75000000000000	40.76201966799999\\
7.00000000000000	43.13837287599999\\
7.25000000000000	45.58533112599999\\
7.50000000000000	48.10646045899999\\
7.75000000000000	50.69865979199999\\
8.00000000000000	53.35864774999999\\
8.25000000000000	56.10019487499999\\
8.50000000000000	58.90267270799998\\
8.75000000000000	61.77117908299999\\
9.00000000000000	64.71222804199999\\
9.25000000000000	67.72304620800000\\
9.50000000000000	70.80581874900000\\
9.75000000000000	73.96992862399999\\
10.00000000000000	74.57093366599999\\
10.25000000000000	75.22621095799998\\
10.50000000000000	75.94695058299997\\
10.75000000000000	76.73244504099998\\
11.00000000000000	77.58979099999998\\
11.25000000000000	78.51679629199998\\
11.50000000000000	79.51810604199999\\
11.75000000000000	80.58765720799998\\
12.00000000000000	81.72954458299998\\
12.25000000000000	82.93814424999998\\
12.50000000000000	84.22588299999998\\
12.75000000000000	85.58266049999999\\
13.00000000000000	87.00609891699999\\
13.25000000000000	88.49746924999999\\
13.50000000000000	90.06216149999999\\
13.75000000000000	91.69914599999998\\
14.00000000000000	93.41704799999998\\
14.25000000000000	95.19605324999998\\
14.50000000000000	97.04538783299998\\
14.75000000000000	98.96961766599998\\
15.00000000000000	100.97026137399997\\
15.25000000000000	103.03083066599997\\
15.50000000000000	105.16951979099997\\
15.75000000000000	107.39111649899996\\
16.00000000000000	109.67484683299996\\
16.25000000000000	112.02508699999996\\
16.50000000000000	114.44198808399996\\
16.75000000000000	116.92836220899996\\
17.00000000000000	119.48759574999995\\
17.25000000000000	122.11916749999996\\
17.50000000000000	124.81831845799996\\
17.75000000000000	127.59873137399997\\
18.00000000000000	130.44242741499997\\
18.25000000000000	133.35836241499996\\
18.50000000000000	136.34625387299997\\
18.75000000000000	139.39998112299998\\
19.00000000000000	142.52925808099997\\
19.25000000000000	145.72626545599996\\
19.50000000000000	149.00177653899996\\
19.75000000000000	152.34417299699996\\
20.00000000000000	153.03431499699997\\
20.25000000000000	153.78972166399996\\
20.50000000000000	154.61906295599996\\
20.75000000000000	155.52726687299997\\
21.00000000000000	156.50080637299996\\
21.25000000000000	157.54702478999997\\
21.50000000000000	158.65966541499998\\
21.75000000000000	159.84183187399998\\
22.00000000000000	161.09870499899998\\
22.25000000000000	162.42072520699998\\
22.50000000000000	163.82214433199997\\
22.75000000000000	165.29323458199997\\
23.00000000000000	166.83422362299996\\
23.25000000000000	168.44521233199995\\
23.50000000000000	170.12516216499995\\
23.75000000000000	171.87749653999995\\
24.00000000000000	173.69517424799994\\
24.25000000000000	175.58763762299995\\
24.50000000000000	177.68537591499995\\
24.75000000000000	179.72469524799996\\
25.00000000000000	181.85432145599995\\
25.25000000000000	184.05671199799994\\
25.50000000000000	186.30782095599994\\
25.75000000000000	188.62275266399993\\
26.00000000000000	191.01806816399994\\
26.25000000000000	193.47898416399994\\
26.50000000000000	196.01066103899993\\
26.75000000000000	198.60981412199993\\
27.00000000000000	201.28157803899992\\
27.25000000000000	204.02315666399991\\
27.50000000000000	206.83336932999993\\
27.75000000000000	209.72550087199991\\
28.00000000000000	212.68619899699991\\
28.25000000000000	215.72754312199990\\
28.50000000000000	218.83574295599990\\
28.75000000000000	222.00455608099992\\
29.00000000000000	225.24444595599991\\
29.25000000000000	228.55254333099990\\
29.50000000000000	231.93578508099989\\
29.75000000000000	235.38776220599991\\
30.00000000000000	236.18871912199990\\
30.25000000000000	237.05977574699989\\
30.50000000000000	238.00255474699989\\
30.75000000000000	239.02186132999989\\
31.00000000000000	240.11003557999990\\
31.25000000000000	241.29083112199990\\
31.50000000000000	242.55242291299990\\
31.75000000000000	243.85636762099989\\
32.00000000000000	245.22037699599989\\
32.25000000000000	246.65845253699987\\
32.50000000000000	248.16688741199988\\
32.75000000000000	249.74751549599989\\
33.00000000000000	251.39692082999989\\
33.25000000000000	253.12321845499989\\
33.50000000000000	254.92760141399989\\
33.75000000000000	256.79222437299990\\
34.00000000000000	258.72506133199988\\
34.25000000000000	260.73060245699986\\
34.50000000000000	262.80582312399986\\
34.75000000000000	264.94621312399988\\
35.00000000000000	267.16426758299986\\
35.25000000000000	269.45595299899986\\
35.50000000000000	271.81646154099985\\
35.75000000000000	274.24866037399983\\
36.00000000000000	276.75786870699983\\
36.25000000000000	279.32957016499984\\
36.50000000000000	281.97456562399981\\
36.75000000000000	284.68944562399980\\
37.00000000000000	287.47724595799980\\
37.25000000000000	290.33403720799981\\
37.50000000000000	293.26597566599980\\
37.75000000000000	296.26627262399978\\
38.00000000000000	299.34433049899980\\
38.25000000000000	302.48937128999978\\
38.50000000000000	305.74730858099980\\
38.75000000000000	309.03262124799983\\
39.00000000000000	312.39892916499980\\
39.25000000000000	315.81952212299979\\
39.50000000000000	319.31174853999977\\
39.75000000000000	322.87818578999975\\
40.00000000000000	323.79140249799974\\
};
\addlegendentry{$N_{\mathrm{remap}} = 40$}

\addplot [color=mycolor6, line width=1.5pt]
  table[row sep=crcr]{%
0.25000000000000	1.15651558300000\\
0.50000000000000	3.30231245800000\\
0.75000000000000	4.16422645800000\\
1.00000000000000	5.09783566600000\\
1.25000000000000	6.03625312400000\\
1.50000000000000	7.19831608300000\\
1.75000000000000	8.15469779200000\\
2.00000000000000	9.11277950000000\\
2.25000000000000	10.13460679200000\\
2.50000000000000	11.22881795900000\\
2.75000000000000	12.39195675100000\\
3.00000000000000	13.62634600100000\\
3.25000000000000	14.93693641800000\\
3.50000000000000	16.30782250100000\\
3.75000000000000	17.75855533400000\\
4.00000000000000	19.27261104200000\\
4.25000000000000	20.85803100100000\\
4.50000000000000	22.51945341800000\\
4.75000000000000	24.25433937600000\\
5.00000000000000	26.05122400100000\\
5.25000000000000	27.92977025100000\\
5.50000000000000	29.87786996000000\\
5.75000000000000	31.90082125200000\\
6.00000000000000	33.98754137700000\\
6.25000000000000	36.14682216900000\\
6.50000000000000	38.37436421100000\\
6.75000000000000	40.67404967000000\\
7.00000000000000	43.04716429500000\\
7.25000000000000	45.50141837800000\\
7.50000000000000	48.01634650300000\\
7.75000000000000	50.59891575300000\\
8.00000000000000	53.25327117000000\\
8.25000000000000	55.98346525300000\\
8.50000000000000	58.78469879400000\\
8.75000000000000	61.65060216900000\\
9.00000000000000	64.59084033600000\\
9.25000000000000	67.60842175300000\\
9.50000000000000	70.69229729500000\\
9.75000000000000	73.84717050400000\\
10.00000000000000	77.06695692100000\\
10.25000000000000	80.36236288000001\\
10.50000000000000	83.72994438000001\\
10.75000000000000	87.16315929700001\\
11.00000000000000	90.66945904700000\\
11.25000000000000	94.24236692200000\\
11.50000000000000	97.90096513000000\\
11.75000000000000	101.61289592200001\\
12.00000000000000	105.40390613100001\\
12.25000000000000	109.27479271400001\\
12.50000000000000	113.20807529800001\\
12.75000000000000	117.20746046400001\\
13.00000000000000	121.27389163000001\\
13.25000000000000	125.41399533800001\\
13.50000000000000	129.62475371300002\\
13.75000000000000	133.90553842200001\\
14.00000000000000	138.26345533900002\\
14.25000000000000	142.69722658900002\\
14.50000000000000	147.19381254700002\\
14.75000000000000	151.76265067200004\\
15.00000000000000	156.40408354700003\\
15.25000000000000	161.11018292200004\\
15.50000000000000	165.90490092200002\\
15.75000000000000	170.77256150500003\\
16.00000000000000	175.69546533900004\\
16.25000000000000	180.69078625500003\\
16.50000000000000	185.75301100500005\\
16.75000000000000	190.88307304600005\\
17.00000000000000	196.08824329600006\\
17.25000000000000	201.36804363000005\\
17.50000000000000	206.70834388000006\\
17.75000000000000	212.12932838000006\\
18.00000000000000	218.11064083900007\\
18.25000000000000	223.68113417200007\\
18.50000000000000	229.32069204700008\\
18.75000000000000	235.02601058800008\\
19.00000000000000	240.79824417100008\\
19.25000000000000	246.63763929600009\\
19.50000000000000	252.55036979600010\\
19.75000000000000	258.53365513000011\\
20.00000000000000	264.56209242200009\\
20.25000000000000	270.68844238000008\\
20.50000000000000	278.34789496300010\\
20.75000000000000	284.62068546300009\\
21.00000000000000	290.95945225500009\\
21.25000000000000	297.36974417100009\\
21.50000000000000	303.84515721300011\\
21.75000000000000	310.40354942100009\\
22.00000000000000	317.01702525400009\\
22.25000000000000	323.70436842100008\\
22.50000000000000	330.47751467100011\\
22.75000000000000	337.30775879600009\\
23.00000000000000	344.21314763000009\\
23.25000000000000	351.19583133900011\\
23.50000000000000	358.24201829700013\\
23.75000000000000	365.35476154700012\\
24.00000000000000	372.54564029700015\\
24.25000000000000	379.79728746300015\\
24.50000000000000	387.12002571300013\\
24.75000000000000	394.52535542100014\\
25.00000000000000	402.00745917100016\\
25.25000000000000	409.59016333800014\\
25.50000000000000	417.38954117100013\\
25.75000000000000	425.33785954600012\\
26.00000000000000	433.11018092100011\\
26.25000000000000	440.93017354600011\\
26.50000000000000	448.81953467100010\\
26.75000000000000	456.81150054600010\\
27.00000000000000	464.84565921300009\\
27.25000000000000	472.95682508800007\\
27.50000000000000	481.13858592100007\\
27.75000000000000	489.38737354600005\\
28.00000000000000	497.71263800500003\\
28.25000000000000	506.12053142200006\\
28.50000000000000	514.58129508900004\\
28.75000000000000	523.10981271399999\\
29.00000000000000	531.70396850600002\\
29.25000000000000	540.40310721399999\\
29.50000000000000	549.16000396400000\\
29.75000000000000	557.98664300600001\\
30.00000000000000	566.87860788099999\\
30.25000000000000	575.82686000599995\\
30.50000000000000	584.84935413099993\\
30.75000000000000	593.95777117299997\\
31.00000000000000	603.13007146500001\\
31.25000000000000	612.40629142399996\\
31.50000000000000	621.76464279899994\\
31.75000000000000	631.16625150699997\\
32.00000000000000	640.61244592399999\\
32.25000000000000	650.17876154900000\\
32.50000000000000	659.76615321600002\\
32.75000000000000	669.47743350799999\\
33.00000000000000	679.21122050800000\\
33.25000000000000	689.00781746700000\\
33.50000000000000	698.87855992499999\\
33.75000000000000	708.82825404999994\\
34.00000000000000	718.82614930000000\\
34.25000000000000	728.93085692499994\\
34.50000000000000	739.10093829999994\\
34.75000000000000	749.33853075799993\\
35.00000000000000	759.63354942499996\\
35.25000000000000	769.99887996699999\\
35.50000000000000	780.44159821699998\\
35.75000000000000	790.96227230099998\\
36.00000000000000	801.53991446700002\\
36.25000000000000	812.15898484200000\\
36.50000000000000	822.85555855099994\\
36.75000000000000	833.63988884299999\\
37.00000000000000	844.51338467599999\\
37.25000000000000	855.43446496699994\\
37.50000000000000	867.61644888399996\\
37.75000000000000	878.71449471799997\\
38.00000000000000	889.86276500999998\\
38.25000000000000	901.08699659299998\\
38.50000000000000	912.38366788500002\\
38.75000000000000	923.73432071800005\\
39.00000000000000	935.15454271800002\\
39.25000000000000	946.63117067600001\\
39.50000000000000	958.20844496699999\\
39.75000000000000	969.84651000899999\\
40.00000000000000	981.54673446799995\\
};
\addlegendentry{NuFi}

\end{axis}
\end{tikzpicture}%